\documentclass[twoside,leqno,fleqn]{article}
\usepackage{amsmath}
\usepackage{amssymb}
\usepackage{array}
\usepackage[shortlabels]{enumitem}
\usepackage[onehalfspacing]{setspace}
\usepackage[T1]{fontenc}
\usepackage{graphicx}
\usepackage{xcolor}
\usepackage[lowtilde]{url}
\usepackage{hyperref}
\usepackage{orcidlink}

\newcommand{\bbbn}{\mathbb{N}}

\newcommand{\qed}{\hfill{$\Box$}} 
\newtheorem {overview}{Overview} 
\newtheorem{satz}{Satz}[section]

\newtheorem{proposition}[satz]{Proposition}
\newtheorem{corollary}[satz]{Corollary}
\newtheorem{theorem}[satz]{Theorem}
\newtheorem{remark}[satz]{Remark}
\newtheorem{example}[satz]{Example}

\newtheorem{agree}[satz]{Agreement}

\newtheorem{consequ}[satz]{Consequences}
\counterwithout{equation}{section}

\begin{document}
\pagestyle{headings} 
\thispagestyle{empty}

\title{Abstract computation over first-order structures. Part II\,b: Moschovakis' operator and other non-determinisms} \author{Christine Ga\ss ner}

\begin{center}

{\Large Abstract Computation over First-Order Structures \vspace{0.2cm}}\\{\Large Part II\,b\\ Moschovakis' Operator and Other Non-Determinisms}
\vspace{0.5cm}\\
{\bf Christine Ga\ss ner} \orcidlink{0009-0005-5163-4681}
 \footnote{In Part III and Part IV, we present a generalization of the results discussed at the CCC in Kochel in 2015. My thanks go to Dieter Spreen, Ulrich Berger, and the other organizers of the conference in Kochel. Further results were discussed at the CCA 2015. I would like to thank the organizers of this conference and in particular Martin Ziegler and Akitoshi Kawamura.} {\vspace{0.1cm}\\University of Greifswald, Germany, 2025\\ gassnerc@uni-greifswald.de}\\\end{center} 

\vspace{0.2cm}

\begin{abstract} 
BSS RAMs were introduced to provide a mathematical framework for characterizing algorithms over first-order structures. Non-deterministic BSS RAMs help to model different non-deterministic approaches. Here, we deal with different types of binary non-determinisms and study the consequences of the decidability of the identity relation and the decidability of finite sets consisting of one or two constants. We compare the binary non-determinism resulting from a non-deterministic branching process, the digital non-determinism resulting from the restriction of guesses to two constants, and some other non-determinisms resulting from the use of Moschovakis' operator applied to oracle sets restricted to tuples of constants. Moreover, we show that the performance capability and the efficiency of individual machines are influenced by the following properties. 1. The identity relation belongs to the underlying structure. 2. The identity is semi-decidable over the underlying structure. 3. Two single-element sets of constants are semi-decidable. 4. A set of two constants is semi-decidable. The order of these properties corresponds to the strength of their influence. In all cases mentioned, the semi-decidability of the sets implies their decidability. 
\end{abstract}

\vspace{0.5cm}


\vspace{1cm}
\newpage

\tableofcontents

\vspace{0.3cm}

\noindent {\large\bf Contents of Part I} 

{\small
\begin{description}[itemindent=-1cm]\parskip -0.5mm

\item {\bf \hspace{0.1cm}1\quad Introduction} 

\item {\bf \hspace{0.1cm}2\quad BSS-RAMs over first-order structures.} 
First-order structures and their signatures. The $\sigma$-programs. Finite machines and more. Infinite machines and their configurations. Computation and semi-decidability by BSS RAMs. Non-deterministic BSS RAMs. Oracle machines. Other non-deterministic BSS RAMs.

\item {\bf \hspace{0.1cm}3\quad Examples.} 
Input and output functions. Computing functions and deciding sets. Evaluating formulas. 
\end{description}

\vspace{0.1cm}

\noindent {\large\bf Contents of Part II\,a} 

{\small
\begin{description}[itemindent=-1cm]\parskip -0.5mm
\item {\bf \hspace{0.1cm}4\quad BSS RAMs with Moschovakis’ operator.}
Moschovakis’ non-deterministic operator $\nu$. Examples. The power of non-deterministic BSS RAMs. Evaluation of $\nu$-instructions for $n$-ary or decidable $Q$’s. 

\item {\bf \hspace{0.1cm}5\quad Simulation of $\nu$-oracle machines.}
Useful subprograms and more. The non-deterministic simulation of $\nu$-oracle machines. Replacing $\nu$-instructions by trying to semi-decide $Q$. More details: The simulation of ${\cal N}_Q$ semi-deciding $Q$.
\end{description}
\noindent For Part I and Part II\,a (including introductions, repetitions, summaries, and outlooks), see \cite{GASS25A} and \cite{GASS25B}.}
}

\newpage
\section*{Introduction and repetition} 
\addcontentsline{toc}{section}{\bf Introduction and repetition}
\markboth{INTRODUCTION AND REPETITION}{INTRODUCTION AND REPETITION}
Different types of abstract machines over any first-order structure can lead to different classes of semi-decidable problems. The decidability and complexity of such problems and their level or degree and their completeness in certain hierarchies of problems also result from the properties of the underlying structures. In particular, the number of specific forms of non-determinisms that can lead to different classes of result functions or halting sets depends on whether the identity relation is available to perform equality tests or how many constants belong to a structure. The following characterization takes into account several types of structures and also includes an analysis of {\sf computation paths determined by the sequences of configurations}. Let us recall and introduce some notions. 

\vspace{0.3cm}

{\bf $\sigma$-Programs.} For any ${\cal A}$-machine ${\cal M}$, its {\sf program} ${\cal P}_{\cal M}$ has the form

\vspace{0.05cm}

$1: {\sf instruction}_{1}; \quad \ldots;\quad\ell _{{\cal P}_{\cal M}}-1:\, {\sf instruction}_{\ell _{{\cal P}_{\cal M}}-1};\quad \,\,\,\ell _{{\cal P}_{\cal M}}:\, {\sf stop}.\hfill \textcolor{blue}{(*)}\,\,$

\vspace{0.05cm}

\noindent
Any substring \,
$i: {\sf instruction}_{i}; \quad \ldots\,;\quad j: {\sf instruction}_j;\,\,\,\,\big[{j+1}: {\sf stop.}\big]$ \, with $1\leq i\leq j<\ell_{{\cal P}_{\cal M}}$
is a {\sf program segment} of ${\cal P}_{\cal M}$. {\sf subpr} can be a {\sf subprogram} of ${\cal P}_{\cal M}$ if \, $\ell: {\sf subpr};$ \, is a program segment of ${\cal P}_{\cal M}$ (which means that any label of an instruction in {\sf subpr} is not the label of an instruction in another part of ${\cal P}_{\cal M}$). 

\vspace{0.1cm}

{\bf The classes ${\rm SDEC}_{\cal A}$ and ${\rm DEC}_{\cal A}$.} For any ${\cal A}$-machine ${\cal M }$, the {\em halting set $H_{\cal M}$ of ${\cal M }$} is given by $H_{\cal M}=\{ \vec x \in U_{\cal A}^{\infty} \mid {\cal M} ( \vec x)\downarrow \}$. $\{ H_{\cal M}\mid {\cal M}\in {\sf M}_{\cal A} \}$ is the class ${\rm SDEC}_{\cal A}$ of {\sf all decision problems $P\subseteq U_{\cal A}^\infty $ that are semi-decidable} by a BSS RAM in ${\sf M}_{\cal A}$ and ${\rm DEC}_{\cal A}$ is the class of {\sf all decidable problems} $P\subseteq U_{\cal A}^\infty $ for which both sets $P$ and its complement $U_{\cal A}^\infty\setminus P$ are in ${\rm SDEC}_{\cal A}$ so that ${\rm DEC}_{\cal A}={\rm SDEC}_{\cal A}\cap \{P\subseteq U_{\cal A}^\infty\mid U_{\cal A}^\infty\setminus P\in {\rm SDEC}_{\cal A}\}$ holds. For further classes, see \cite{GASS25A}. 

\vspace{0.1cm}

{\bf Configurations.} For any infinite 1-tape ${\cal A}$-machine ${\cal M}$, its {\sf overall state} at a given point of time is determined by a {\sf configuration} $(\ell, \nu_1,\ldots,\nu_{k_{\cal M}}, u_1,u_2,\ldots)$ that can also be denoted by $(\ell\,.\,\vec \nu\,.\,\bar u)$ where $\ell \in {\cal L}_{\cal M}$, $\vec \nu=( \nu_1,\ldots,\nu_{k_{\cal M}})\in \mathbb{N}_+^{k_{\cal M}}$, and $\bar u=( u_1,u_2,\ldots)\in U_{\cal A}^\omega$ hold. The transformation of one configuration $con_t$ of ${\cal M}$ into the next configuration $con_{t+1}$ of ${\cal M}$ is defined by the {\sf transition system} $({\sf S}_{\cal M},\to_{\cal M})$ that determines how an instruction can {\sf be processed }or {\sf executed} by ${\cal M}$. The modification of --- a finite number of components of --- $con_t\in {\sf S}_{\cal M}$ by applying $\to_{\cal M}$ is defined for deterministic machines in \cite[pp.\,583--584]{GASS20} and for various types of machines in \cite[Overviews 20 and 21]{GASS25A} and \cite[Overview 6]{GASS25B}. 

\vspace{0.1cm}

{\bf Sets ${\sf Con}_ {{\cal M}, {\sf i}}$ of sequences of configurations.} For any 1-tape ${\cal A}$-machine ${\cal M}$ and any ${\sf i}\in {\sf I}_{\cal M}$, let  ${\sf Con}_ {{\cal M}, {\sf i}}$ be the set all infinite sequences $(con_t)_{t\geq 1}$ of configurations given by $con_t=(\ell_t\,.\,\vec \nu^{\,(t)}\,.\, \bar u^{(t)})\in {\sf S}_{\cal M}$ and \big[non-\big]deterministically generated by repeatedly applying the total single-valued or multi-valued function $\to_{\cal M}$ such that, according to the types of BSS RAMs considered here,

\vspace{0.1cm}

\hspace{0.5cm} \begin{tabular}{lrcllclcl} 
&($con_1$&$ \!\!\textcolor{blue}{=} \!\!$&${\rm Input}_{\cal M} ({\sf i})$&\,\,and\quad &$con_t$&$ \!\!\textcolor{red}{ \to_{\cal M} } \!\!$&$con_{t+1}$)\\or\quad &
($ ({\sf i},con_1)$&$ \!\!\textcolor{blue}{\in} \!\!$&${\rm Input}_{\cal M}$&\,\,and\quad &$con_t$&$ \!\!\textcolor{red}{\to_{\cal M}} \!\!$&$con_{t+1}$)\\or \quad &
($con_1$&$ \!\!\textcolor{blue}{=} \!\!$&${\rm Input}_{\cal M} ({\sf i})$&\,\,and\quad &$con_t$&$ \!\!\textcolor{red}{\genfrac{}{}{0pt}{2}{\longrightarrow} {\longrightarrow}_ {\cal M}} \!\!$&$con_{t+1}$)
\end{tabular}

\vspace{0.1cm}

\noindent  hold.  Any   $((\ell_t\,.\,\vec \nu^{\,(t)}\,.\, \bar u^{(t)}))_{t\geq 1} \in {\sf Con}_ {{\cal M}, {\sf i}}$   in turn determines a finite {\em computation path} $(\ell_t)_{t=1..s} $ with $\ell_t\not=\ell_{{\cal P}_{\cal M}}$ for $t<s$ and $\ell_s=\ell_{{\cal P}_{\cal M}}$ or an infinite {\em computation path} $(\ell_t)_{t\geq 1} $ with $\ell_t\not=\ell_{{\cal P}_{\cal M}}$ for all $t\geq 1$.  The transition systems  of  BSS RAMs  in \textcolor{brown}{${\sf M}_{\cal A}^{(d)}$}  or \textcolor{brown}{$({\sf M}_{\cal A}^{(d)})^{\rm ND}$} and other   \text{\big[non-\big]}deterministic  \textcolor{brown}{$d$}-tape  machines (for definitions, see \cite[pp.\,591--592]{GASS20} and \cite[pp.\,5--6]{GASS25B}) and ${\sf Con}_ {{\cal M}, {\sf i}}$ are defined by using the corresponding adapted functions for changing configurations $(\ell_t,(\vec \nu^{\,(t,j)}\,.\,\bar u^{(t,j)})_{j=1..d})$.

\vspace{0.2cm}

{\bf Computation paths $(\ell_1,\ell_2,\ldots \big[,\ell_s\big]) \in {\cal L}_{\cal M}^\infty \cup{\cal L}_{\cal M}^\omega$.} For any  machine ${ \cal M}$ and any  $\,{\sf i}\in {\sf I}_{\cal M}$, any infinite sequence of configurations in ${\sf Con}_ {{\cal M}, {\sf i}}$  leads to  exactly one maximal sequence of labels that form a {\sf computation path} of ${\cal M}$. These paths are in general denoted by $B^*, B_0^*, \ldots$ (derived from the notion {\sf B}erechnungspfad). For a deterministic machine ${\cal M}$ and any $\,{\sf i}\in {\sf I}_{\cal M}$, the infinite sequence of configurations in ${\sf Con}_ {{\cal M}, {\sf i}}$ is uniquely determined and determines uniquely one computation path, i.e., exactly one sequence whose components are the labels of the instructions executed in this order by ${\cal M}$ on ${\sf i}$ until the stop instruction is reached. Less formally, we say that ${\cal M}$ {\em goes through this computation path for ${\sf i}$} or {\em traverses it for ${\sf i}$} or {\em ${\sf i}$ traverses the path of ${\cal M}$}. Sometimes, we also use other similar formulations. For two sequences $(con_t)_{t\geq 1}\in {\sf Con}_ {{\cal M}, {\sf i}}$  and $(con_t')_{t\geq 1} \in {\sf Con}_ {{\cal M}, {\sf i}'}$  resulting from the computation of a 1-tape machine ${\cal M}$ on two inputs ${\sf i}$ and ${\sf i}'$ and being  given by  $con_t=(\ell_t\,.\,\vec \nu^{\,(t)}\,.\, \bar u^{(t)})$ and $con_t'=(\ell_ t'\,.\,\vec \mu^{\,(t)}\,.\, \bar z^{(t)})$, both sequences $(\ell_t)_{t\geq 1} $ and $(\ell_t')_{t\geq 1} $ can match. For many inputs, ${\cal M}$ can traverse the same computation path. All computation paths of ${\cal M}$ together form a tree with label 1 as root.  For a  non-deterministic machine ${\cal M}$ and  a single  input ${\sf i}\in {\sf I}_{\cal M}$,  the  guessing of individuals or labels can lead to different  sequences in ${\sf Con}_ {{\cal M}, {\sf i}} $ and  to  different  computation paths of ${\cal M}$ that could be  traversed by ${\sf i}$. Such an ${\sf i}$ belongs to $H_{\cal M}$ if one of these computation paths finally contains $\ell_{{\cal P}_{\cal M}}$. We will give more details later. The guessing of labels leads in general to a {\sf derivation} or {\sf transformation} tree of configurations (for a subtree of such a tree, see   Figure \ref{Ableitung}).

The evaluation of computation paths plays a major role in different areas of computability theory and helps to classify computational problems and algorithms. Because of their great importance for theoretical approaches in computer science and for dealing with the foundations of logic programming, classical and other logics, and mathematics, there will be an introduction with more details --- with adapted, less informal, and more precise definitions --- in a {\sf special part}. The evaluation of {\sf computation} or {\sf decision tree models} includes the investigation of similar paths. They are often used for evaluating the time complexity of algorithms at a higher level of abstraction. This kind of computational models can be found --- for instance --- in practice-oriented introductions to computational geometry (such as in \cite{PS85}) and --- to mention further applications --- in the evaluation of the BSS model and Koiran's weak model over the real numbers (cf.\,\,e.g.\,\,\cite{CSS94}). Examples for the evaluation of the capabilities of various ${\cal A}$-machines and their complexity with the help of program or computation paths can also be found in several proofs in \cite{GASS97, GASS01, GASS08A, GASS08C,GASS09A, GASS09B, GASS10, GASS13, GASS17}. 

\vspace{0.2cm}

{\bf Relativized identity.}  We think that the {\sf identity} between objects should imply equality and uniqueness and that the {\sf equality} means that the objects are equal with respect to (a set of) their attributes, components, and properties. In this sense, the {\sf identity} describes an absolute equality and the notion {\sf equality} stands only for an equivalence related to an underlying universe. We consider each equality relation denoted by = as a form of identity within our theories, even if it cannot be assumed that this holds in any other \big[meta-\big]theory or world. 

\vspace{0.2cm}

{\bf Structures with identity relation.} The class ${\sf Struc}$ contains all first-order structures of the form $(U; (c_i)_{i\in N_1}; (f_i)_{i\in N_2}; (r_i)_{i\in N_3}) $. For ${\cal A}\in {\sf Struc}$, the identity relation ${\rm id}_{U_{\cal A}}$ is given by ${\rm id}_{U_{\cal A}}=\{(x,x)\mid x\in U_{\cal A}\}$ (for this definition, see \cite[p.\,7]{Mendelson}). Here, for $x,y\in U_{\cal A}$, we say that {\em $x$ is equal to $y$} {\sf if and only if} $x$ and $y$ are {\em identical} ({\em over ${\cal A}$}) which means that $(x,y) \in {\rm id}_{U_{\cal A}}$ holds. {\sf The character $=$ is used as symbol} for ${\rm id}_{U_{\cal A}}$ and {\sf as name} and {\sf notation} for ${\rm id}_{U_{\cal A}}$. By \cite[Agreement 2.1, p.\,6]{GASS25A}, ${\cal A}$-machines may perform an equality test only if the underlying structure ${\cal A}$ contains ${\rm id}_{U_{\cal A}}$. If $(U_{\cal A};c_{\alpha_1},\ldots,c_{\alpha_{n_1}};f_1,\ldots, f_{n_2}; r_1,\ldots,r_{n_3-1},=)$ is a reduct of ${\cal A}$ whose signature $\sigma$ is of the form $(n_1;m_1, \dots,m_{n_2}; k_1, \dots,k_{n_3-1},2)$, $r_{n_3}$ is the relation ${\rm id}_{U_{\cal A}}$, and $k_{n_3}=2$ holds, then  we also use (1) to (3).

(1) The memory locations for storing elements of $U_{\cal A}$ are not a feature of these individuals within ${\cal A}$. Two registers of an ${\cal A}$-machine can contain two objects that are identical over ${\cal A}$. For more, see Agreement \ref{agree}.

(2) An ${\cal A}$-machine ${\cal M}$ whose program is a $\sigma$-program can evaluate an equality test only if an instruction of the form {\sf if $r_{n_3}^{2}(Z_{j_1},Z_{j_{2}})$ then goto $\ell _1$ else goto $\ell _2$}
 belongs to ${\cal P}_{\cal M}$. Alternatively, the {\sf symbol} $\textcolor{red}{\,=\,}$ can be used instead of $r_{n_3}^{2}$ and the infix notation $Z_{j_1}=Z_{j_{2}}$ for the first-order literal $r_{n_3}^{2}(Z_{j_1},Z_{j_{2}})$.

(3) Let $\sigma^{\circ}$ be the signature $(m_1, \dots,m_{n_2}; k_1, \dots,k_{n_3-1})$ derived from $\sigma$ (cf.\,\,\cite[pp.\,14\,ff.]{A75}). The axioms of identity theory for first-order logic summarized in Overview \ref{IdentAxioms} hold in ${\cal A}$. ${\cal A}$ enables the {\sf $\sigma^{\circ}$-compliant interpretation} of the relation and operation symbols and the {\sf assignments} of individuals in $U_{\cal A}$ to the variables $X,Y,\ldots$. Each symbol $c_i^0$ (in a program) can be interpreted as the {\sf constant} $c_{\alpha_i}$ ($i\leq n_1$) and each constant in ${\cal A}$ is, by definition, an {\sf individual} that can also be assigned to $X,Y,\ldots$. The symbol $=$ can only be interpreted as ${\rm id}_{U_{\cal A}}$.
\begin{overview}[The system $axi^{\sigma^{\circ}}$ of axioms, cf. e.g. Asser {\cite[p.\,48]{A75}}]\label{IdentAxioms}

\hfill

\nopagebreak 

\noindent \fbox{\parbox{11.8cm}{
\noindent {\small 
\begin{tabular}{l} 
$\textcolor{blue}{X}=\textcolor{blue}{X}$,\vspace{0.1cm} \hfill {\scriptsize (A1)}\\
$\textcolor{blue}{X}=\textcolor{magenta}{Y} \to (\textcolor{blue}{X}=V\to \textcolor{magenta}{Y}=V)$, \hfill {\scriptsize (A2)}\vspace{0.1cm}\\
$\textcolor{blue}{X}=\textcolor{magenta}{Y} $\\\hspace{0.4cm}$\to f_i^{m_i}(X_1,\ldots,X_j,\textcolor{blue}{X},X_{j+1},\ldots,X_{m_i})= f_i^{m_i}(X_1,\ldots,X_j,\textcolor{magenta}{Y},X_{j+1},\ldots,X_{m_i})$,\\ \hfill {\scriptsize (A3/$i$/$j$)}\\
$\textcolor{blue}{X}=\textcolor{magenta}{Y}$\\\hfill
$\to (r_i^{k_i}(X_1,\ldots,X_j,\textcolor{blue}{X},X_{j+1},\ldots,X_{k_i})\to r_i^{k_i}(X_1,\ldots,X_j,\textcolor{magenta}{Y},X_{j+1},\ldots,X_{k_i}))$. \\\hfill{\scriptsize (A4/$i$/$j$)}\\
\hline {\small \hfill Here, $X,Y,V,X_1,X_2,\ldots $ are variables for individuals in $U_{\cal A}$ and we have \hspace*{0.18cm}}\\

 \hfill {\small $j\leq m_i-1$ and $j\leq k_i-1$ for all $i\leq n_2$ and $i\leq n_3-1$, respectively.} \end{tabular}}
}}
\end{overview} 
We assume that the corresponding rules and principles {\sf also hold for the identity relation ${\rm id}$ in our metatheory}. Two tuples are identical if their components are pairwise equal. Inputs of a deterministic BSS RAM that are identical with respect to ${\rm id}_{U_{\cal A}}$ lead to the same initial configuration, and so on. 

\vspace{0.2cm}

{\bf The class ${\sf Struc}_{c_1,c_2}$ and the machine constants $c_1$ and $c_2$.} Let {\sf $c_1$ and $c_2$} be two different objects which means that $c_1\not=c_2 $ and thus $(c_1,c_2)\not \in {\rm id}$ hold and let ${\sf Struc}_{c_1,c_2}$ be the subclass of ${\sf Struc}$ consisting of all structures of the form $(U; (c_i)_{i\in N_1}; (f_i)_{i\in N_2}; (r_i)_{i\in N_3})$ that contain $c_1$ and $c_2$ as constants with $(c_1,c_2)\in U_{\cal A}^2\setminus {\rm id}_{\cal A}$ (which means $\{1,2\} \subseteq N_1$).   For interpreting  $c_1^0$ in the program of any ${\cal A}$-machine ${\cal M}$  by means of a reduct ${\cal B}_{\cal M}$ of ${\cal A}\in {\sf Struc}_{c_1,c_2}$, we use only $c_1$   which means  that we use $\alpha_1=1$ and $c_1$ is the machine constant $c_{\alpha_1}$  in ${\cal B}_{\cal M}$.  For  interpreting $c_2^0$, we use only $\alpha_2=1$ and  $c_{\alpha_2}$.  For discussions on the usefulness of two distinguishable constants, see also \cite{{H98},GASS01, Prun02, P95}, and others. 

\vspace{0.2cm}

{\bf The functions $\bar \chi_P$ and $\chi_P$.} Let ${\cal A}\in {\sf Struc}$ and $P\subseteq U_{\cal A}^\infty$. If ${\cal A}$ contains the constant $c_1 $, then the {\sf partial characteristic function} $\bar \chi_P: \, \subseteq U_{\cal A}^\infty\to \{c_1\}$ and, for ${\cal A}\in {\sf Struc}_{c_1,c_2}$, the {\sf characteristic function} $\chi_P:U_{\cal A}^\infty\to \{c_1,c_2\}$ are defined by 

\vspace{0.1cm}

$\begin{array}{rcrcl}\begin{array}{ll}\bar \chi_P(\vec x)= c_1 & \mbox{ if }\vec x \in P,\\\bar \chi_P(\vec x)\uparrow& \mbox{ otherwise} \end{array}\hspace{0.8cm}\mbox{ and }\hspace{0.8cm}
\chi_P(\vec x) =\left \{\begin{array}{ll} c_1 & \mbox{ if }\vec x \in P,\\ c_2& \mbox{ otherwise, } \end{array}\right.\\
\end{array}$

\vspace{0.1cm}

\noindent respectively, for all $\vec x\in U_{\cal A}^\infty$. Such a $\bar \chi_P$ is computable over ${\cal A}$ (i.e.\,\,{\sf  by  means of a BSS RAM ${\cal M}$ in ${\sf M}_{\cal A}$}) if $P\in {\rm SDEC}_{\cal A}$ holds. The computation of $\chi_P$ by an ${\cal M} \in {\sf M}_{\cal A}$ means that $\vec x\in P$ is accepted by ${\cal M}$ and $\vec x\in U_{\cal A}^\infty \setminus P$ is rejected. 

\vspace{0.1cm}

In Section \ref{SectTwoConst}, we discuss the power of non-deterministic BSS RAMs over structures ${\cal A}\in {\sf Struc}_{c_1,c_2}$. To characterize the significance and usefulness of $c_1$ and $c_2$, we also deal with some typical relationships between the decidability of sets $P \subseteq U_{\cal A}^\infty$ and the computability of their characteristic functions $\chi_P$ over ${\cal A}$. For individual structures without identity relation, the strength of different non-determinisms can depend on the question of whether and in which way it is possible to distinguish between the constants. We compare a non-determinism caused by non-deterministic branchings, the digital non-determinism, and non-determinisms resulting from the use of the $\nu$-operator applied to oracle sets of tuples of constants for structures with decidable or semi-decidable identity or recognizable constants (where we use the term {\sf recognizable} in a less formal way), and the like. Finally, we give a summary, a small outlook, and some flowcharts. 

\setcounter{section}{5}

\section{Two constants and consequences}\label{SectTwoConst}

\subsection{Two constants as outputs and the decidability}\label{Section_6_1}

We will prove some useful relationships mentioned in \cite[pp.\,13--14]{GASS25A} for characterizing the decidability of problems over ${\cal A}\in {\sf Struc}_{c_1,c_2}$ (cf.\,\,Overview \ref{DecAndChar}). 

\subsubsection {On sufficient conditions for the decidability of problems} 
We will specify  some conditions under which  an algorithm for  computing  $\chi_P$ can be used to decide a problem $P$.
\begin{proposition}[Computability and identity imply decidability]\label{Propos1} \hfill 

\noindent Let  ${\cal A}\in {\sf Struc}_{c_1,c_2}$, let  ${\cal A}$ contain the identity relation ${\rm id}_{\cal A}$, and let $P\subseteq U_{\cal A}^\infty$. 

Then, the computability of $\chi_P$ by a BSS RAM in ${\sf M}_{\cal A}$ implies $P\in {\rm DEC}_{\cal A}$. 
\end{proposition}
{\bf Proof.} Let us assume that ${\cal M}\in {\sf M}_{\cal A}$ computes $\chi_P$ for $P\subseteq U_{\cal A}^\infty$. Let the program ${\cal P}_{\cal M}$ be of the form \textcolor{blue}{$(*)$} and let ${\cal P}_1$ be the following program segment. 
\begin{overview}[The program segment \textcolor{blue}{${\cal P}_1$}]\label{FirstInstr}

\hfill 

\nopagebreak 

\noindent \fbox{\parbox{11.8cm}{
\hspace{0.83cm} $\ell _{{\cal P}_{\cal M}}\,\,: Z_2:=c_1^0;$ 

\quad $\ell _{{\cal P}_{\cal M}}+1\,\,:$ {\sf if $Z_1=Z_2$ then goto \textcolor{blue}{$\ell _{{\cal P}_{\cal M}}+2$} else goto $\textcolor{blue}{\ell _{{\cal P}_{\cal M}}+1};$} 

\quad $\ell _{{\cal P}_{\cal M}}+2\,\, :$ {\sf stop}.}}
\end{overview}
The symbol $=$ is used as described above. Let ${\cal P}_2$ result from 
replacing the instruction labeled by $\ell _{{\cal P}_{\cal M}}+1$ in Overview \ref{FirstInstr} by the instruction in Overview \ref{SecondInstr}.

\begin{overview}[The second instruction of the program segment \textcolor{blue}{${\cal P}_2$}]\label{SecondInstr}

\hfill 

\nopagebreak 

\noindent \fbox{\parbox{11.8cm}{

\hspace{1.83cm} {\sf if $Z_1=Z_2$ then goto \textcolor{blue}{$\ell _{{\cal P}_{\cal M}}+1$} else goto $\textcolor{blue}{\ell _{{\cal P}_{\cal M}}+2}$}}}
\end{overview}
 For $\textcolor{blue}{i}\in \{1,2\}$, let

\vspace{0.1cm}

\quad 1 : {\sf instruction}$_1$;\quad\ldots; \quad $\ell _{{\cal P}_{\cal M}}-1$ : {\sf instruction}$_{\ell _{{\cal P}_{\cal M}}-1}$; \quad ${\cal P}_i$ \hfill \textcolor{blue}{$(*\langle {\cal P}_i)$}. 

\vspace{0.2cm}

\noindent be the program ${\cal P}_{{\cal M}_{\textcolor{blue}{i}}}$ of a BSS RAM ${\cal M}_i$ that results from replacing the program segment \fbox{$\ell _{{\cal P}_{\cal M}}$\,\,: {\sf stop}.} in \textcolor{blue}{$(*)$} by ${\cal P}_{\textcolor{blue}{i}}$. Then, $\ell _{{\cal P}_{{\cal M}_{\textcolor{blue}{i}}}}=\ell _{{\cal P}_{\cal M}}+2$ holds. 
${\cal M}_1$ semi-decides the decision problem $P$. ${\cal M}_2$ semi-decides $U_{\cal A}^\infty \setminus P$. \qed

\vspace{0.2cm}

For the use of program segments such as ${\cal P}_1$ and ${\cal P}_2$ considered above, see also the last instructions in Example \ref{CounterEx1} The next proposition is more general (cf.\,\,Remark \ref{fromCtoD}) and also of practical relevance (cf.\,\,Remark \ref{Pract}). 

\begin{remark}\label{Pract}In practical programming, it can be difficult to check the equality of certain objects and it cannot be taken for granted that the equality tests can be realized (for example, by means of $\,==$ or the like). In particular, the attempt to test the equality of objects such as real numbers must be avoided. For the use of the finite precision tests $<_k$ for comparing the real numbers, see also \cite{BH98}.\end{remark}

\begin{proposition}[Computability $\!$\&$\!$ recognizability imply decidability]\label{Propos2}Let ${\cal A}\in {\sf Struc}_{c_1,c_2}$, $\{c_1\}\in {\rm SDEC}_{\cal A}$, and $\{c_2\}\in {\rm SDEC}_{\cal A}$, and $P\subseteq U_{\cal A}^\infty$.

Then, the computability of $\chi_P$ by a BSS RAM in ${\sf M}_{\cal A}$ implies $P\in {\rm DEC}_{\cal A}$. 
\end{proposition}
{\bf Proof.} Let us assume that ${\cal M}\in {\sf M}_{\cal A}$ computes $\chi_P$ and let ${\cal P}_{\cal M}$ be given by \textcolor{blue}{$(*)$}. Moreover, let the programs of two BSS RAMs ${\cal M}'\in {\sf M}_{\cal A}$ and ${\cal M}'' \in {\sf M}_{\cal A}$ semi-deciding $\{c_1\}$ and $\{c_2\}$, respectively, be given by \textcolor{blue}{$(*')$} and \textcolor{blue}{$(*'')$}.

\begin{tabbing} 
 \hspace{0.7cm}$1$\=: {\sf instruction}$_{1}''$; \quad \ldots;\quad\=$\ell _{{\cal P}_{{\cal M}''}}-1$\=: {\sf instruction}$_{\ell _{{\cal P}_{{\cal M}''}}-1}$;\quad \=$\ell _{{\cal P}_{{\cal M}''}}$\=: {\sf stop}.\quad\,\,\=(i'') \= \kill
 \hspace{0.7cm}$1$\>: {\sf instruction}$_{1}'$; \quad \ldots;\quad\>$\ell _{{\cal P}_{{\cal M}'}}-1$\>: {\sf instruction}$_{\ell _{{\cal P}_{{\cal M}'}}-1}'$;\quad \>$\ell _{{\cal P}_{{\cal M}'}}$\>: {\sf stop}.\>\,\textcolor{blue}{$(*')$}\\
 \hspace{0.7cm}$1$\>: {\sf instruction}$_{1}''$; \quad \ldots;\quad\>$\ell _{{\cal P}_{{\cal M}''}}-1$\>: {\sf instruction}$_{\ell _{{\cal P}_{{\cal M}''}}-1}''$;\quad \>$\ell _{{\cal P}_{{\cal M}''}}$\>: {\sf stop}.\> \textcolor{blue}{$(*'')$}\end{tabbing} 

\noindent We want to define the programs of two BSS RAMs ${\cal M}_1\in {\sf M}_{\cal A}$ and ${\cal M}_2\in {\sf M}_{\cal A}$ semi-deciding $P\subseteq U_{\cal A}^\infty$ and $U_{\cal A}^\infty\setminus P$, respectively. 
For this purpose, we derive two program segments \textcolor{blue}{${\cal P}_1$} and \textcolor{blue}{${\cal P}_2$} from \textcolor{blue}{$(*')$} and \textcolor{blue}{$(*'')$} such that, for $i\in \{1,2\}$, the program ${\cal P}_{{\cal M}_i}$ has the form \textcolor{blue}{$(*\langle{\cal P}_i)$} and results from replacing the stop instruction in ${\cal P}_{\cal M }$ --- as in the proof of Proposition \ref{Propos1} --- by using $\textcolor{blue}{{\cal P}_i}$. Let ${\cal P}_{{\cal M}_1}$ be the program $\textcolor{blue}{{\cal P}_0{\cal P}_1}$ defined in Overview \ref{FromCompToSemi}. The first segment \textcolor{blue}{${\cal P}_0$} contains the first instructions of ${\cal P}_{\cal M}$. Then, $\textcolor{blue}{{\cal P}_0{\cal P}_1}$ allows to compute first the output of ${\cal M}$ in order to partially check then whether this value belongs to $\{c_1\}$ by using the segment $\textcolor{blue}{{\cal P}_1}$ for semi-deciding the set $\{c_1\}$ (restricted to the outputs of ${\cal M}$). This means that we want to apply $\textcolor{blue}{{\cal P}_1}$ to the outputs in $\{c_1,c_2\}$ of ${\cal M}$ as inputs.

\begin{overview}[The program $\,\textcolor{blue}{{\cal P}_0{\cal P}_{\rm init}' {\cal P}'}$ for semi-deciding $P$ by ${\cal M}_1$]\label{FromCompToSemi}

\hfill 

\nopagebreak 

\noindent \fbox{\parbox{11.8cm}{
\vspace{0.2cm}

$(1)$ Let the first program segment \textcolor{blue}{${\cal P}_0$} of ${\cal P}_{{\cal M}_1}$ be the program segment 

\vspace{0.1cm}

\hspace{1cm} $1: {\sf instruction}_{1}; \quad \ldots;\quad\ell _{{\cal P}_{\cal M}}-1:\, {\sf instruction}_{\ell _{{\cal P}_{\cal M}}-1};$

\vspace{0.2cm}

\hspace{0.55cm} of ${\cal P}_{\cal M}$ that leads to an intermediate result $z_1=c(Z_{c(I_1)})\in\{c_1,c_2\}$.

\vspace{0.2cm}

$(2')$ Let the second program segment \textcolor{blue}{${\cal P}_1$} of ${\cal P}_{{\cal M}_1}$ be divided into ${\cal P}_{\rm init}'$ and ${\cal P}'$.
\begin{itemize}
\item 
Let \textcolor{blue}{${\cal P}_{\rm init}'$} be given by 

\vspace{0.1cm}

\begin{tabular}{l}

\hspace{1.95cm} $\ell _{{\cal P}_{\cal M}}:\hspace{0.2cm} I_2:=1; \quad \ldots\quad ; \quad  \ell _{{\cal P}_{\cal M}} + k_{{\cal M}'}:\hspace{0.2cm}I_{k_{{\cal M}'}+2}:=1;$ \\

\noindent\hspace{0.3cm} $\ell _{{\cal P}_{\cal M}} + k_{{\cal M}'}+1:\hspace{0.2cm} {\sf init}(Z_{I_{k_{{\cal M}'+1}}});$\\
\end{tabular}

\vspace{0.1cm}

for providing the first part of ${\rm In}_{{\cal M}'}(z_1)$ for $z_1 \in \{c_1,c_2\}$.

\item Let \textcolor{blue}{${\cal P}'$} be derived from \textcolor{blue}{$(*')$} in such a way that 

\begin{itemize}\parskip -0.7mm
\item \parskip -0.7mm each label \fbox{$\ell$} $\in \{1,\ldots, \ell _{{\cal P}_{{\cal M}'}}\}$ in ${\cal P}_{{\cal M}'}$ is replaced by 

\vspace{0.2cm}

\hspace{0.12cm} $\ell _{{\cal P}_{\cal M}}+k_{{\cal M}'}+3+\ell$ \,, 

\vspace{0.1cm}

\item each {\sf instruction}$_k'$ of the form \fbox{$I_j:=I_j+1$} is extended  by  adding

\vspace{0.2cm}
\hspace{0.12cm}  $; \, {\sf init}(Z_{I_{k_{{\cal M}'+1}}})$ \,. 
\end{itemize}\end{itemize}
}} 
\end{overview}
\begin{overview}[The behavior of \textcolor{blue}{${\cal M}_1$} semi-deciding $P$]
\hfill

\nopagebreak 

\noindent \fbox{\parbox{11.8cm}{ 

\begin{tabular}{rcccl}Input& &$\!\!\!$Intermediate result$\!\!\!$& &Final behavior\\ 
$P \ni \vec x$& \quad {\LARGE $\longrightarrow$} &$ c_1$ &{\LARGE $\longrightarrow$} \quad & ${\cal M}_1$ {\sf halts}\\
$U_{\cal A}^\infty \setminus P\ni \vec x $& \quad {\LARGE $\longrightarrow$} & $c_2$ & {\LARGE $\longrightarrow$} \quad & ${\cal M}_1$ {\sf loops\, forever}\\
&{\scriptsize$\quad(1)$}&&{\scriptsize$ (2')\,\,$}\\
\end{tabular}
}}
\end{overview}
 \textcolor{blue}{${\cal P}_1$} is the program segment ${\cal P}_{\rm init}'{\cal P}'$. 
For getting the most important parts ${\rm In}_{{\cal M}'}(c_1)$ and ${\rm In}_{{\cal M}'}(c_2)$ of the initial configurations of ${\cal M}'$ that also contain the intermediate result $c_1$ or $c_2$, we introduce the program segment ${\cal P}_{\rm init}'$ (whose first label is $\ell _{{\cal P}_{\cal M}}$) as given in Overview \ref{FromCompToSemi}. Its definition has the following consequences. By using ${\cal P}_0$, i.e., by using all instructions of ${\cal P}_{\cal M}$ apart from the last instruction, the execution of $\textcolor{blue}{{\cal P}_{\rm init}'}$ can start with $c(I_1)=1$ and $c(Z_1)=z_1 \in \{c_1,c_2\}$ and provides moreover, by its instructions, $c(I_2)=\cdots =c(I_{k_{\cal M}'})=1$, $c(I_{k_{{\cal M}'}+1})=c(I_{k_{{\cal M}'}+2})=2$, and $c(Z_2)=z_1 $ before the execution of further instructions in ${\cal P}'$ for the simulation of the execution of ${\cal P}_{{\cal M}'}$ follows. The program segment $\textcolor{blue}{{\cal P}'}$ is derived from ${\cal P}_{{\cal M}'}$ as given in Overview \ref{FromCompToSemi}. The subprogram ${\sf init}(Z_{I_{k_{{\cal M}'+1}}})$ contains 3 instructions and helps step by step to replace the application of the input procedure of ${\cal M}'$. For its definition and its meaning, see \cite[p.\,22]{GASS25B}. Therefore, $c(Z_k)=z_1 $ holds for any $Z$-register $Z_{k}$ of ${\cal M}_1$ before its value is used in executing an instruction of ${\cal P}'$ for the first time (note, that this kind of extension was also done, for example, in Overview 19 in Part II\,a). Consequently, ${\cal M}_1$ halts on an input $\vec x\in U_{\cal A}^\infty$ only if ${\cal M}$ outputs $c_1$ for this input $\vec x$ which means that ${\cal M}_1$ semi-decides $P$. 

Let the program ${\cal P}_{{\cal M}_2}$ be divided into ${\cal P}_0$ and ${\cal P}_2$ and derived from ${\cal P}_{\cal M}$ and ${\cal P}_{{\cal M}''}$ in the same way as ${\cal P}_{{\cal M}_1}$ where \textcolor{blue}{${\cal P}_{\rm init}''$} and \textcolor{blue}{${\cal P}''$} result from replacing each prime $'$ in ${\cal P}_{\rm init}'$ and ${\cal P}'$ considered in Overview \ref{FromCompToSemi} by a double prime $''$.
Consequently, $\,\textcolor{blue}{{\cal P}_0{\cal P}_{\rm init}'' {\cal P}''}\,$ is the program of ${\cal M}_2$ for semi-deciding $U_{\cal A}^\infty\setminus P$. \qed
\begin{overview}[The behavior of \textcolor{blue}{${\cal M}_2$} semi-deciding $U_{\cal A}^\infty \setminus P$]
\hfill

\nopagebreak 

\noindent \fbox{\parbox{11.8cm}{ 

\begin{tabular}{rcccl}Input& &$\!\!\!$Intermediate result$\!\!\!$& &Final behavior\\ 
$P \ni \vec x$& \quad {\LARGE $\longrightarrow$} &$ c_1$ &{\LARGE $\longrightarrow$} \quad & ${\cal M}_2$ {\sf loops\, forever}\\
$U_{\cal A}^\infty \setminus P\ni \vec x $& \quad {\LARGE $\longrightarrow$} & $c_2$ & {\LARGE $\longrightarrow$} \quad & ${\cal M}_2$ {\sf halts}\\
&{\scriptsize$\quad (1)$}&&{\scriptsize$ (2'')\,\,$}\\
\end{tabular}
}}
\end{overview}

\vspace{0.1cm}

\begin{remark}\label{fromCtoD} Let ${\cal A}\in {\sf Struc}_{c_1,c_2}$ and $i\in \{1,2\}$. If $\,{\rm id}_{U_{\cal A}}$ is semi-decidable by an ${\cal M}\in{\sf M}_{\cal A}$ whose program is of the form $\textcolor{blue}{(*)}$, then $\{c_i\}$ is semi-decidable by means of a program that results from ${\cal P}_{\cal M}$ by adding, in front of ${\cal P}_{\cal M}$, the segment\, 

\vspace{0.1cm}

$1':$ {\sf if $I_1=I_2$ then goto $2'$ else goto $1';$} \quad $2':\, I_1:=I_1+1; \quad Z_2:=c_i^0;$ .

\vspace{0.1cm}

\noindent For more, see Proposition \ref{BasicProp}. For the background, see also Agreement \ref{agree}. \end{remark}

\begin{corollary}[Recognizable identity implies decidability]\label{CorollPropos2}\hfill

\noindent Let ${\cal A}$ be in $ {\sf Struc}_{c_1,c_2}$, ${\rm id}_{U_{\cal A}}\in {\rm SDEC}_{\cal A}$, and $P\subseteq U_{\cal A}^\infty$. 

Then, the computability of $\chi_P$ by a BSS RAM in ${\sf M}_{\cal A}$ implies $P\in {\rm DEC}_{\cal A}$. 
\end{corollary}

\subsubsection {Computable characteristic functions for decidable problems}

\begin{figure}[h]
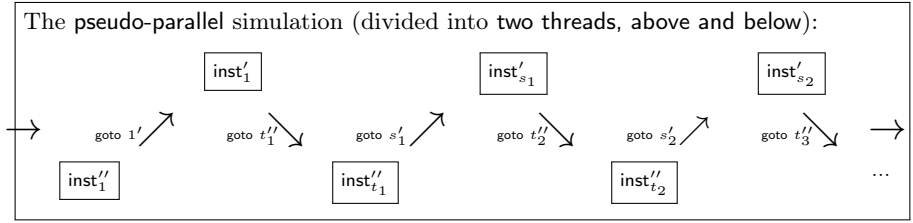

\noindent\hspace{0.4mm}\fbox{\parbox{11.67cm}{
{\small The {\sf pseudo-parallel} simulation (divided into {\sf two threads, above and below}):}

\vspace{0.2cm}

{\scriptsize

 \noindent \hspace{2.4cm}\fbox{{\sf inst}$'_1$} \hspace{2.8cm}\fbox{{\sf inst}$'_{s_1}$} \hspace{2.7cm}\fbox{{\sf inst}$'_{s_2}$} \hspace*{0.8cm}

\vspace{0.2cm}

\noindent{\Large$\!\!\!\rightarrow$}$\!\!$ \hspace{0.7cm}{\tiny {\sf goto} $1'\!\!$}{\Large $\nearrow$} \hspace{0.59cm}{\tiny {\sf goto} $t_1''\!\!\!\!$}{\Large$\searrow$} \hspace{0.59cm}{\tiny {\sf goto} $s_1'$}{\Large $\nearrow$} \hspace{0.59cm}{\tiny {\sf goto} $t_2''\!\!\!\!$}{\Large$\searrow$} \hspace{0.59cm}{\tiny {\sf goto} $s_2'$}{\large $\nearrow$} \hspace{0.59cm}{\tiny {\sf goto} $t_3''\!\!\!\!$}{\Large$\searrow$}\hspace{0.46cm} {\Large$\!\!\rightarrow$}$\!\!\!\!\!\!\!\!\!\!\!$ 

\vspace{0.2cm}

\hspace{0.38cm} \fbox{{\sf inst}$''_1$} \hspace{2.73cm}\fbox{{\sf inst}$''_{t_1}$} \hspace{2.73cm}\fbox{{\sf inst}$''_{t_2}$}} \hspace{2.46cm}$^ {\hdots}$

\vspace{0.1cm} }}

\caption{From deciding whether $\vec x \in P$ or not to the computation of $\chi_P(\vec x)$}
\label{ParaSim}
\end{figure}

The possibility to simulate two machines in a pseudo-parallel way (see Figure \ref{ParaSim}) leads to the following proposition. For proving this, we will derive the program of a 3-tape BSS RAM ${\cal M}$ from the programs of two BSS RAMs ${\cal M}'$ and ${\cal M}''$.

\begin{proposition}[Decidability implies computability]\label{Propos3}\hfill

\noindent Let ${\cal A}\in\! {\sf Struc}_{c_1,c_2}$. 

For any problem $P$ in $ {\rm DEC}_{\cal A}$, $\chi_P$ is computable by a BSS RAM in $ {\sf M}_{\cal A}$.
\end{proposition}
{\bf Proof.} Let $P\subseteq U_{\cal A}^\infty$ and let ${\cal M}'$ and ${\cal M}''$ be two 1-tape BSS RAMs in $ {\sf M}_{\cal A}$ that semi-decide the decision problems $P$ and $U_{\cal A}^\infty\setminus P$, respectively. Let the programs ${\cal P}_{{\cal M}'}$ and ${\cal P}_{{\cal M}''}$ be of the form \textcolor{blue}{$(*')$} and \textcolor{blue}{$(*'')$}, respectively, as given in the proof of Proposition \ref{Propos2} and let $d'=1$ and $d''=2$. Under these assumptions, let ${\cal M}$ be a deterministic 3-tape BSS RAM in ${\sf M}^{(3)}_{\cal A}$ using the first tape $Z_{d',1},Z_{d',2}, \ldots $ to simulate ${\cal M}'$ and the second tape $Z_{d'',1},Z_{d'',2}, \ldots $ to simulate ${\cal M}''$. For any input $\vec x\in U_{\cal A}^\infty$, we want to get the answer $c_1$ or $c_2$ (for {\sf yes} or {\sf no}) as output given by ${\cal M}$ that says in this way whether $\vec x$ belongs to $P$ or not. Only one of the machines ${\cal M}'$ and ${\cal M}''$ halts on $\vec x$. Consequently, it is not possible to simulate first only ${\cal M}'$ since it could be that ${\cal M}'(\vec x)\uparrow$ holds. A possible solution to the problem is to use a procedure that allows the simulation of the execution of exactly one instruction of one of the machines before the then current instruction of the other machine is processed next. 
For this purpose, two index registers $I_{3,d'}$ and $I_{3,d''}$ can be used for storing all information about which instruction of ${\cal P}_{{\cal M}'}$ and ${\cal P}_{{\cal M}''}$, respectively, is the next that is to be processed after the simulation of the execution of an instruction of the other program.
For $\vec x\in U_{\cal A}^\infty$, the start configuration ${\rm Input}_{\cal M}(\vec x)$ is  given by $(1,(\vec \nu^{\,(j)}\,.\,\bar u^{(j)})_{j=1..3})$ with $\vec \nu^{\,(1)} = (n,1,\ldots, 1)$, $\bar u^{(1)} = (x_1, \ldots ,x_n, x_n,\ldots )$, and $\vec \nu^{\,(j)}=(1,\ldots , 1) $ and $\bar u^{(j)} =(x_n,x_n,\ldots )$ for $j\in\{2,3\}$. The output (after reaching the instruction labeled by $\ell _{{\cal P}_{\cal M}}$) is given by ${\rm Out}_{\cal M}((\vec \nu^{\,(j)}\,.\,\bar u^{(j)})_{j=1..3})= (u_{1,1}, \ldots,u_{1,\nu_{1,1}})$. If, for input $\vec x$, ${\cal M}'$ traverses a computation path given by $(1,s_1,s_2,\ldots)$ and executes a sequence 

\vspace{0.05cm}

{\sf instruction}$'_1$, {\sf instruction}$'_{s_1}$, {\sf instruction}$'_{s_2}$, $\ldots$ 

\vspace{0.05cm}

\noindent of instructions and ${\cal M}''$ executes another sequence of instructions 

\vspace{0.05cm}

{\sf instruction}$''_1$, {\sf instruction}$''_{t_1}$, {\sf instruction}$''_{t_2}$, $\ldots$ 

 \vspace{0.05cm}

 \noindent for some computation path $(1,t_1,t_2,\ldots)$, then the program ${\cal P}_{\cal M}$ of ${\cal M}$ defined in Overviews \ref{CharactchiP1}, \ref{CharactchiP2}, and \ref{CharactchiP3} allows the {\sf pseudo-parallel simulation} of ${\cal M}'$ and ${\cal M}''$ on $\vec x$ as illustrated in Figure \ref{ParaSim}. All subprograms ${\sf inst}'_1, {\sf inst}'_2, \ldots$ and ${\sf inst}''_{1}, \ldots$ of ${\cal P}_{\cal M}$ marked with a prime $'$ or a double prime $''$ and following new labels in $\{1\,',\ldots,(\ell _{{\cal P}_{{\cal M}'}} - 1)\,'\} \cup \{1\,'',\ldots, (\ell _{{\cal P}_{{\cal M}''}} - 1)\,''\}$ are derived from ${\cal P}_{{\cal M}'}$ and ${\cal P}_{{\cal M}''}$ and intended to simulate ${\cal M}'$ and ${\cal M}''$ step by step for computing $\chi_P$. \qed
\begin{overview}
[The program ${\cal P}_{\cal M}$ incl. \textcolor{blue}{${\cal P}_1$} and \textcolor{blue}{${\cal P}_2$} for computing $\chi_P$]\label{CharactchiP1}

\hfill

\nopagebreak 

\noindent \fbox{\parbox{11.8cm}{\begin{tabular}{rlrl}

$1 :$&$\!\!\!(Z_{2,1},\ldots,Z_{2, I_{2,1}}) := (Z_{1,1},\ldots,Z_{1, I_{1,1}});$\\

\textcolor{blue}{$\tilde 1\,':$}&$\!\!\!${\sf if $I_{3,2}=i$ then goto \textcolor{blue}{$i''$}$;$}\\

\textcolor{blue}{$\tilde 1\,'':$}&$\!\!\!${\sf if $I_{3,1}=i$ then goto \textcolor{blue}{$i'$}$;$}\\

\textcolor{blue}{$1\,' :$}&$\!\!\!${\sf inst}$'_1;\, \,\ldots;\,\hspace{0.34cm}$\textcolor{blue}{$( \ell _{{\cal P}_{{\cal M}'}}-1)\,'$}$\, :$ {\sf inst}$'_{ \ell _{{\cal P}_{{\cal M}'}}-1};$ $\,\,\,\, \ell\,' _{{\cal P}_{{\cal M}'}}\,: Z_{1,1} := c_1^0;\,\,I_{1,1}= 1;$\\
$1\,''':$ &$\!\!\!${\sf goto} $\ell _{{\cal P}_{\cal M}};$\\ 

\textcolor{blue}{$1\,'':$}&$\!\!\!${\sf inst}$''_1; \,\,\ldots;\,\textcolor{blue}{( \ell _{{\cal P}_{{\cal M}''}}-1)\,''}: $ {\sf inst}$''_{ \ell _{{\cal P}_{{\cal M}''}}-1}; \hspace{0.03cm}\,\,\,\ell\,'' _{{\cal P}_{{\cal M}''}}: Z_{1,1} := c_2^0; \,\,I_{1,1}= 1;$\\

$\ell _{{\cal P}_{\cal M}}:$&$\!\!\!${\sf stop}.\hfill {\small For pseudo instructions, see also \cite{GASS25B}.}

\end{tabular}
}}
\end{overview}

\begin{overview}
[Program segment ${\cal P}_1$ in ${\cal P}_{\cal M}$ for computing $\chi_P$]\label{CharactchiP2}

\hfill

\nopagebreak 

\noindent \fbox{\parbox{11.8cm}{

1) Let ${\cal P}_{{\cal M}'}$ be the program of the form \textcolor{blue}{$(*')$} for semi-deciding $P$ and let \textcolor{blue}{${\cal P}_1$} be 

\vspace{0.1cm}

\hspace{2cm}\textcolor{blue}{$1\,'$} $:$ {\sf inst}$'_1;$\,\, \textcolor{blue}{$2$}$\,':$ {\sf inst}$'_2; \,\ldots;$\hspace{0.22cm}\textcolor{blue}{$( \ell _{{\cal P}_{{\cal M}'}}-1)\,'$}$:$ {\sf inst}$'_{ \ell _{{\cal P}_{{\cal M}'}}-1};$

\vspace{0.1cm}

where each subprogram {\sf inst}$'_k$ results from {\sf instruction}$_k'$ by
\parskip -0.7mm
\begin{itemize} 
\item\parskip -0.7mm replacing
\begin{itemize} 
\item\parskip -0.7mm 
 each symbol $Z_i$ by $Z_{d',i}$ , 
\item each symbol $I_i$ by $I_{d',i}$ 
\end{itemize} 
 \item and adding
 $;$ and the subprogram \,\,$I_{3,d'}:= k+1;$ {\sf goto \textcolor{blue}{$\tilde 1'$}}
\end{itemize} 

if {\sf instruction}$_k'$ is not of type (4) or (5), and otherwise it is given as follows.
 
\parskip -0.9mm 
\hfill {\small Here, $d'=1$.}
}}

\nopagebreak 

\noindent \fbox{\parbox{11.8cm}{

If {\sf instruction}$'_k$ has the form  \fbox{\sf if $cond$ then goto $\ell _1$ else goto $\ell _2$}, then let {\sf inst}$'_k$ be the subprogram

\quad {\sf \hspace{1.2cm} if $cond\,'$ then goto $(k')^{*}$ else goto $(k')^{**};$}

\quad $(k')^{*} \hspace{0.17cm}: \,$ $I_{3,d'}:= \ell_1;$ \hspace{0.15cm} {\sf goto \textcolor{blue}{$\tilde 1'$}$;$}

\quad $(k')^{**} \hspace{0.03cm}: \,$ $I_{3,d'}:= \ell_2;$ \hspace{0.15cm} {\sf goto \textcolor{blue}{$\tilde 1'$}}

where $\,cond\,'$ results from $\,cond\,$ by replacing $Z_i$ by $Z_{d',i}$ and $I_i$ by $I_{d',i}$.

\vspace{0.1cm}
}}
\end{overview}
\begin{overview}
[Program segment ${\cal P}_2$ in ${\cal P}_{\cal M}$ for computing $\chi_P$]\label{CharactchiP3}
\hfill

\nopagebreak

\noindent \fbox{\parbox{11.8cm}{
2) Let ${\cal P}_{{\cal M}''}$ be the program of the form \textcolor{blue}{$(*'')$} for semi-deciding $U_{\cal A}^\infty\setminus P$ and let 

\vspace{0.1cm}
\hspace{2cm}\textcolor{blue}{$1$}$\,'':$ {\sf inst}$''_1;$\, \textcolor{blue}{$2\,''$}$: $ {\sf inst}$''_2; \,\,\ldots;\, \,$\textcolor{blue}{$( \ell _{{\cal P}_{{\cal M}''}}-1)\,''$}$: $ {\sf inst}$''_{ \ell _{{\cal P}_{{\cal M}''}}-1};$

\vspace{0.1cm}
be the program segment \textcolor{blue}{${\cal P}_2$} in ${\cal P}_{\cal M}$ obtained when each {\sf inst}$''_k$ is defined analogously to {\sf inst}$'_k$ where each $\,'$ is replaced by a double prime\, $''$ .\hspace{0.3cm} {\small Here, $d''=2$.}
}}\end{overview}

\subsection{Two constants and their recognizability}

The considerations in Section \ref{Section_6_1} suggest that the computational power of BSS RAMs can depend on whether the identity (relation) can be used or not. In Section \ref{BasicPropSec}, we describe some consequences of the \big[semi-\big]decidability of ${\rm id}_{\cal A}$ in more detail. Sections \ref{FromSemiToDecSecI} and \ref{FromSemiToDecSecII} deal with the decidability of some semi-decidable sets. For a summary, see Overviews \ref{BasicImplicI} and \ref{BasicImplicII}.

\subsubsection{On sufficient conditions for recognizing the constants}\label{BasicPropSec}

\begin{proposition}[Basic properties: Identity and constants]\label{BasicProp}\hfill

\noindent
Let ${\cal A}$ be a structure in $ {\sf Struc}_{c_1,c_2}$ and the properties $(a)$ to $(e)$ be given.

\quad \begin{tabular}{cllr}
$(a)$&${\cal A}$ contains ${\rm id}_{U_{\cal A}}$,\\
 $(b)$&${\rm id}_{U_{\cal A}}\in {\rm DEC}_{\cal A}$, &\qquad \,\qquad$(d)$&$\{c_1\}, \{c_2\}\in {\rm SDEC}_{\cal A}$,\\
 $(c)$&${\rm id}_{U_{\cal A}}\in {\rm SDEC}_{\cal A}$, &\qquad \,\qquad$(e)$&$\{c_1, \, c_2\}^\infty\in {\rm SDEC}_{\cal A}$.\\
\end{tabular}

\vspace{0.15cm}
 Then, $(a)$ implies $(b)$, $(b)$ implies $(c)$,  $(c)$ implies $(d)$, and  $(d)$ implies $(e)$.
\end{proposition}
{\bf Proof.}
 Here, we only look at the last implication in more detail. If ${\cal M}'$ and ${\cal M}''$ are machines in ${\sf M}_{\cal A}$ semi-deciding $\{c_1\}$ and $\{c_2\}$, respectively, then a program defined analogously to the program ${\cal P}_{\cal M}$ in Overview \ref{CharactchiP1} allows to semi-decide $\{c_1,c_2\}^\infty$. In fact, no changes need to be made. However, the four instructions in ${\cal P}_{\cal M}$ for determining the output $c_1$ or $c_2$ of length 1 and the labels $ \ell\,' _{{\cal P}_{{\cal M}'}}$ and $ \ell\,'' _{{\cal P}_{{\cal M}''}}$ can be omitted. It is only important that ${\cal M}$ halts if one of the labels $\ell\,'_{{\cal P}_{{\cal M}'}}$ or $\,\ell\,''_{{\cal P}_{{\cal M}''}}$ could be reached and sufficient that (instead of these labels) one of the labels $1\,'''$ and $\,\ell_{{\cal P}_{\cal M}}$, respectively, is immediately and directly reachable. For the relationship between $(c)$ and $(d)$, see also Remark \ref{fromCtoD}. 
\qed

\vspace{0.cm}

We will now give examples to show that three of the reverse implications do not have to be right
 and that the performance of machines over ${\cal A}$ without \big[semi-decidable\big] identity can thus follow different standards.
\begin{example}[A structure with $(b)$ in which $(a)$ does not hold]\label{CounterEx1}\hfill 

\noindent  Let $a$ and $b$ be two different objects which means $a\not = b$ and $(a,b)\not \in {\rm id}$ and let $ \{a,b\}^*$ be the smallest set $U$ with the following properties. (1) $U$ contains the empty string $\varepsilon$. (2) If $ w\in U$ and $x\in \{a,b\}$ hold, then the uniquely determined result $ wx$ of the concatenation of $w$ and $x$ belongs to $U$. $U$ is thus the set of all finite strings (including the empty string) over the alphabet $\{a,b\}$. 

Let ${\cal A}_1$ be the structure $(\{a,b\}^*;\,a,b,\varepsilon\,; sub_{\rm l}, sub_{\rm r};\, r_1)$ of signature $(3;1,1;2)$ in ${\sf Struc}_{a,b}$ that contains three constants $a$, $b$, and $\varepsilon$, the binary relation $r_1$ that is a subset of the identity relation ${\rm id}_{U_{{\cal A}_1}}$ and given by $r_1=\{(\varepsilon,\varepsilon),(a,a),(b,b)\}$, and two unary operations $sub_{\rm l}$ and $sub_{\rm r}$ defined by $sub_{\rm l}(wx)=x$ and $sub_{\rm r}(wx)=w$ for all $w\in \{a,b\}^*$ and $x\in\{a,b\}$, $sub_{\rm l}(\varepsilon)=\varepsilon$, and $sub_{\rm r}(\varepsilon)=\varepsilon$. The identity relation ${\rm id}_{U_{{\cal A}_1}}$ given by ${\rm id}_{U_{{\cal A}_1}}={\rm id} \cap (\{a,b\}^*\times \{a,b\}^*)$ is the smallest binary relation $r$ with the following properties. (1) $r$ contains $(\varepsilon, \varepsilon)$. (2) If, for any two strings $ w_1, w_2\in \{a,b\}^*$, $( w_1, w_2)\in r$ holds, then $( w_1x, w_2x)\in r$ holds for any $x\in \{a,b\}$. The identity relation does not belong to ${\cal A}_1$, but ${\rm id}_{U_{{\cal A}_1}}$ is semi-decidable by a BSS RAM ${\cal M}_1\in {\sf M}_{{\cal A}_1}$ using ${\cal P}_{{\cal M}_1}$ and the complement $U_{{\cal A}_1}^\infty\setminus {\rm id}_{U_{{\cal A}_1}}$ is semi-decidable by a BSS RAM ${\cal M}_2\in {\sf M}_{{\cal A}_1}$ using ${\cal P}_{{\cal M}_2}$.

\vspace{0.2cm}

\noindent {\sf 
\begin{tabular}{lll}
${\cal P}_{{\cal M}_1}$: &\quad$\ell_1: \, \, I_2:=I_2+1;$ \,\, $\ell_1+1: $ if $I_1=I_2$ then goto $\ell_2$ else goto $\ell_1+1;$\\
&\quad$\ell_2: \, \, Z_5:= c_3^0;$ \\ 
&\textcolor{red}{\quad $\ell_3 : \,\,Z_3:=f_1^1(Z_1);$ \,\, $ Z_4:=f_1^1(Z_2); $} \\
&\textcolor{red}{\quad $\ell_4 : \,\,Z_1:=f_2^1(Z_1);$ \,\, $ Z_2:=f_2^1(Z_2); $} \\
&\quad$\ell_5 : \,\,$ if $r_1^{2}(Z_3,Z_4)$ then goto \textcolor{blue}{$\ell_6$} else goto $\textcolor{blue}{\ell_5};$ \\
&\quad$\ell_6 : \,\,$ if $r_1^{2}(Z_4,Z_5)$ then goto $\textcolor{blue}{\ell_7}$ else goto $\textcolor{blue}{\ell_3};$ \\
&\quad$\ell_7: \, \,$ stop. \\
\end{tabular}}

\vspace{0.1cm}

\noindent {\sf 
\begin{tabular}{lll}
${\cal P}_{{\cal M}_2}$:&\quad$\ell_1: \, \, I_2:=I_2+1;$ \,\, if $I_1=I_2$ then goto $\ell_2$ else goto $\ell_7;$\\
&\quad$\ell_2: \, \, Z_5:= c_3^0;$ \\ 
&\textcolor{red}{\quad $\ell_3 : \,\,Z_3:=f_1^1(Z_1);$ \,\, $ Z_4:=f_1^1(Z_2); $} \\
&\textcolor{red}{\quad $\ell_4 : \,\,Z_1:=f_2^1(Z_1);$ \,\, $ Z_2:=f_2^1(Z_2); $} \\
&\quad$\ell_5 : \,\,$ if $r_1^{2}(Z_3,Z_4)$ then goto $\textcolor{blue}{\ell_6}$ else goto $\textcolor{blue}{\ell_7};$ \\
&\quad$\ell_6 : \,\,$ if $r_1^{2}(Z_4,Z_5)$ then goto $\textcolor{blue}{\ell_6}$ else goto $\textcolor{blue}{\ell_3};$ \\
&\quad$\ell_7: \, \,$ stop. \\
\end{tabular}}

\vspace{0.2cm}

\noindent
The \textcolor{red}{red} program segment is intended to decompose the inputted strings. The subprograms labeled by $\ell_5$ and $\ell_6$, respectively, enable the evaluation of the symbols $c(Z_3)$ and $c(Z_4)$ resulting from splitting $c(Z_1)$ and $c(Z_2)$ and cutting off the last symbols on the right-hand side. If both $c(Z_3)$ and $c(Z_4)$ are $\varepsilon$, then ${\cal M}_1$ halts and accepts the input. $(baa,aaa)$ is not in ${\rm id}_{U_{{\cal A}_1}}$. For the input $(baa,aaa)$, ${\cal M}_1$ traverses the infinite computation path $(\ell_1,\ell_1+1,\ell_2,\ell_3,\ell_3+1,\ell_4,\ell_4+1,\ell_5,\ell_6, \ell_3,\linebreak\ell_3+1,\ell_4,\ell_4+1,\ell_5,\ell_6, \ell_3,\ell_3+1,\ell_4,\ell_4+1,\ell_5,\textcolor{blue}{\ell_5},\textcolor{blue}{\ell_5},\textcolor{blue}{\ldots})$ determined by  the sequence in ${\sf Con}_ {{\cal M}_1, (baa,aaa)}$ which contains the following subsequence of configurations.

\vspace{0.1cm}

\begin{tabular}{l}
\hspace{0.96cm}$(\ell_1\,.\, (2,1) \,.\, (baa,aaa,aaa,aaa,aaa,aaa,aaa,\ldots))$\\
$\to_ {{\cal M}_1}^* (\ell_2\,.\, (2,2) \,.\, (baa,aaa,aaa,aaa,aaa,aaa,aaa,\ldots))$\\
$\to _{{\cal M}_1}^* (\ell_5\,.\, (2,2) \,.\, (ba,aa,a,a,\varepsilon,aaa,aaa,\ldots))$\\
$\to _{{\cal M}_1}^* (\ell_5\,.\, (2,2) \,.\, (b,a,a,a,\varepsilon,aaa,aaa,\ldots))$\\
$\to _{{\cal M}_1} (\ell_6\,.\, (2,2) \,.\, (b,a,a,a,\varepsilon,aaa,aaa,\ldots))$\\
$\to _{{\cal M}_1} (\ell_3\,.\, (2,2) \,.\, (b,a,a,a,\varepsilon,aaa,aaa,\ldots))$\\
\end{tabular}

$\to _{{\cal M}_1}^* (\textcolor{blue}{\ell_5}\,.\, (2,2) \,.\, (\varepsilon,\varepsilon,b,a,\varepsilon,\ldots)) \to _{{\cal M}_1}\! (\textcolor{blue}{\ell_5}\,.\, (2,2) \,.\, (\varepsilon,\varepsilon,b,a,\ldots))\! \to _{{\cal M}_1}\textcolor{blue}{\cdots}$

\vspace{0.1cm}

\noindent
 $\to_ {{\cal M}_1}^*$ is the transitive extension of $\to_ {{\cal M}_1}$. For input $(ab,ab)$, ${\cal M}_1$ traverses the finite computation path $(\ell_1,\ell_1+1,\ell_2,\ell_3,\ell_3+1,\ell_4,\ell_4+1,\ell_5,\ell_6, \ell_3, \ell_3+1,\ell_4,\ell_4+1,\linebreak\ell_5,\ell_6, \ell_3,\ell_3+1,\ell_4,\ell_4+1,\ell_5,\ell_6,\textcolor{blue}{\ell_7})$ and thus $(ab,ab)$ can be accepted.

\vspace{0.15cm}

\begin{tabular}{l}
\hspace{0.96cm}$(\ell_1\,.\, (2,1) \,.\, (ab,ab,ab,ab,ab,ab,ab,\ldots))$\\
$\to_ {{\cal M}_1}^* (\ell_2\,.\, (2,2) \,.\, (ab,ab,ab,ab,ab,ab,ab,\ldots))$\\
$\to _{{\cal M}_1}^* (\ell_5\,.\, (2,2) \,.\, (a,a,b,b,\varepsilon,ab,ab,\ldots))$\\
$\to _{{\cal M}_1}^* (\ell_5\,.\, (2,2) \,.\, (\varepsilon,\varepsilon,a,a,\varepsilon,ab,ab,\ldots))$\\
$\to _{{\cal M}_1}^* (\ell_5\,.\, (2,2) \,.\, (\varepsilon,\varepsilon,\varepsilon,\varepsilon,\varepsilon,ab,ab,\ldots))$\\
$\to _{{\cal M}_1}^* (\textcolor{blue}{\ell_7}\,.\, (2,2) \,.\, (\varepsilon,\varepsilon,\varepsilon,\varepsilon,\varepsilon,ab,ab,\ldots))$\\
\end{tabular}

\vspace{0.15cm}

\noindent For input $(\varepsilon, \varepsilon)$, ${\cal M}_2$ traverses the computation path
$(\ell_1,\ell_1+1,\ell_2,\ell_3,\ell_3+1,\ell_4,\linebreak\ell_4+1,\ell_5,\ell_6, \textcolor{blue}{\ell_6}, \textcolor{blue}{\ell_6}, \textcolor{blue}{\ldots})$. 

\noindent Moreover, $\chi_{{\rm id}_{U_{{\cal A}_1}}}$ can be computed, for example, by using the following program.

\vspace{0.1cm}

{\sf 
\begin{tabular}{lll}
${\cal P}_{\cal M}$:&\quad $\ell_1: \, \, I_2:=I_2+1;$ \,\, if $I_1=I_2$ then goto $\ell_1'$ else goto $\ell_1'';$\\
&\quad $\ell_1': \, \, I_2:=I_2+1;$ \,\, $\textcolor{magenta}{I_1:=1;}$ \,\, $Z_{I_2}:=Z_{I_1};$ \,\,\textcolor{magenta}{$ Z_1:=c_1^0;$} \,\, goto $\ell_2;$ \\ 
&\quad $\ell_1'': \, \, \textcolor{magenta}{I_1:=1;}$ \,\, $\textcolor{magenta}{Z_1:=c_2^0;}$ \,\, goto $\ell_7;$ \\ 
&\quad $\ell_2: \, \, Z_6:= c_3^0;$ \\ 
&\quad \textcolor{red}{$\ell_3 : \,\,Z_4:=f_1^1(Z_2);$ \,\, $ Z_5:=f_1^1(Z_3); $} \\
&\quad \textcolor{red}{$\ell_4 : \,\,Z_2:=f_2^1(Z_2);$ \,\, $ Z_3:=f_2^1(Z_3); $} \\
&\quad $\ell_5 : \,\,$ if $r_1^{2}(Z_4,Z_5)$ then goto $\textcolor{blue}{\ell_6}$ else goto $\textcolor{blue}{\ell_1''};$ \\
&\quad $\ell_6 : \,\,$ if $r_1^{2}(Z_5,Z_6)$ then goto $\textcolor{blue}{\ell_7}$ else goto $\textcolor{blue}{\ell_3};$ \\
&\quad$\ell_7: \, \,$ stop. \\
\end{tabular}}

\vspace{0.1cm}\noindent If we replace the program segment \fbox{{\sf $\ell_7: \, \,$ stop.}} in ${\cal P}_{\cal M}$ --- similarly as proposed in the proof of Proposition \ref{Propos1} --- by 

\hspace{1.7cm }{\sf $\ell_7: \, \, Z_2:=c_1^0;$ if $r_1^{2}(Z_1,Z_2)$ then goto $\ell_8$ else goto $\ell_7;$ 

\hspace{1.7cm}$\ell_8:\,\,$ stop.}

\noindent then the resulting program also allows to semi-decide $\,{\rm id}_{U_{{\cal A}_1}}$ over ${\cal A}_1$.

\vspace{0.2cm}\noindent If we replace the program segment \fbox{{\sf $\ell_7: \, \,$ stop.}} in ${\cal P}_{\cal M}$ by 

\hspace{1.7cm }{\sf $\ell_7: \, \, Z_2:=c_2^0;$ if $r_1^{2}(Z_1,Z_2)$ then goto $\ell_8$ else goto $\ell_7;$ 

\hspace{1.7cm}$\ell_8:\,\,$ stop.}

\noindent then the resulting program allows to semi-decide the complement of $\,{\rm id}_{U_{{\cal A}_1}}$.
\end{example}

\begin{example}[A structure with $(d)$ in which $(c)$ does not hold]\label{CounterEx2}\hfill 

\noindent Let ${\cal A}_2$ be the finite structure $(\{a,b,c,d\} ;a,b ;; r_a,r_b)$ of signature $(2;;1,1)$ in ${\sf Struc}_{a,b}$ that contains four elements, the constants $a$ and $b$, and two unary relations $r_a$ and $r_b$ given by $r_a=\{a\}$ and $ r_b= \{b\}$.
 The semi-decidability of $\{a\}$ over ${\cal A}_2$ (i.e., the semi-decidability of $\{a\}$ by a machine in ${\sf M}_{{\cal A}_2}$) results from the execution of a program such as

\vspace{0.1cm}

{\sf 
\begin{tabular}{ll}
\quad\qquad $1: \,\,\,$ if $I_1=I_2$ then goto $2$ else goto $1;$\\
\quad\qquad $2 : \,\,\,$ if $r_1^{1}(Z_1)$ then goto $3$ else goto $2;$ \\
\quad\qquad$3:\,\,\,$ stop. \\
\end{tabular}}

\vspace{0.1cm}

\noindent If the literal \fbox{$r_1^{1}(Z_1)$} is replaced by $r_2^{1}(Z_1)$, then the new program can be used for semi-deciding $\{b\}$.
The semi-decidability of the complement $U_{{\cal A}_2}^\infty\setminus\{a\}$ results from changing some of the labels in the {\sf goto}-commands.

\vspace{0.1cm}

{\sf 
\begin{tabular}{ll}
\quad\qquad $1: \,\,\,$ if $I_1=I_2$ then goto $2$ else goto $3;$\\
\quad\qquad $2 : \,\,\,$ if $r_1^{1}(Z_1)$ then goto $2$ else goto $3;$ \\
\quad\qquad$3:\,\,\,$ stop. \\
\end{tabular}}

\vspace{0.1cm}

\noindent If the literal \fbox{$r_1^{1}(Z_1)$} is replaced by $r_2^{1}(Z_1)$, then the new program can be used for semi-deciding $U_{{\cal A}_2}^\infty\setminus\{b\}$.
 Moreover, $\chi_{\{a\}}$ can be computed by the BSS RAM in $ {\sf M}_{{\cal A}_2}$ using the following program.

{\sf 
\begin{tabular}{ll}
\quad\qquad $\ell_1: \,\,$ if $I_1=I_2$ then goto $\ell_2$ else goto $\ell_3;$\\
\quad\qquad $\ell_2 : \,\,$ if $r_1^{1}(Z_1)$ then goto $\ell_4$ else goto $\ell_3;$ \\
\quad\qquad$\ell_3: \, \, I_1:=1; \, \, Z_1:=c_2^0;$ \\ 
\quad\qquad$\ell_4: \,\,$ stop. \\
\end{tabular}}

\vspace{0.1cm}

\noindent $\chi_{\{b\}}$ can be computed by the BSS RAM in $ {\sf M}_{{\cal A}_2}$ using the following program.

\vspace{0.1cm}

{\sf 
\begin{tabular}{ll}
\quad\qquad $\ell_1: \, \, I_2:=I_2+1;$ \,\, $Z_{I_2}:=Z_{I_1};$ \,\, $Z_1:=c_2^0;$ \\ 
\quad\qquad $\ell_2: \,\,$ if $I_1=I_3$ then goto $\ell_3$ else goto $\ell_5;$\\
\quad\qquad $\ell_3 : \,\,$ if $r_2^{1}(Z_2)$ then goto $\ell_4$ else goto $\ell_5;$ \\
\quad\qquad$\ell_4: \, \, Z_1:=c_1^0;$ \\ 
\quad\qquad$\ell_5: \, \,I_1:=1;$ stop. \\
\end{tabular}}

\vspace{0.4cm}

\noindent {\bf Why is $(c)$ not true in ${\cal A}_2$?} Let us assume that there is an ${\cal M}\in {\sf M}_{{\cal A}_2}$ that does not halt on input $(d,c)\not\in {\rm id}_{U_{{\cal A}_2}}$ and that halts on input $(c,c)\in {\rm id}_{U_{{\cal A}_2}}$. For these inputs, the initial configurations $con_{1}$ and $con'_{1}$ of ${\cal M}$ are given by 

\vspace{0.1cm}

\qquad $con_{1} =(1\,.\,(2,1,\ldots,1)\,.\,(d,c,c,\ldots))$, 

\qquad $con'_{1} =(1\,.\,(2,1,\ldots,1)\,.\,(c,c,c,\ldots))$.

\vspace{0.2cm}

\noindent For both resulting sequences  $(con_{t})_{t\geq 1}\in {\sf Con}_{{\cal M},(d,c)}$ and $(con'_{t})_{t\geq 1}\in {\sf Con}_{{\cal M},(c,c)}$ given by 

\vspace{0.1cm}

\qquad $con_{t}= (\ell_t\,.\,\vec\nu^{(t)}\,.\,\bar u^{(t)})$ \hfill \textcolor{gray}{\small ($\ell_t\in {\cal L}_{\cal M}$, $\vec\nu^{(t)}\in \bbbn_+^{k_{\cal M}}$, $\bar u^{(t)}\in U_{{\cal A}_2}^{\omega}$)}, 

\qquad $con'_{t}= (\ell'_t\,.\,\vec\mu^{(t)}\,.\,\bar z^{(t)})$ \hfill \textcolor{gray}{\small($\ell_t'\in {\cal L}_{\cal M}$, $\vec\mu^{(t)}\in \bbbn_+^{k_{\cal M}}$, $\bar z^{(t)}\in U_{{\cal A}_2}^{\omega}$)}, 

\vspace{0.2cm}

\noindent and $\bar u^{(t)}=( u_{t,1}, u_{t,2},\ldots)$ and $\bar z^{(t)}=( z_{t,1}, z_{t,2},\ldots)$ for all $t\geq 1$, we can prove the following properties $(i)$ to $(iii)$ by induction on $t$. 

\vspace{0.1cm}

 \hspace*{0.36cm}$\begin{array}{lclllcllc} 
 \ell'_{t}&=&\ell_{t}, \hspace*{1.2cm}&&\vec\mu^{(t)}&=&\vec\nu^{(t)}, &&(i)\\
 u_{t,j}&\in& \{a,b,c\} &\Rightarrow \hspace*{0.6cm} & z_{t,j}&=&u_{t,j}
 &(\mbox{\rm for all }j\geq 1),&(ii)\\
 u_{t,j}&=&d & \Rightarrow & z_{t,j}&=&c &(\mbox{\rm for all }j\geq 1).\hspace*{0.4cm}&(iii)\\
\end{array}$

\vspace{0.2cm}

\noindent Note, that the symbol $\Rightarrow$ is here used to denote an implication in our metatheory. For the initial configurations, these relationships hold if we assume that $t=1$ is valid in this case. Therefore, the assumption holds for the base case $t=1$. If, for a $t\geq 1$, the relationships $(i)$ to $(iii)$ and, thus,

\vspace{0.2cm}

 \noindent\quad$\!z_{t,j}\in r_a$ \, $\Leftrightarrow$ \, $u_{t,j}\in r_a$ \,\,\,\,{\em and} \,\, $z_{t,j}\in r_b$ \, $\Leftrightarrow$ \, $u_{t,j}\in r_b$ 

\vspace{0.2cm}

\noindent hold (where $\Leftrightarrow$ refers to an equivalence), then a further step of the execution of ${\cal M}$ leads to the relationships 

\vspace{0.1cm}

$\begin{array}{lclllcll} 
 \ell'_{t+1}&=&\ell_{t+1}, \hspace*{0.86cm} &&\vec\mu^{( t+1)}&=&\vec\nu^{( t+1)}, \\
 u_{t+1,j}&\in& \{a,b,c\} &\Rightarrow \hspace*{0.6cm} & z_{t+1,j}&=&u_{t+1,j}
 &(\mbox{\rm for all }j\geq 1),\\
 u_{t+1,j}&=&d & \Rightarrow & z_{t+1,j}&=&c &(\mbox{\rm for all }j\geq 1).\\
\end{array}$
\vspace{0.2cm}

 \noindent This means by induction that there holds $\ell_{t}=\ell'_{t}$ for any $t\geq 1$ which implies that both resulting computation paths are finite or both paths are infinite. Because this contradicts our assumption, ${\rm id}_{U_{{\cal A}_2}}$ is not in ${\rm SDEC}_{{\cal A}_2}$.
\end{example}

\begin{example}[A structure with $(e)$ in which $(d)$ does not hold]\label{CounterEx3}\hfill 

\noindent Let ${\cal A}_3$ be the finite structure $(\{a,b,c,d\} ;c,d ;; r_a,r_b)$ of signature $(2;\!;1,1)$ in ${\sf Struc}_{c,d}$ that contains four elements, the constants $c$ and $d$, and two unary relations $r_a$ and $r_b$ given by $r_a=\{a\}$ and $ r_b= \{b\}$. $\{c,d\}$ is semi-decidable by a BSS RAM ${\cal M}_{c,d}\in {\sf M}_{{\cal A}_3}$ using the following program ${\cal P}_{{\cal M}_{c,d}}$.

\vspace{0.1cm}

{\sf 
\begin{tabular}{ll}
\qquad\quad $\ell_1: \, \, I_2:=I_2+1;$ \, $\ell_1+1:$\, if $I_1=I_2$ then goto $\ell_2$ else goto $\ell_1+1;$\\
\qquad\quad $\ell_2 : \,\,$ if $r_1^{1}(Z_2)$ then goto $\ell_2$ else goto $\ell_3;$ \\
\qquad\quad $\ell_3 : \,\,$ if $r_2^{1}(Z_2)$ then goto $\ell_3$ else goto $\ell_4;$ \\
\qquad\quad$\ell_4: \, \, $ stop. \\
\end{tabular}}

\vspace{0.1cm}

\noindent Note, that $U_{{\cal A}_3}^\infty\setminus \{c,d\}$ is also semi-decidable by a BSS RAM. ${\cal M}_{a,b}^*\in {\sf M}_{{\cal A}_3}$ semi-decides $\{a,b\}\cup (U_{{\cal A}_3}^\infty\setminus U_{{\cal A}_3})$ by using the following program ${\cal P}_{{\cal M}_{a,b}^*}$.

\vspace{0.1cm}

{\sf 
\begin{tabular}{ll}
\qquad\quad $\ell_1: \, \, I_2:=I_2+1;$ \, \, if $I_1=I_2$ then goto $\ell_2$ else goto $\ell_4;$\\
\qquad\quad $\ell_2 : \,\,$ if $r_1^{1}(Z_2)$ then goto $\ell_4$ else goto $\ell_3;$ \\
\qquad\quad $\ell_3 : \,\,$ if $r_2^{1}(Z_2)$ then goto $\ell_4$ else goto $\ell_3;$ \\
\qquad\quad$\ell_4: \, \, $ stop. \\
\end{tabular}}

\vspace{0.4cm}

\noindent {\bf Why is $(d)$ not true in ${\cal A}_3$?} Let us assume that there is a machine ${\cal M}\in {\sf M}_{{\cal A}_3}$ that halts for input $c\in \{c\}$ and that does not halt for the input $d\not\in \{c\}$. For the constant $c$ and the constant $d$ as inputs, the corresponding initial configurations $con_{1}$ and $con'_{1}$ are given by $con_{1}=(1\,.\,(1,1,\ldots,1)\,.\,(c,c,c,\ldots))$ and $con'_{1}=(1\,.\,(1,1,\ldots,1)\,.\,(d,c,c,\ldots))$, respectively. For these initial configurations and both resulting sequences  $(con_t)_{t\geq 1} \in {\sf Con}_{{\cal M},c}$ and  $(con_t')_{t\geq 1} \in  {\sf Con}_{{\cal M},d}$  given, for $t\geq 2$, by the configurations $con_{t}$ and $con'_{t}$ with 
 $con_{t}= (\ell_{t}\,.\,\vec\nu^{(t)}\,.\,\bar u^{(t)})$ and $con'_{t}= (\ell'_{t}\,.\,\vec\mu^{(t)}\,.\,\bar z^{(t)})$, we have $\ell_t=\ell_t'$, $\vec\nu^{(t)}=\vec\mu^{(t)}$, and

\vspace{0.1cm}

 \hspace*{0.53cm}$u_{t,j}\not \in r_a$ \, {\rm and} \, $z_{t,j}\not \in r_a$ \, {\rm and} \,
$u_{t,j}\not \in r_b$ \, {\rm and} \, $z_{t,j}\not \in r_b$ 

\vspace{0.1cm}

\noindent for all $j\geq 1$ and all $t\geq 1$ which implies, among others, $\ell_{t+1}=\ell'_{t+1}$ for any $t\geq 1$ and, consequently, the result is again that both corresponding computation paths end with $\ell_{{\cal P}_{\cal M}} $ or both paths are infinite. Because this is a contradiction to our assumption, the set $\{c\}$ is not in ${\rm SDEC}_{{\cal A}_3}$.
\end{example}
\subsubsection{The decidability of semi-decidable single-element sets}\label{FromSemiToDecSecI} 
The structure ${\cal A}_2$ (in Example \ref{CounterEx2}) was in particular considered to give an example for single-element sets containing one constant and being semi-decidable. However, it is also easy to see that the sets of this kind are even decidable. The next proposition states that this is no coincidence and holds for any single-element set. Statements such as the following can be proved by using basic knowledge of classical first-order logic. For more theoretical details, see our {\sf special part on program and computation paths}.

\begin{proposition}[Recognizing an individual implies decidability]\label{SemiDeci}\hfill 

\noindent Let ${\cal A}$ be a structure in $ {\sf Struc}$ and $x_0\in U_{\cal A}$.

 The semi-decidability of $\{x_0\}$ by a BSS RAM in ${\sf M}_{\cal A}$ implies the semi-decidability of the complement of $\{x_0\}$  over ${\cal A}$ and thus  $\{x_0\}\in {\rm DEC}_{\cal A}$. 
\end{proposition}
{\bf Proof.} Let ${\cal M}$ be a BSS RAM in ${\sf M}_{\cal A}$ that semi-decides $\{x_0\}$ which means that $H_{\cal M}=\{x_0\}$ holds. Then, ${\cal M}$ goes, for input $x_0\in U_{\cal A}$ and only for $x_0$, through a finite computation path $(\ell_{x_0,t})_{t=1..s} $ which means that $\ell_{x_0,t}\in \{1,\ldots, \ell_{{\cal P}_{\cal M}}\}$ and $\ell_{x_0,s}=\ell_{{\cal P}_{\cal M}}$ hold. Let this path be determined by the sequence $(con_{x_0,t})_{t\geq 1} \in {\sf Con}_{{\cal M},x_0}$ given by  $con_{x_0,t}= (\ell_{x_0,t}\,.\,\vec\nu^{(x_0,t)}\,.\,\bar u^{(x_0,t)})$ and $\vec\nu^{(x_0,t)} \in \bbbn_+^{k_{\cal M}}$ for all $t\geq 1$. Consequently, there holds $con_{x_0,1}=(1\,.\,(1,\ldots,1)\,.\,(x_0,x_0,\ldots))$ where $(1,\ldots,1)$ stands for $\vec\nu^{(x_0,1)} \in \{1\}^{k_{\cal M}}$. For all other inputs $\vec x\in U_{\cal A}^\infty$ (which means $\vec x\not =x_0$), ${\cal M}$ does not halt and $\vec x$ thus traverses an infinite computation path $(\ell_{\vec x,t})_{t\geq 1} $ determined by $(con_{\vec x,t})_{t\geq 1} \in {\sf Con}_{{\cal M},\vec x}$ with $con_{\vec x,t}= (\ell_{\vec x,t}\,.\,\vec\nu^{(\vec x,t)}\,.\,\bar u^{(\vec x,t)})$ and $\ell_{x_0,t}\not=\ell_{{\cal P}_{\cal M}}$ for all $t\geq 1$. This means that there is a $t$ with $\ell_{\vec x,t}\not = \ell_{x_0,t}$. Let $t_{\vec x}$ be the smallest integer $t$ of this kind. Consequently, $\ell_{\vec x,t_{\vec x}-1} = \ell_{x_0,t_{\vec x}-1}$ holds and $\ell_{x_0,t_{\vec x}-1}$ is the label of an instruction in ${\cal P}_{\cal M}$ that is of type (4) or (5). This implies that the complement $U_{\cal A}^\infty \setminus\{x_0\}$ can be semi-decided by using a BSS RAM ${\cal M}^{\rm c}_{x_0}\in {\sf M}_{\cal A}$ whose program is given in Overview \ref{SemiOver1} and includes the program segment ${\cal P}_{(s+1)}$ defined in Overview \ref{SemiOver0}. For the use of $ {\cal L}_{{\cal M},{\rm T}}, {\cal L}_{{\cal M},{\rm H}_{\rm T}}, \ldots$ consisting of all the labels for instructions of type $(4), (5), \ldots$ , see \cite{GASS20} and \cite{GASS25A}.
\begin{overview}
[From $\{x_0\}\in {\rm SDEC}_{\cal A}$ to segments ${\cal P}_{(s+1)}$ and ${\cal P}_{(s+2)}$]\label{SemiOver0}

\hfill

\nopagebreak 

\noindent \fbox{\parbox{11.8cm}{

Let the program ${\cal P}_{\cal M}$ of ${\cal M}$ semi-deciding $\{x_0\}$ be of the form \textcolor{blue}{$(*)$} and let

\vspace{0.1cm}

\quad $\ell_{x_0,1}:$ {\sf instruction}$_{\ell_{x_0,1}}; \,\ldots;\quad \ell_{x_0,s-1}:$ {\sf instruction}$_{\ell_{x_0,s-1}}; \,\,\textcolor{red}{\ell_{x_0,s}}\,: {\sf stop}.$ \,\,$\textcolor{blue}{(*_{x_0})}$

\vspace{0.1cm}

be the sequence of labeled instruction resulting from 
 $(\ell_{x_0,t}) _{t=1.. s}$ and 
 ${\cal P}_{\cal M}$. 

\vspace{0.1cm}

Let, for $i \in \{1,2\}$, ${\cal P}_{(s+i)}$ be the program segment

\vspace{0.1cm}

\hspace{1.19cm} $1:$ {\sf instruction}$_{1}'; \quad \ldots;\quad s-1:$ {\sf instruction}$_{s-1}';$

\vspace{0.1cm}

derived from $\textcolor{blue}{(*_{x_0})}$ such that each of the new labeled instructions

\vspace{0.1cm}

\hspace{1.19cm} $t:$ {\sf instruction}$_{t}'$ \hfill (for $t\in \{ 1,\ldots, s-1\}$)\, 

\vspace{0.05cm}

results

\vspace{0.1cm}

\begin{itemize}\parskip -0.5mm
\item from \fbox{$\ell_{x_0,t} :{\sf instruction}_{\ell_{x_0,t}}$} by replacing only

\vspace{0.1cm}
\begin{itemize}
\item the label \fbox{$\ell_{x_0,t}$} by $t$
\hfill \textcolor{gray}{\small ({\sf instruction}$_{t}' $ and {\sf instruction}$_{\ell_{x_0,t}}$ match)}
\end{itemize}
{\sf if} ${\ell_{x_0,t}}\not \in {\cal L}_{{\cal M},{\rm T}}\cup {\cal L}_{{\cal M},{\rm H}_{\rm T}}$ ,
\end{itemize}

or

\begin{itemize}\parskip -0.5mm
\item from \fbox{\sf $\ell_{x_0,t} : \,$ if $cond_{\ell_{x_0,t}}$ then goto $\ell_{x_0,t+1} $ else goto \textcolor{blue}{$\ell^{(\ell_{x_0,t})}_2$} }
or

\vspace{0.1cm}

\hspace{0.85cm}\fbox{\sf $\ell_{x_0,t} : \,$ if $cond_{\ell_{x_0,t}}$ then goto \textcolor{blue}{$\ell^{(\ell_{x_0,t})}_1$} else goto $\ell_{x_0,t+1}$} by replacing

\vspace{0.1cm}

\begin{itemize}
\parskip -0.5mm 
\item the label \fbox{$\ell_{x_0,t}$} by $t$ ,
\item the label \fbox{$\ell_{x_0,t+1}$} by $t+1$ ,
\item the label \fbox{$\textcolor{blue}{\ell^{(\ell_{x_0,t})}_j}$} by \textcolor{blue}{$ s+i$}\hfill {\small (for $j=1$ and $j=2$, resp.)},
\end{itemize}
{\sf otherwise}. 
\end{itemize}
}}
\end{overview}

For each instruction $instruction_{\ell_{x_0,t}}$ of type (4) or (5) whose label $\ell_{x_0,t}$ belongs to the computation path $(\ell_{x_0,t}) _{t=1.. s}$ traversed by ${\cal M}$, the condition $cond_{\ell_{x_0,t}}$ has the form $r_{i_\ell}^{k_{i_\ell}}(Z_{j_{\ell,1}},\ldots, Z_{j_{\ell,k_{i_\ell}}})$ or $I_{j_\ell}=I_{k_\ell}$ for certain indices $i_\ell, j_\ell,k_\ell,j_{\ell,1},\ldots$ that depend on the label $\ell\in {\cal L}_{{\cal M},{\rm T}}\cup {\cal L}_{{\cal M},{\rm H}_{\rm T}}$ that satisfies $\ell_{x_0,t}=\ell$. 

\begin{overview}[From $\{x_0\}\in {\rm SDEC}_{\cal A}$ to the co-semi-decision of $\{x_0\}$]\label{SemiOver1}

\hfill

\nopagebreak 

\noindent \fbox{\parbox{11.8cm}{

Let the program of the machine ${\cal M}_{x_0}^{\rm c}\in {\sf M}_{\cal A}$ be given by

\vspace{0.1cm}

\hspace{1.39cm} ${\cal P}_{ \textcolor{blue}{(s+1)}}$ {\sf\,\, $\textcolor{red}{s: {\sf goto } \,s};$ \,\, \textcolor{blue}{$s+1: $ stop}}.  .
}}
\end{overview}

\vspace{0.1cm}

{\bf Why does ${\cal M}_{x_0}^{\rm c}$ semi-decide the complement of $\{x_0\}$?} We have to show that, for any input $\vec x\in U_{\cal A}^\infty\setminus \{x_0\}$, the execution of ${\cal M}_{x_0}^{\rm c}$ on $\vec x$ does not lead to an infinite computation path $(\ell'_t)_{t\geq 1}$ with $\ell_t'\in \{1,\ldots,s\}\setminus\{s+1\}$. For $t_{\vec x}$ defined above, let us consider the first $t_{\vec x}-2$ labels $\ell_{x_0,t}$ (with $t\leq t_{\vec x}-2$) that belong to both computation paths traverses by ${\cal M}$ for $x_0$ and for $\vec x$. If, for $t\leq t_{\vec x}-2$, the label $\ell_{x_0,t}$ belongs to $ {\cal L}_{{\cal M},{\rm T}}\cup {\cal L}_{{\cal M},{\rm H}_{\rm T}}$, then we can assume, by definition of $t_{\vec x}$, that the evaluation of the conditions $cond_{\ell_{x_0,t}}$ (whose values depend on the configurations $con_{\vec x,t}$ and $con_{ x_0,t}$, respectively) leads to the same truth values and to the same next label in the initial part $\ell_{x_0,1}, \ldots, \ell_{x_0,t_{\vec x}-1}$ of both computation paths traversed by ${\cal M}$ for $x_0$ and $\vec x$. Moreover, we have ${\ell_{x_0,t_{\vec x}-1}\in {\cal L}_{{\cal M},{\rm T}}\cup {\cal L}_{{\cal M},{\rm H}_{\rm T}}}$ and the evaluation of the condition $cond_{\ell_{x_0,t_{\vec x}-1}}$ provides different truth values. Because the computed contents of the index and $Z$-registers of ${\cal M}$ and ${\cal M}_{x_0}^{\rm c}$ match for input $\vec x$ until this point and the same holds for both machines on $x_0$, the evaluation of $cond_{\ell_{x_0,t_{\vec x}-1}}$ provides also different truth values for ${\cal M}_{x_0}^{\rm c}$ on $x_0$ and $\vec x$. This means that $(\ell_{x_0,1},\ldots, \ell_{x_0,t_{\vec x}-1},\textcolor{blue}{\ell_i^{( \ell_{x_0,t_{\vec x}-1}) }}, \ldots)$ (with $i=1$ or $i=2$) is the first part of the computation path traversed by ${\cal M}$ for the input $\vec x$ (for $\vec x\not = x_0$) and thus the path $(1,\ldots, t_{\vec x}-1,\textcolor{blue}{s+1})$ of length $t_{\vec x}$ is the computation path traversed by ${\cal M}_{x_0}^{\rm c}$ for $\vec x$ which implies that ${\cal M}_{x_0}^{\rm c}$ halts for $\vec x$. Since the computation path traversed by ${\cal M}$ for $ x_0$ is $(\ell_{x_0,1},\ldots, \ell_{x_0,s-1}, \textcolor{red}{\ell_{{\cal P}_{\cal M}}})$, ${\cal M}_{x_0}^{\rm c}$ traverses $(1,\ldots, s-1, \textcolor{red}{s},\textcolor{red}{s},\textcolor{red}{\ldots})$ for $ x_0$. Thus, ${\cal M}_{x_0}^{\rm c}$ loops forever only for the input $ x_0$. 
\qed

\vspace{0.3cm}

 The construction idea used in the proof of Proposition \ref{SemiDeci} is (even) more interesting when we know a suitable machine semi-deciding $\{x_0\}$.  For the finite signature $\sigma$ given in the introduction, it is easy to encode the $\sigma$-programs in such a way that all these codes can be effectively enumerated and the resulting set ${\rm CODES}_{\sigma}^{\rm Goedel}$ of all codes is decidable (for suitable G\"odel numbers storable in index registers, see \cite{GASS20}). 
 Let ${\cal A}$ be a structure of signature $\sigma$ and contain at least two individuals, a constant $c_i$, and ${\rm id}_{U_{\cal A}}$. Let us consider the single-element set $\{c_i\}$ and let ${\rm HALT}_{\cal A}^{\{c_i\}}$, ${\rm HALT}_{\cal A}^{\{c_i\}\subseteq}$, and ${\rm HALT}_{\cal A}^{\subseteq\{c_i\}}$ be the subsets of $ \mbox{{\rm CODES}}_{\sigma}^{\rm Goedel}$ containing the codes of all BSS RAMs ${\cal M}$ in ${\sf M}_{\cal A}$ whose halting sets $H_{\cal M}$  satisfy $H_{\cal M}= \{c_i\}$, $\{c_i\} \subseteq H_{\cal M}$, and $H_{\cal M} \subseteq \{c_i\}$, respectively. Then, ${\rm HALT}_{\cal A}^{\{c_i\}\subseteq}$ belongs to the class ${\cal A}\mbox{-}{\textcolor{blue}{\Sigma_{1 }^{\rm ND}}}$ of all problems non-deterministically semi-decidable by a BSS RAM in ${\sf M}^{\rm ND}_{\cal A}$ (for the definition of such classes, see \cite[pp.\,30--31]{GASS25A}). For showing ${\rm HALT}_{\cal A}^{\{c_i\}\subseteq}\in {\cal A}\mbox{-}{\textcolor{blue}{\Sigma_{1 }^{\rm ND}}}$, the program of a universal machine (given in \cite{GASS20}) can be used as a program segment and applied to an inputted code of a machine for simulating a suitable number (encoded in form of an additional guessed part of an input) of steps of the encoded machine on $c_i$.
Moreover, it is also possible to non-deterministically accept the codes of all machines ${\cal M }$ in ${\rm CODES}_{\sigma}^{\rm Goedel}$ that do not belong to ${\rm HALT}_{\cal A}^{\subseteq\{c_i\}}$ by guessing a tuple in $U_{\cal A}^\infty\cap ( H_{\cal M} \setminus \{c_i\})$ and applying the program of a universal machine as a program segment for simulating ${\cal M}$ on this tuple. The non-deterministic semi-decidability of problems such as the complement of ${\rm HALT}_{\cal A}^{\subseteq\{c_i\}}$ by BSS RAMs in $ {\sf M}^{\rm ND}_{\cal A}$ can be done by simulating the machines whose codes are the inputs  on a guessed tuple in $U_{\cal A}^\infty \setminus \{c_i\}$ whose length can be determined in a way analogous to {\it (algo 1)} in the proof of \cite[Theorem 4.4]{GASS25B} or {\it (algo 2)} in \cite[Remark 4.5]{GASS25B}. 
By ${\rm CODES}_{\sigma}^{\rm Goedel}\setminus{\rm HALT}_{\cal A}^{\subseteq\{c_i\}}\in {\cal A}\mbox{-}{\textcolor{blue}{\Sigma_{1 }^{\rm ND}}}$ which implies ${\rm HALT}_{\cal A}^{\subseteq\{c_i\}}\in {\cal A}\mbox{-}{\textcolor{blue}{\Pi_1^{\rm ND}}}$, we obtain ${\rm HALT}_{\cal A}^{\{c_i\}}\in{\cal A}\mbox{-}{\textcolor{blue}{\Delta_{2 }^{\rm ND}}}$. More details later.
\begin{remark}[The computability of $\chi_{\{x_0\}}$]\label{SemiChar} For ${\cal A}\in {\sf Struc}_{c_1,c_2}$, the computation of the characteristic function $\chi_{\{x_0\}}$ can be performed by a BSS RAM ${\cal M}^\chi_{x_0} \in {\sf M}_{\cal A}$ using ${\cal P}_{ \textcolor{blue}{(s+2)}}$ given in Overview \ref{SemiOver0} as described in Overview \ref{SemiOver2}. ${\cal M}^\chi_{x_0}$ traverses a computation path of the form $(1,\ldots, ,t_{\vec x}-1,\textcolor{blue}{s+2},s+3,s+4)$ for any input $\vec x \in U_{\cal A}^{\infty}\setminus \{x_0\}$ and $(1,\ldots, s-1, \textcolor{red}{s},s+1, s+3,s+4)$ only for $ x_0$. 
\end{remark}
\begin{overview}[From $\{x_0\}\in {\rm SDEC}_{\cal A}$ to the compution of $\chi_{\{x_0\}}$, I]\label{SemiOver2}

\hfill

\nopagebreak 

\noindent \fbox{\parbox{11.8cm}{

Let the program of the machine ${\cal M}^\chi_{x_0} \in {\sf M}_{\cal A}$ be given by

 \vspace{0.1cm}

\hspace{1.39cm} ${\cal P}_{ \textcolor{blue}{(s+2)}}$ \,\, {\sf $\textcolor{red}{s: Z_1:= c_1^0};$ \,\, $s+1:$ goto $s+3;$ \,\, $ \textcolor{blue}{s+2: Z_1:= c_2^0};$

\hspace{1.89cm} \,\, $ s+3: I_1:= 1;$ \hspace{0.33cm} $s+4: $ stop.} .
}}
\end{overview}

\begin{remark}[$P\in {\rm DEC}_{\cal A}$ for finite sets $P\in {\rm SDEC}_{\cal A}$]\label{SemiRem} The ideas for proving the latter proposition can also be used for showing that any finite semi-decidable set is decidable and has a computable characteristic function. For the use of other algorithms, see also  Consequences \ref{Conse1} and Overview \ref{SemiOver3}.
\end{remark}

\subsubsection{The decidability of semi-decidable identity relations}\label{FromSemiToDecSecII} We will prove the statement that $(c)$ implies $(b)$. For an undecidable  identity relation with  semi-decidable complement, see Example \ref{semi_co_semi}. Let ${\cal A}$ be a structure such that ${\rm id}_{\cal A}$ is semi-decidable by a BSS RAM  ${\cal M}\in {\sf M}_{\cal A}$ whose program ${\cal P}_{\cal M}$ is a  $\sigma$-program. We assume that the axioms in $axi^{\sigma^{\circ}}$ in Overview \ref{IdentAxioms} hold for ${\rm id}_{U_{\cal A}}$ in ${\cal A}$ and that the corresponding axioms for the identity relation ${\rm id}$, its restriction ${\rm id}_{U_{\cal A}}$, and the identity relation on $\bbbn_+$ --- considered for indices stored in index registers, for all labels, and so on --- are valid in our metatheory. Recall, the transformation of one configuration into the next configuration of a deterministic machine in ${\sf M}_{\cal A}$ can be realized by applying functions defined by using operations and relations that belong to ${\cal A}$ or ${\cal A}_{\bbbn}$  (\cite{GASS25A}). Consequently, the equality of two configurations of ${\cal M}$ implies the equality of their directly following successors. Let us also note that we do not use multisets which can contain many copies of an object. This means that many copies of one individual form a set containing only one element. Moreover, we use the following agreement.
\begin{agree}\label{agree} Our copy instructions enable to copy the content of one register into another register. After the copy process, both registers have the same content. 
 A copy of an individual  that results from executing  a copy instruction or another copy command and the original are identical within ${\cal A}$ and in our world of transition systems.\footnote{There are similarities to the concept of pointers in programming, but it is not the same. We do not really use any concept of pointers known from programming. The indirect addressing can only be used in executing instructions of type (4). A modification of an object does not automatically mean a change of its copies.} \end{agree}
This means that, for $c(I_j)=2$ and $c(I_k)=1$, the execution of $Z_{I_j}:=Z_{I_k}$ leads to $c(Z_2)=c(Z_{1})$. A machine ${\cal M}\in {\sf M}_{\cal A}$ semi-deciding ${\rm id}_{U_{\cal A}}$ accepts  $(x_1,x_2)\in U_{\cal A}^2$  if $x_1=x_2$ holds and the pair $(x_1,v)\in U_{\cal A}^2$ {\sf if} $v$ is a copy of $x_1$. In both cases,  ${\rm Input}_{\cal M}$ provides the same configurations  and, consequently,   ${\sf Con}_ {{\cal M},(x_1,x_2)}$ and $ {\sf Con}_ {{\cal M},(x_1,v)}$  contain the same sequence  of configurations by the identity axioms of our metatheory. These observations leads to the following theorem. 
\begin{theorem}[Semi-deciding identity implies decidability]\label{IDTheo}\hfill 

\noindent Let ${\cal A}$ be a structure in $ {\sf Struc}$.

 The semi-decidability of ${\rm id}_{U_{\cal A}}$  by a BSS RAM in ${\sf M}_{\cal A}$ implies the semi-decidability of the complement of $\,{\rm id}_{U_{\cal A}}$  over ${\cal A}$ and thus  ${\rm id}_{U_{\cal A}}\in {\rm DEC}_{\cal A}$. 
\end{theorem} 
{\bf Proof.} Let ${\cal M}\in {\sf M}_{\cal A}$ semi-decide ${\rm id}_{U_{\cal A}}$, $(x_1,x_2)\in U_{\cal A}^2$, and $z$ be a copy of $x_2$. For $(x_1,x_2)$, we try to find out whether the behavior of ${\cal M}$ would be the same for the inputs $(x_1,x_2)$ and $(z,x_2)$  or not (which would imply $(x_1,x_2) \in {\rm id}_{U_{\cal A}}$ and $(x_1,x_2) \not\in {\rm id}_{U_{\cal A}}$, respectively). We select these inputs since, by \textcolor{blue}{{\it (ini)}}, ${\rm Input}_{\cal M}(z,x_2)$ is easier to evaluate  than   ${\rm Input}_{\cal M}(x_1,v)$ for a copy $v$ of $x_1$. The application of ${\cal P}_{\cal M}$ on ${\rm Input}_{\cal M}(z,x_2)$ leads to a stop configuration. The sequences  in ${\sf Con}_ {{\cal M},(x_1,x_2)}$ and ${\sf Con}_ {{\cal M},(z,x_2)}$  are the same if and only if $x_1=x_2$. In such a case, configuration by configuration we get, due to $axi^{\sigma^{\circ}}$, the same next configuration for both inputs. Thus, for $x_1=x_2$, both inputs traverse the same computation path of ${\cal M}$. Since ${\cal M}$ semi-decides ${\rm id}_{\cal A}$, this holds only if $x_1=x_2$. Thus, for deciding ${\rm id}_{U_{\cal A}}$, we can first simulate ${\cal M}$ on $(z,x_2)$ and count additionally each step that ${\cal M}$ would execute for input $(z,x_2)$. Because the execution of ${\cal M}$ on $(z,x_2)$ ends after a finite number $t_{x_2}$ of steps, the computation of ${\cal M}$ stops for $ {\rm Input}_{\cal M}(x_1,x_2)$ also after $t_{x_2}$ steps if $x_1=x_2$ holds. Thus, it is sufficient to simulate only $t_{x_2}$ steps of ${\cal M}$ on $(x_1,x_2)$ to recognize whether $x_1=x_2$ holds or not. If $t_{x_2}$ steps of ${\cal M}$ on $(x_1,x_2)$ do not end with the stop instruction of ${\cal M}$, then such an input can be rejected after simulating $t_{x_2}$ steps of ${\cal M}$. For implementing this algorithm for a structure  ${\cal A}$ in $ {\sf Struc}_{c_1,c_2}$, we give the program of a BSS RAM ${\cal N}_{\cal M}^{(2)}$ in ${\sf M}_{\cal A}^{(2)}$  in Overview \ref{IDOver}. Let ${\cal N}_{\cal M}^{(2)}$ be equipped with  $k_{\cal M}+1$ index registers for each tape. Then, we can start with 

\vspace{0.05cm}

${\rm Input}_{{\cal N}_{\cal M}^{(2)}}(x_1,x_2)=(1,\,(\, ((\vec \nu\,.1)\,.\, (x_1,x_2)\,.\, \bar u),\,((1,\ldots,1)\,.\,(z,x_2)\,.\, \bar u)\, ))$ \hfill \textcolor{blue}{\it (ini)}

\vspace{0.05cm}

\noindent  for $\vec \nu =(2,1,\ldots, 1)\in \bbbn_+^{k_{\cal M}}$,   $\bar u=(x_2,x_2,\ldots)$, and the copy  $z$ of $x_2$  in $Z_{2,1}$  which is very useful since $ {\rm Input}_{\cal M}(x_1,x_2)= (1\,.\, \vec \nu\,.\, (x_1,x_2)\,.\, \bar u)$ and ${\rm Input}_{\cal M}(z,x_2)= (1\,.\,\vec \nu \,.\,(z,x_2)\,.\,\bar u)$  hold.

\begin{overview}[From ${\rm id}_{U_{\cal A}}\in {\rm SDEC}_{\cal A}$ to the computation of $\chi_{{\rm id}_{U_{\cal A}}}$, I]\label{IDOver}

\hfill

\nopagebreak 

\noindent \fbox{\parbox{11.8cm}{

Let ${\cal M}\in {\sf M}_{\cal A}$ semi-decide ${\rm id}_{U_{\cal A}}$ and let ${\cal P}_{\cal M}$ be the form \textcolor{blue}{$(*)$}. Then, let

\hspace{1cm}$1: I_{2,1}:=I_{2,1}+1;$ \,$2:$\, {\sf if $I_{1,1}=I_{2,1}$ then goto $1'' $ else goto $\tilde \ell;$}

\hspace{0.76cm}$1\,'' :$ {\sf inst}$'' _1;\,$\,\,\, $2\,'':$ {\sf inst}$'' _2;\, \,\ldots;\,\hspace{0.22cm}( \ell _{{\cal P}_{\cal M}}-1)\,'' :$ {\sf inst}$'' _{ \ell _{{\cal P}_{\cal M}}-1};\,\, \hspace{0.22cm}\ell _{{\cal P}_{\cal M}}'' : $ {\sf subpr}$'';\,$ 

\qquad$1\,'\,\,:$ {\sf inst}$'_1;$\,\,\,\, $2\,'\,:$ {\sf inst}$'_2; \,\,\,\ldots;\hspace{0.22cm}\,( \ell _{{\cal P}_{\cal M}}-1)'\,\,:\,$ {\sf inst}$'_{ \ell _{{\cal P}_{\cal M}}-1};\, 
\hspace{0.22cm} {\cal P}$

\vspace{0.1cm}

be the program ${\cal P}_{{\cal N}_{\cal M}^{(2)}}$ where we use $d\,\textcolor{blue}{''}=\textcolor{blue}{2}$ and $d\, \textcolor{blue}{'}= \textcolor{blue}{1}$ and 

\vspace{0.1cm}

1) \,{\sf subpr}$''$ stands for \qquad $I_{d'', k_{\cal M}+1}:= I_{d'', k_{\cal M}+1}+1$ , 

\quad\,\,{\sf subpr}$'$ stands for \hspace{0.8cm} $I_{d', k_{\cal M}+1}:= I_{d', k_{\cal M}+1}+1$ , 

\vspace{0.1cm}

2) each subprogram {\sf inst}$'' _k$ results from {\sf instruction}$_k$ by \parskip -0.7mm 
\begin{itemize} \parskip -2mm 
\item
 replacing \textcolor{brown}{each}  \, $Z_i$ by $Z_{d'',i}$ , \, $I_i$ by $I_{d'',i}$ , \, $\ell$ by $\ell''$ \,\,and
\item adding \, {\sf subpr}$'';$ \,  in front of the instruction after the label,
\end{itemize}  \parskip -2mm 
3) each subprogram {\sf inst}$'_k$ results from {\sf instruction}$_k$ 
by\begin{itemize} \parskip -2mm 
\item replacing  \textcolor{brown}{each  \, $Z_i$ by $Z_{d',i}$ , \, $I_i$ by $I_{d',i}$ , \, $\ell$ by $\ell'$ \,\,and } \item \textcolor{brown}{adding \,\,{\sf subpr}$';$  {\sf if $I_{d'', k_{\cal M}+1}=I_{d', k_{\cal M}+1}$ then goto $\tilde \ell;$}
\,\, in front of the instruction after the label}, \end{itemize} \parskip -1.8mm 

4) the program segment ${\cal P}$ stands for

\vspace{0.1cm}

\hspace{0.89cm}$\ell _{{\cal P}_{\cal M}}': \,I_{1,1}:=1;\, Z_{1,1}:= c_1^0; \,\, $ {\sf goto $\tilde \ell^{*} ;$}

\hfill {\small \textcolor{gray}{($I_{d'', k_{\cal M}+1}\!=\!I_{d', k_{\cal M}+1}+1$ is satisfied.)}}

\vspace{0.15cm}

\hspace{1.39cm}$\!\!\tilde \ell\,\,: \,I_{1,1}:=1; \, Z_{1,1}:= c_2^0; \,\, $

\vspace{0.05cm}

\hspace{1.17cm}$\tilde \ell^* \,:\,$ {\sf stop}.
}}
\end{overview}
The program  of  ${\cal N}_{\cal M}^{(2)}$  given in Overview \ref{IDOver}  enables to compute $\chi_{{\rm id}_{U_{\cal A}}}$. For input $(x_1,x_2)$, ${\cal N}_{\cal M}^{(2)}$ can start with the initial values $(x_1,x_2,x_2, \ldots)$ stored on its tape 1 and the values $(x_2,x_2,x_2, \ldots)$ stored on its tape 2. $I_{2, k_{\cal M}+1}$ and $I_{1, k_{\cal M}+1}$ can be used for counting the computation steps of ${\cal M}$ simulated for its inputs $(z,x_2)$ and $(x_1,x_2)$, respectively. A program for semi-deciding the complement of ${\rm id}_{U_{\cal A}}$ can be derived from ${\cal P}_{{\cal N}_{\cal M}^{(2)}}$ by using 

\vspace{0.1cm}

\hspace{0.35cm}{$\ell _{{\cal P}_{\cal M}}'$ \,:\, 
\textcolor{brown}{ goto $\ell _{{\cal P}_{\cal M}}'$}};\,\, \hfill {\small \textcolor{gray}{($I_{d'', k_{\cal M}+1}\!=\!I_{d', k_{\cal M}+1}+1$ is satisfied.)}}

\hspace{0.86cm}\textcolor{brown}{$\,\tilde \ell$} : {\sf stop}.
\qed

\begin{consequ}[A further program for $\chi_{\{c_{\alpha_i}\}}$ if $\{c_{\alpha_i}\}\in {\rm SDEC}_{\cal A}$]\label{Conse1}\hfill

\noindent  If ${\cal P}_{\cal M}$ has the form \textcolor{blue}{$(*)$} and allows to semi-decide a single-element set $\{c_{\alpha_i}\}\subseteq U_{\cal A}^\infty$ and we delete the instruction \fbox{$ I_{2,1}:=I_{2,1}+1$} labeled by $1$ in Overview \ref{IDOver} and replace it by $Z_{2,1}:=c_i^0 $, then the resulting program defined analogously to ${\cal P}_{{\cal N}_{\cal M}^{(2)}}$ in Overview \ref{IDOver} allows to compute $\chi_{\{c_{\alpha_i}\}}$ (in the case that ${\cal A}\in {\sf Struc}_{c_1,c_2}$). 
\end{consequ}

\noindent {\bf Further consequences.} If $s_0$ is the length of the computation path traversed by an ${\cal M}\in {\sf M}_{\cal A}$ semi-deciding $\{x_0\}$ for input $x_0$ (which means ${\cal M}$ executes $s_{\textcolor{brown}{0}}-1$ steps on $x_0$), then, for ${\cal A}\in {\sf Struc}_{c_1,c_2}$, the program in Overview \ref{SemiOver3} allows to compute $\chi_{\{x_0\}}$ where \fbox{$I_{ k_{\cal M}+1}:=s_0$} stands for a subprogram containing $s_0-1$ times the instruction $I_{ k_{\cal M}+1}:=I_{ k_{\cal M}+1}+1$. It is possible to transfer the ideas for computing $\chi_{\{\vec x^{(1)},\ldots,\vec x^{(k)} \}}$. 
For details, see our {\sf special part} on computation paths.
\begin{overview}
[From $\{x_0\}\in {\rm SDEC}_{\cal A}$ to the compution of $\chi_{\{x_0\}}$, II]\label{SemiOver3}

\hfill

\nopagebreak 

\noindent \fbox{\parbox{11.8cm}{

Let ${\cal M}\in {\sf M}_{\cal A}$ semi-decide $\{x_0\}$ and let ${\cal P}_{\cal M}$ be the form \textcolor{blue}{$(*)$}. Then, let

\vspace{0.1cm}

\hspace{1.3cm}$1: $ {\sf if $I_{1}=I_{2}$ then goto $1^* $ else goto $\tilde \ell;$} \,\, $1^*:\,I_{ k_{\cal M}+1}:=s_0;$ 

\hspace{1.06cm}$1\,'\,:$ {\sf inst}$'_1;$\,\,\,\, $2\,'\,:$ {\sf inst}$'_2; \,\,\ldots;\hspace{0.22cm}\,( \ell _{{\cal P}_{\cal M}}-1)'\,:\,$ {\sf inst}$'_{ \ell _{{\cal P}_{\cal M}}-1};\, 
\hspace{0.22cm} {\cal P}$

\vspace{0.1cm}

be the program of a 1-tape machine ${\cal N}_{\cal M}\in {\sf M}_{\cal A}$ where

1) {\sf subpr} stands for \quad $I_{ k_{\cal M}+2}:=I_{ k_{\cal M}+2}+1$ ,

\vspace{0.1cm}

2) each subprogram {\sf inst}$' _k$ results from {\sf instruction}$_k$ by
\begin{itemize} \parskip -1.0mm 
\item
 replacing each label $\ell$ by $\ell'$ and 

adding   \textcolor{brown}{\sf \,\, if $I_{ k_{\cal M}+1}=I_{k_{\cal M}+2}$ then goto $\tilde \ell;$} \, {\sf subpr}$;$ \,\,  in front of the instruction after the label,
\end{itemize}\parskip -1.5mm 
3) the program segment ${\cal P}$ stands for

\vspace{0.05cm}

\hspace{0.89cm}$\ell _{{\cal P}_{\cal M}}': \,I_{1}:=1;\, Z_{1}:= c_1^0; \,\, $ {\sf goto $\tilde \ell^{*} ;$}
\hfill {\small \textcolor{gray}{($I_{k_{\cal M}+2}\!=\!I_{k_{\cal M}+1}$ is satisfied.)}}

\vspace{0.15cm}

\hspace{1.39cm}$\!\!\tilde \ell\,\,: \,I_{1}:=1; \, Z_{1}:= c_2^0; $ \quad $\tilde \ell^* \,:\,$ {\sf stop}.

\vspace{0.05cm}
}}
\end{overview}

\vspace{0.3cm}

\noindent
{\bf From ${\rm id}_{U_{\cal A}} \in {\rm SDEC}_{\cal A}$ to the computation of $\chi_{{\rm id}_{U_{\cal A}}}$, II.} \,
For ${\cal A} \in {\sf Struc}_{c_1,c_2}$, $\chi_{{\rm id}_{U_{\cal A}}}$ can also be computed by a 3-tape BSS RAM ${\cal N}_{\cal M}^{(3)}$ that simulates the execution of ${\cal M}$ on inputs $(x_1,x_2)$ and $(x_2,x_2)$ pseudo-parallel on the first two tapes. ${\cal P}_{{\cal N}_{\cal M}^{(3)}}$ has a similar form as the program in Overview \ref{CharactchiP1}. In case of  ${\cal P}_{{\cal N}_{\cal M}^{(3)}}$, both machines ${\cal M }'$ and ${\cal M}''$ stand for ${\cal M }$ semi-deciding ${\rm id}_{U_{\cal A}}$ and the subprogram labeled by 1 is removed.
Such an ${\cal N}_{\cal M}^{(3)}$ halts if one of both computation threads reaches the label $\ell_{{\cal P}_{\cal M}}'$ or $\ell_{{\cal P}_{\cal M}}''$. ${\cal N}_{\cal M}^{(3)}$ outputs $c_1$ if the stop instructions are reached for both threads at the same time, i.e. after the same number of steps. Otherwise, the output is $c_2$. A solution is attached, cf.\,\,Overview \ref{CharactchiId}.
For any ${\cal A}\in {\sf Struc}$, the decidability of $\,{\rm id}_{U_{\cal A}}$ can be derived as above.

\begin{example}[A co-semi-decidable and undecidable identity]\label{semi_co_semi}\hfill

\noindent
Let ${{\cal A}_4}$ be the structure $(\{a,b\}^\omega\cup \{a,b\};a,b; head, tail; r_1)$ of signature $(2;1,1;2)$ in ${\sf Struc}_{a,b}$ that contains two constants $a$ and $b$, the relation $r_1$ given by $r_1=\{(a,a),(b,b)\}$, and two operations $head$ and $tail$ defined by $head(\bar v)=v_1$ and $tail(\bar v)=(v_2,v_3,\ldots)$ for all sequences $\bar v\in \{a,b\}^\omega$ given by $\bar v=(v_1,v_2,\ldots)$ and $head(x)=x$ and $tail(x)=x$ for $x\in \{a,b\}$. A program suitably generated by deleting the second line and the six$^{\rm th}$ line in ${\cal P}_{{\cal M}_2}$ considered in Example \ref{CounterEx1} enables to semi-decide $U_{{\cal A}_4}^\infty\setminus {\rm id}_{U_{{\cal A}_4}}$ (the resulting program is below attached). However, it can be shown that, for any input $(\bar v^{(1)},\bar v^{(2)})$ in $\{a,b\}^\omega\!\times\! \{a,b\}^\omega$, any ${\cal M}\in {\sf M}_{{\cal A}_4}$ can check only the first components of two initial parts $(v_{1,1},\ldots, v_{1,s_1})$ and $(v_{2,1},\ldots, v_{2,s_2})$ of $\bar v^{(1)}$ and $\bar v^{(2)}$ for equality within a finite number $s$ of transformation (or computation) steps (where $s_1,s_2<s$). The acceptance of $(\bar v^{(1)},\bar v^{(2)})$ by ${\cal M}$ implies the acceptance of pairs that are not belong to ${\rm id}_{U_{{\cal A}_4}}$. Thus, ${\rm id}_{U_{{\cal A}_4}}$ is not semi-decidable by a BSS RAM in ${\sf M}_{{\cal A}_4}$. For similar investigations, see papers dealing with Type-2-Turing machines on the basis of \cite{WEIHRAUCH}. 
\end{example}

\subsection{Binary and other non-determinisms restricted to constants}

 If no constants are available (or there is only one constant as in the case of groups, cf.\,\cite{GASS01}), then the non-deterministic behavior of a machine can be generated by guessing labels which enable binary non-deterministic branchings. 

\begin{overview}[Instructions for binary non-deterministic branchings]

\hfill

\nopagebreak 

\noindent \fbox{\parbox{ 11.8cm}{

Instructions for guessing labels

\qquad $(11)$ \quad {\sf $\ell :\,$ goto $\ell_1$ or goto $\ell_2$}
}}
\end{overview}

\begin{example} The set $\{\vec x \in \mathbb{R}^\infty\mid (\exists \vec \alpha\in \{0,1\}^\infty )(\sum_{i=1}^n \alpha_ix_i=\pi)\}$ is, for any input $\vec x \in \mathbb{R}^n$, non-determin\-istically semi-decidable over $(\mathbb{R}; 0; +;\{\pi\})$ by executing at most $n-1$ additions and $n$ non-deterministic branchings. 
\end{example}

For ${\cal A}\in {\sf Struc}_{c_1,c_2}$, it is possible to restrict the domain of guesses to the constants $c_1$ and $c_2$ or the domain of permitted oracle sets $Q$ for applications of $\nu$ by using $Q=\{c_1,c_2\}^\infty$. 
Recall, ${\sf M}_{\cal A}^{\rm NDB}$, ${\sf M}_{\cal A}^{\nu}(Q)$, and $ {\sf M}_{\cal A}^{{\rm DND}}$ are the classes of all non-deterministic BSS RAMs ${\cal M}$ with the following properties.
\begin{itemize} 
\item \parskip -0.7mm 
${\cal M}\in {\sf M}_{\cal A}^{\rm NDB}$ can only execute instructions of types $(1)$ to $(8)$ and (11) 

and its single-valued input function is defined by \textcolor{blue}{(I1)}. 
\item ${\cal M}\in {\sf M}_{\cal A}^{\nu}(Q)$ can only execute instructions of types $(1)$ to $(8)$ and (10) 

and its single-valued input function is defined by \textcolor{blue}{(I1)}.
\item ${\cal M}\in {\sf M}_{\cal A}^{{\rm DND}}$ can only execute instructions of types $(1)$ to $(8)$ 

and its multi-valued input function is given by \textcolor{blue}{(I3$'$)}.
\end{itemize}
The input functions and output procedures are given by 

\vspace{0.1cm}

\noindent ${\rm In}_{\cal M}= \{(\vec x,(\,(n, 1,\ldots, 1)\,.\,\vec x\,.\, (x_n,x_n,x_n,\,\ldots)\,))\mid \vec x\in U_{\cal A}^\infty\}$\hfill \textcolor{blue}{(I1)}.

\vspace{0.1cm}
 
\noindent ${\rm In}_{\cal M}=\{(\vec x,(\,\underbrace{(n, 1,\ldots, 1)}_{\in(\mathbb{N}_+)^{k_{\cal M}}}.\underbrace{\vec x\,.\,\,\vec y\,\,.\,\,(x_n,x_n, \ldots)}_{\in U_{\cal A}^{\omega}}\,))\mid \vec x\in U_{\cal A}^\infty\,\,\&\,\,\vec y \in \{c_1,c_2\}^\infty\}$\hfill\textcolor{blue}{(I3$'$)}

\vspace{0.1cm}
\noindent 
and 
${\rm Output}_{\cal M}(\ell\,.\,\vec \nu\,.\,\bar u)= (u_1,\ldots, u_{\nu_1})$ for any configuration $(\ell\,.\,\vec \nu\,.\,\bar u)\in {\sf S}_{\cal M}$. In any case, we have ${\rm Res}_{\cal M}: U_{\cal A}^\infty\to {\mathfrak P}(U_{\cal A}^\infty)$. For the relationships between the several classes of result functions, see Section \ref{SecDifferent}.

\vspace{0.1cm}

Whereas $\{{\rm Res}_{\cal M}\mid {\cal M}\in {\sf M}_{\cal A}\}\not = \{{\rm Res}_{\cal M}\mid {\cal M}\in {\sf M}_{\cal A}^{\rm NDB}\}$ follows from the definition of the result functions, ${\rm SDEC}_{\cal A}={\rm SDEC}_{\cal A}^{\rm NDB}$ holds since it is possible to systematically go through all computation paths of a BSS RAM in ${\sf M}_{\cal A}^{\rm NDB}$ by selecting exactly one of the corresponding labels $\ell_1$ or $\ell_2$ for all instructions of type $(11)$ at any given time and to do this in a deterministic way by a machine in ${\sf M}_{\cal A}$ for all these instructions.
The proof of Theorem \ref{BinaryBranchingDeterm} is based on the analysis of program paths introduced in a {\sf special part}. To better understand the differences between computation and program paths, see also the flowcharts for illustrating programs in the last section. The program paths are the sequences of labels that describe the longest paths in a flowchart. 

\begin{theorem}[Replacing the guessing of labels] \label{BinaryBranchingDeterm}\hfill

\noindent Let ${\cal A}\in {\sf Struc}$.

 Then, there holds \quad ${\rm SDEC}_{\cal A}={\rm SDEC}_{\cal A}^{\rm NDB}$.
\end{theorem}
{\bf Proof.} 
 We have only to prove the inclusion ${\rm SDEC}_{\cal A}^{\rm NDB}\subseteq {\rm SDEC}_{\cal A}$.
Let ${\cal M}\in {\sf M}_{\cal A}^{\rm NDB}$. 
For any input, ${\cal M}$ can execute a sequence of instructions given by 

\vspace{0.1cm}

\quad {\sf $\ell_{i^*_\alpha}:$ goto $\ell_{j^*_\alpha}$ or goto $\ell_{k^*_\alpha};$} \hfill ($\alpha\in\{1,2,\ldots\big[,r_0\big]\}$).\hspace*{0.5cm}

\vspace{0.1cm}

\noindent This sequence also depends on the non-deterministic choices of $\ell_ {j^*_{\alpha'}}$ or $\ell_ {k^*_{\alpha'}}$ by ${\cal M}$ for $\alpha'\in\{1,\ldots,\alpha-1\}$. Let $\vec x\in H_{\cal M}\subseteq U_{\cal A}^\infty$ and $B_0^*$ be a finite computation path $(\ell_1^{(B_0^*)},\,\ldots, \ell_{t_0}^{(B_0^*)})\in {\cal L}_{\cal M}^{t_0}$ traversed by ${\cal M}$ for input $\vec x$. Moreover, let us now assume that a sequence $(\ell_{i^*_1},\ldots, \ell_{i^*_{r_0}})$ as given above is the subsequence of this computation path such that $\{\ell_{i^*_\alpha}\mid 1\leq \alpha \leq r_0\}= {\cal L}_{{\cal M},{\rm B}}\cap \{\ell_i^{(B_0^*)}\mid 1\leq i\leq t_0\}$ and $r_0\leq t_0$ hold. Then, a deterministic machine in ${\sf M}_{\cal A}$ is able to perform the accepting process done by ${\cal M}$ on $\vec x$. Let us consider  a  BSS RAM  ${\cal N}_{\cal M}\in {\sf M}_{\cal A}^{(2)}$  whose  program is given in Overview \ref{NDBtoSemiOver0} for realizing  the accepting process 

\vspace{0.1cm}

\noindent \quad\quad\,{\sf for $t:=1,2,\ldots$ do \{ for $s:= 2^t,\ldots, 2^{t+1}-1$ do 

\noindent \quad\qquad \{ evaluate $s$ and simulate $t$ steps of ${\cal M}$ by using \textcolor{blue}{${\cal P}$} \} \}}  

\vspace{0.1cm}

\noindent  and show that ${\cal N}_{\cal M}$ accepts any $\vec x \in H_{\cal M}$.
For all lengths $t+1\in \{2,3,\ldots\}$ of initial parts of program paths of ${\cal P}_{\cal M}$ and each $s\in\{2^t, \ldots, 2^{t+1}-1\}$, ${\cal N}_{\cal M}$ can use the digits $b_1,\ldots,b_t$ of the binary code $b_0\cdots b_t$ of $s$ given by $s=b_02^t+b_12^{t-1}+\cdots + b_t2^0$ with $b_0=1$. In the case of $t=2$, for example, we are interested in $s\in \{4,\ldots, 7\}$ and $b_0\cdots b_t\in \{100,\ldots, 111\}$, and so on. Corresponding to these digits of $s$, on its input $\vec x$, ${\cal N}_{\cal M}$  tries to simulate $t$ possible steps of ${\cal M}$ on $\vec x$ by executing the program segment \textcolor{blue}{${\cal P}$} defined in Overview \ref{ForSimul_NDB}. This means that ${\cal N}_{\cal M}$ executes ${\cal P}$ for all initial parts of all possible program paths including the subsequences $(\ell_{i_1},\ell_{i_2},\ldots)$ of all labels belonging to ${\cal L}_{{\cal M},{\rm B}}$ and the resulting subsequences $( \ell_{i_1}', \ell_{i_2}',\ldots)$ of successors $\ell_{i_\alpha}' \in G_{\ell_ {i_\alpha}}=\{\ell_ {j_\alpha},\ell_ {k_\alpha}\} \subseteq {\cal L}_{\cal M}$ with $(\ell_ {i_\alpha},\ell_ {i_\alpha}')\in G$ (for $G$ defined in \cite[Overview\,19]{GASS25A}). The evaluation of $s$ can be done by means of additional index registers of ${\cal N}_{\cal M}$. For suitable pseudo instructions, see also \cite{GASS20}.

\begin{overview}[The use of index registers $I_{1,2},I_{1,3},\ldots$ by ${\cal N}_{\cal M}$]\label{SimNuInstValues}

\hfill

\nopagebreak 

\noindent \fbox{\parbox{11.8cm}{

\setlength{\extrarowheight}{-2pt} 
\begin{tabular}{lll}
& {\rm Meaning}&{\rm Values} \\
$ I_{1,2}\!\!$& index $t$ for counting the outer loop passes&$1,2,\ldots$ \\
$ I_{1,3}\!\!$&index for determining the upper bound of $s\!\!\!$&$2^{t+1}-1$\\
$ I_{1,4}\!\!$& index $s$ for counting the inner loop passes&$ 2^t\,..\,2^{t+1}-1$ \\
$ I_{1,5}\!\!$& index for evaluating $s$ by using $c(I_{1,5})\,{\rm div}\,2$&$\sum_{j=0}^ib_j2^{i-j}$ ($i=0\,..\,t$)\\
$ I_{1,6}\!\!$& index  for evaluating $s$ by using $c(I_{1,5})\,{\rm mod}\,2\!\!$&$b_i$ \hfill ($i=1\,..\,t$)\\
$ I_{1,7}\!\!$& index used for checking whether $b_i=\!1$ holds\!&$1$\\

\end{tabular}
}}
\end{overview}

\begin{overview}[From an ${\rm NDB}$-machine to the semi-decision by ${\cal N}_{\cal M}$]\label{NDBtoSemiOver0}

\hfill

\nopagebreak 

\noindent \fbox{\parbox{11.8cm}{

Let ${\cal P}_{{\cal N}_{\cal M}}$ be the  program 

{\sf 
\hspace{1.88cm} ${\sf subpr}_{\rm init};$ for $ I_{1,2}:= 1,2,\ldots$ do \{ $ I_{1,3}:=2^{ I_{1,2}+1}-1;$

\hspace{3.83cm} for $ I_{1,4}:= 2^{ I_{1,2}},\ldots, I_{1,3}$ do 
 \{ $ I_{1,5}:= I_{1,4};$ 
\textcolor{blue}{${\cal P}$} \} \} 

\hspace{0.9cm} $\ell_{{\cal P}_{\cal M}}': {\sf stop}$.}  .
}}

\nopagebreak \noindent \fbox{\parbox{11.8cm}{

${\sf subpr}_{\rm init}$ is here the following subprogram for copying the input $\vec x$ of ${\cal N}_{\cal M}$. 

\hspace{2.18cm} $(Z_{2,1},\ldots, Z_{2,I_{2,1}}):= (Z_{1,1},\ldots, Z_{1,I_{1,1}})$}}
\end{overview}

\begin{overview}[The program segment ${\cal P}$ of ${\cal P}_{{\cal N}_{\cal M}}$]\label{ForSimul_NDB}

\hfill

\nopagebreak 

\noindent \fbox{\parbox{11.8cm}{

Let ${\cal M}\in {\sf M}_{\cal A}^{\rm NDB}$ and let ${\cal P}_{\cal M}$ be of the form \textcolor{blue}{$(*)$}. Then, let

\vspace{0.1cm}

\hspace{1cm} $1':$ {\sf inst}$_1; \quad\ldots;\quad (\ell_{{\cal P}_{\cal M}}-1)':$ {\sf inst}$_{\ell_{{\cal P}_{\cal M}}-1}; \quad 1^*: {\sf subpr}_{\rm init}; $

\vspace{0.1cm}

be the program segment \textcolor{blue}{${\cal P}$} of ${\cal P}_{{\cal N}_{\cal M}}$ defined as follows.

\vspace{0.1cm}

1) {\sf subpr}$_{\rm div}$ is the subprogram

\hspace{1cm} $ I_{1,5}:= I_{1,5} \,\, {\rm div} \,\, 2;$ \,\, {\sf if $ I_{1,5}= I_{1,7}$ then goto $1^*$} .

\vspace{0.2cm}

2) If {\sf instruction}$_k$ is neither of type (4), (5), (7) nor (11), then {\sf inst}$_k$ results 

 \hspace{0.3cm} from {\sf instruction}$_k$ by 

\parskip -0.7mm 
\begin{itemize} 
\item\parskip -0.7mm replacing each $I_j$ by $I_{2,j}$ and $Z_j$ by $Z_{2,j}$ and
\item
adding
 \,\, $; \,\,$ {\sf subpr$_{\rm div}$ } .
\end{itemize} 
3) If {\sf instruction}$_k$ is of type (4) or (5), then {\sf inst}$_k$ results from {\sf instruction}$_k$ by 
\begin{itemize} 
\item\parskip -0.7mm replacing $I_j$ by $I_{2,j}$, $Z_j$ by $Z_{2,j}$, and each label $\ell$ by \textcolor{blue}{$\ell'$} and
\item adding \,\, {\sf subpr}$ _{\rm div};$ \,\,
in front of the new branching instruction.
\end{itemize} 
4) If {\sf instruction}$_k$ has the form \,\fbox{goto $\ell_ i$ or goto $\ell_ j$}, then {\sf inst}$_k$ consists of
\vspace{0.15cm}

\hspace{1cm} $ I_{1,6}:= I_{1,5} \,\, {\rm mod} \,\, 2;$ \,
 {\sf subpr$_{\rm div};$ \, if $ I_{1,6}= I_{1,7}$ then goto \textcolor{blue}{$\ell_ i'$} else goto \textcolor{blue}{$\ell_ j'$}} .

\vspace{0.2cm}

5) If {\sf instruction}$_k$ is \,\fbox{$I_j:=I_j+1$}, then {\sf inst}$_k$ consists of

\vspace{0.15cm}

\hspace{1cm} $I_{2,j}:=I_{2,j}+1;$ \,\, ${\sf init}(Z_{2,I_{2,k_{\cal M}+1}});$ \,\, {\sf subpr$_{\rm div}$}  .

\vspace{0.2cm}
{\scriptsize 
\begin{tabular}{c}\hline \\Here: $(b_02^i+\!\cdots\! + b_i2^0) \,\, {\rm div} \,\, 2 = b_02^{i-1}+\!\cdots\! + b_{i-1}2^0$\,\, and \,\,$(b_02^i+\!\cdots\! + b_i2^0) \,\, {\rm mod} \,\, 2 = b_i$. 
\end{tabular}}}}\end{overview}

Under these conditions, the execution of ${\cal P}$ would be a simulation of ${\cal M}$ if ${\cal M}$ would go through an initial part $path(s)$ of length $t+1$ in  some program path $B^*$ that includes the first labels $\ell_{i_\alpha}$ that belong to ${\cal L}_{{\cal M},{\rm B}}$ and to $path(s)$ and the successors $ \ell_{i_\alpha}'$ deterministically determined by evaluating $s$. In such a case, the simulation of the computation of ${\cal M}$  causes, among other things, the generation of an initial part $path(s)$ given by $(\ell_1^{ (B^*)},\ldots,\ell_{t+1}^{ (B^*)})\in ({\cal L}_{\cal M}\setminus {\cal L}_{{\cal M},{\rm S}})^t\times {\cal L}_{\cal M} $ with $\ell_1^{ (B^*)}=1$ as follows. For all $i\leq t$, the successor $\ell_{i+1}^{ (B^*)}$ of a label $\ell_{i}^{ (B^*)}$ in $path(s)$ is uniquely determined. All digits of an $s$ are only important if {\sf all} instructions considered are of type $(11)$. If  ${\sf instruction}_{\ell_{i}}$ is of type (4), then the value of $\ell_{i+1}^{ (B^*)}$ only depends on the truth value of the condition in ${\sf instruction}_{\ell_{i}}$ which in turn depends on the current configuration. If $\ell_{i}^{ (B^*)}\not \in {\cal L}_{{\cal M},{\rm H}_{\rm T}}\cup{\cal L}_{{\cal M},{\rm T}}\cup {\cal L}_{{\cal M},{\rm B}}\cup {\cal L}_{{\cal M},{\rm S}}$, then $\ell_{i+1}^{ (B^*)}=\ell_i^{ (B^*)}+1$. If $\ell_{i}^{ (B^*)}$ is a label $\ell_{i_\lambda}$ in $ {\cal L}_{{\cal M},{\rm B}}$ and $G_{ \ell_{i_\lambda}}=\{ \ell_{j_\lambda }, \ell_{k_\lambda}\}$ holds, then $(1-b_i)\cdot \ell_{j_\lambda }+ b_i \cdot \ell_{k_\lambda} \in G_{ \ell_{i_\lambda}} $ (derived from $s=b_02^t+b_12^{t-1}+\cdots + b_t2^0$) is to be used after the evaluation of the instruction labeled by $ \ell_ i^{(B^*)}$ as successor $ \ell_ {i+1}^{(B^*)}$. 

The loops determined by increasing $t$ and $s$  cause only deterministic jumps by ${\cal N}_{\cal M}$. Any sequence of instructions simulated for a $t$ and an $s$ could correspond to an initial part of or a whole computation path of ${\cal M}$. Finally,  after simulating $t_0'$ possible steps of ${\cal M}$ on  $\vec x\in H_{\cal M}$  for some $t_0'$, one of the considered sequences of instructions with labels in  $B'$ given by $(\ell_ 1^{(B')},\ldots, \ell_{t_0'}^{(B')})\in {\cal L}_{\cal M}^{t_0'}$  ends with the stop instruction of ${\cal M}$.  Such a  $B'$ contains a subsequence in ${ \cal L}_{{\cal M},{\rm B}}^{r_0'}$ (as considered above) and is (like $B_0^*$)  a finite computation path of ${\cal M}$  (that can also be equal to $B_0^*$) which means that $\vec x$ is then accepted by ${\cal M}$ and also by ${\cal N}_{\cal M}$. 

 By definition of ${\cal P}_{{\cal N}_{\cal M}}$,  ${\cal N}_{\cal M}$ cannot accept an input  if it does not simulate the stop instruction of ${\cal M}$ for this input.
\qed

\subsection{Different non-determinisms and relationships}\label{SecDifferent}

 We will see that ${\rm SDEC}_{\cal A} \subseteq ({\rm SDEC}_{\cal A}^\nu)^Q={\rm SDEC}_{\cal A}^{\rm DND}\subseteq{\rm SDEC}_{\cal A}^{\rm ND}$ is valid if one of the properties $(a)$ to $(e)$ in Proposition \ref{BasicProp} and $Q=\{c_1,c_2\}^\infty$ or $Q=\{c_1,c_2\}^2$ hold.  It is, for example, easy to show that $({\rm SDEC}_{\cal A}^\nu)^Q={\rm SDEC}_{\cal A}^{\rm NDB}$ holds if ${\cal A}\in{\sf Struc}_{c_1,c_2}$ contains the identity relation. 
\begin{proposition}[Basic properties: Semi-decidable problems]\label{BasicPropSemiD}
\hfill

\noindent Let ${\cal A}$ $\in{\sf Struc}_{c_1,c_2}$, $\{c_1,c_2\}\in {\rm SDEC}_{\cal A}$, $Q=\{c_1,c_2\}^\infty$ or $Q=\{c_1,c_2\}^2$, $P\subseteq U_{\cal A}^\infty$, and 

\vspace{0.2cm}
\qquad \begin{tabular}{lcllcl}
 \quad &{\rm (i)}& $P\in {\rm SDEC}_{\cal A}$,
& {\rm (iii)}& $ P\in {\rm SDEC}_{\cal A}^{\rm DND}$, \\ 
& {\rm (ii)}& $P\in ({\rm SDEC}_{\cal A}^\nu)^Q$, \hspace{0.7cm}
& {\rm (iv)}& $ P\in {\rm SDEC}_{\cal A}^{\rm ND}$.
\end{tabular} 

\vspace{0.3cm}

 Then, \quad {\rm (i)} implies {\rm (ii)}, \quad {\rm (ii)} implies {\rm (iii)}, \quad and \quad {\rm (iii)} implies {\rm (iv)}.
\end{proposition}

Instead of proving this proposition directly, we compare the classes of result functions for the machines in ${\sf M}_{\cal A}^{\rm NDB}$, ${\sf M}_{\cal A}^{\nu}(Q)$, $ {\sf M}_{\cal A}^{{\rm DND}}, \ldots$ for structures ${\cal A}\in{\sf Struc}_{c_1,c_2}$. In this way, we also provide the proofs for \cite[Theorem 2.8]{GASS25A} and \cite[Theorem 2.10]{GASS25A} and characterize in particular Moschovakis' operator applied to oracle sets $Q=\{c_1,c_2\}^\infty$ or $Q=\{c_1,c_2\}^2$.
\begin{theorem}[Digitally non-deterministic computability]\label{TheoDigitOracle}\hfill

\noindent Let ${\cal A}$ be a structure in ${\sf Struc}_{c_1,c_2}$, let  $Q=\{c_1,c_2\}^\infty$, and let $Q\in {\rm SDEC}_{\cal A}$. Then, there holds

\vspace{0.3cm}

 \quad  $\{{\rm Res}_{\cal M}\mid {\cal M}\in {\sf M}_{\cal A}^\nu(Q)\}=\{{\rm Res}_{\cal M}\mid {\cal M}\in {\sf M}_{\cal A}^{\rm DND}\}$.
\end{theorem}
{\bf Proof.} This equation can be derived analogously to the equation in Theorem 4.7 in Part II\,a\, by using the ideas of the proofs of Theorems 4.4 and 4.6. \qed

\vspace{0.3cm}

Since $\{c_1,c_2\}^\infty\in {\rm SDEC}_{\cal A}$ is valid if and only if $\{c_1,c_2\}\in {\rm SDEC}_{\cal A}$ is valid, we have the following. 

\begin{corollary} [Digitally non-deterministic computability]\label{CoroDigitOracle} \hfill 

\noindent  Let ${\cal A}$ be  in $ {\sf Struc}_{c_1,c_2}$, $Q_0=\{c_1,c_2\}^2$, and $\{c_1,c_2\}\in {\rm SDEC}_{\cal A}$. Then, there holds

\vspace{0.2cm}

 \quad $\{{\rm Res}_{\cal M}\mid {\cal M}\in {\sf M}_{\cal A}^\nu(Q_0)\}=\{{\rm Res}_{\cal M}\mid {\cal M}\in {\sf M}_{\cal A}^{\rm DND}\}$.
\end{corollary}

\begin{theorem}[Restriction of guesses]\label{TheoDND_ND}\hfill

\noindent Let ${\cal A}\in {\sf Struc}_{c_1,c_2}$ and $\{c_1,c_2\}\!\in {\rm SDEC}_{\cal A}$. Then, there holds

\vspace{0.3cm}

 \quad $\{{\rm Res}_{\cal M}\mid {\cal M}\in {\sf M}_{\cal A}^{\rm DND}\}\subseteq \{{\rm Res}_{\cal M}\mid {\cal M}\in {\sf M}_{\cal A}^{\rm ND}\}$. 
\end{theorem}
{\bf Proof.} Any ${\cal M}\in {\sf M}_{\cal A}^{\rm DND}$ can be simulated by a machine ${\cal N}_{\cal M}\in {\sf M}_{\cal A}^{\rm ND}$ that executes a subprogram for recognizing whether a value is in $\{c_1,c_2\}$ before one of the guesses is used. ${\cal N}_{\cal M}$ loops forever if such a guessed value does not belong to $\{c_1,c_2\}$.\qed

\begin{corollary}\label{CoroDND_ND} Let ${\cal A}\in {\sf Struc}_{c_1,c_2}$ and let one of the properties $(a)$ to $(e)$ in Proposition \ref{BasicProp} hold. Then, there holds

\vspace{0.3cm}

 \quad ${\rm SDEC }_{\cal A}^{\rm DND}\subseteq {\rm SDEC }_{\cal A}^{\rm ND}$. 
\end{corollary}

\begin{theorem}[Binary non-deterministic computability]\label{BinaryOracle}\hfill

\noindent Let ${\cal A}$ be in ${\sf Struc}_{c_1,c_2}$ and $Q=\{c_1,c_2\}^\infty$. Then, there holds

\vspace{0.3cm}

 \quad  $\{{\rm Res}_{\cal M}\mid {\cal M}\in {\sf M}_{\cal A}^\nu(Q)\}\subseteq \{{\rm Res}_{\cal M}\mid {\cal M}\in {\sf M}_{\cal A}^{\rm NDB}\}$. 
\end{theorem}
{\bf Proof.} Any ${\cal M}\in {\sf M}_{\cal A}^\nu(Q) $ can be simulated by an ${\cal N}_{\cal M}$ in ${\sf M}_{\cal A}^{\rm NDB}$ that simulates the execution of each program segment \fbox{$\ell_i:\,\, Z_j:=\nu[{\cal O}](Z_1,\ldots,Z_{I_1})$;} in ${\cal P}_{\cal M}$ (with $1\leq i\leq s<\ell_{{\cal P}_{\cal M}}$)
by executing 

\vspace{0.2cm}

{\sf \hspace{0.3cm} $\ell_i:\, $ {\sf goto $\ell_i'$ or goto $\ell_i'';$}

\hspace{0.3cm} $ \ell_i':\, Z_j:= c_1^0;$ goto $\ell_i+1$ or goto $\ell_i+1; $ 

\hspace{0.24cm} $ \ell_i'':\, Z_j:= c_2^0;$ goto $\ell_i+1$ or goto $\ell_i+1;$}

\vspace{0.3cm}

\noindent where we again assume that $\ell_i'$ and  $\ell_i''$ are new labels for all $i\in \{1,\ldots, s\}$ such that $\bigcap_{i=1}^{s}\{\ell_i',\ell_i''\}\cap {\cal L}_{\cal M}=\emptyset$ holds. 
\qed

\begin{theorem}[Identity and oracle instructions]\label{IdentBinary} \hfill

\noindent Let ${\cal A}\in {\sf Struc}_{c_1,c_2}$, let ${\cal A}$ contain ${\rm id}_{\cal A}$, and let $Q=\{c_1,c_2\}^\infty$. Then, there holds

\vspace{0.3cm}

\quad $\{{\rm Res}_{\cal M}\mid {\cal M}\in {\sf M}_{\cal A}^{\rm NDB}\} \subseteq \{{\rm Res}_{\cal M}\mid {\cal M}\in {\sf M}_{\cal A}^\nu(Q)\}$. 
\end{theorem}
{\bf Proof.} Let ${\cal M}$ be a BSS RAM in $ {\sf M}_{\cal A}^{\rm NDB}$. The execution of an instruction given by \fbox{$\ell:\, ${\sf goto $\ell_1$ or goto $\ell_2$}$;$} in ${\cal P}_{\cal M}$ can be simulated by an ${\cal N}_{\cal M}$ in ${\sf M}_{\cal A}^\nu(Q)$ as follows. For simplicity, we here assume that $Z_j,Z_{j+1}, Z_{j+2}$ are not used until this time which means that $c(Z_j)=c(Z_{j+1})=c(Z_{j+2})=x_n$ holds if $\vec x \in U_{\cal A}^n$ is the input for ${\cal M}$ and ${\cal N}_{\cal M}$. Therefore, let $c(I_{k_{\cal M}+1})=j, \, \ldots$  hold before executing  pseudo-instructions such as $Z_j:= Z_{j+2}$ in the  process  of simulating  the branching process. Moreover, let the labels $\ell',\ell'', \ell_1', \ell_2'$ of the following subprograms be new labels depending on $\textcolor{red}{\ell}$. 

\vspace{0.2cm}

{\sf \hspace{0.4cm} $\textcolor{red}{\ell}:\, Z_j:= c_1^0; $

\hspace{0.3cm} $ \ell' : \, Z_{j+1}:={\nu}[{\cal O}](Z_j); $

\hspace{0.22cm} $ \ell'' : \,$ if $Z_j=Z_{j+1}$ then goto $\ell_1'$ else goto $\ell_2'; $

\hspace{0.23cm} $ \ell_1':\, Z_j:= Z_{j+2}; \,\, Z_{j+1}:= Z_{j+2};$ \,\, goto \textcolor{blue}{$\ell_1$}$;$ 

\hspace{0.23cm} $ \ell_2':\, Z_j:= Z_{j+2}; \,\, Z_{j+1}:= Z_{j+2};$ \,\, goto \textcolor{blue}{$\ell_2$}$;$}

\vspace{0.3cm}

\noindent Before the simulation is continued with an instruction labeled by $\ell_1$ and $\ell_2$, respectively, $Z_j$ and $Z_{j+1}$ obtain their old value $x_n$.

\begin{figure}[t]
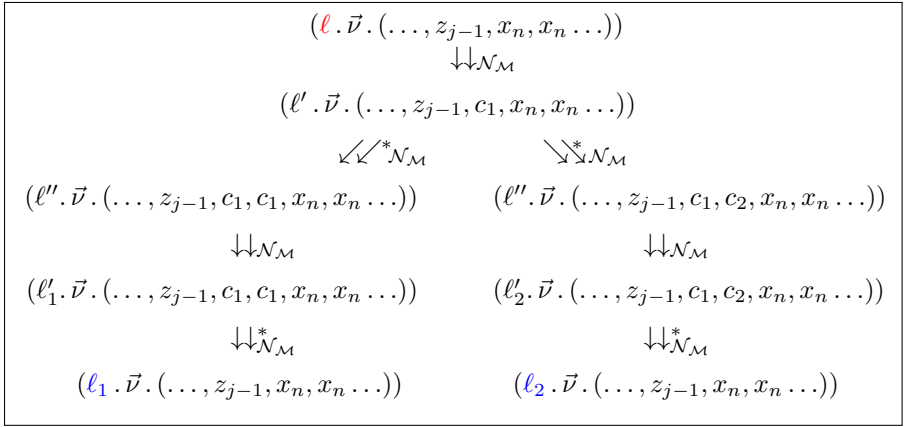

\noindent \,\fbox{\vspace{0.2cm}
\begin{tabular}{c}{\scriptsize}
\,\,$(\textcolor{red}{\ell}\,.\,\vec \nu\,.\,(\ldots, z_{j-1},x_n,x_n\ldots))$\\\vspace{0.2cm}
$\qquad \downarrow\downarrow _{{\cal N}_{\cal M}}$\\ \vspace{0.2cm}
$(\ell'\,.\,\vec \nu\,.\,(\ldots, z_{j-1},c_1, x_n,x_n\ldots))$\\ \vspace{0.2cm}
$\qquad\swarrow \!\!\swarrow\!\!^*\!_{{\cal N}_{\cal M}} \qquad\qquad \searrow\!\!\searrow\!\!\!\!\!^*\, \hspace{0.01cm}_{{\cal N}_{\cal M}}$\\ \vspace{0.2cm}
$\!\!\!(\ell''.\,\vec \nu\,.\,(\ldots, z_{j-1},c_1,c_1, x_n,x_n\ldots)) $ \qquad \, $(\ell''.\,\vec \nu\,.\,(\ldots, z_{j-1},c_1,c_2, x_n,x_n\ldots) )\!\!$\\ \vspace{0.2cm}
$\quad \downarrow\downarrow _{{\cal N}_{\cal M}}$\hspace{4.6cm}$ \downarrow\downarrow _{{\cal N}_{\cal M}}$\\ \vspace{0.2cm} 
$\!\!\!(\ell_1'.\,\vec \nu\,.\,(\ldots, z_{j-1},c_1,c_1, x_n,x_n\ldots)) $ \qquad\, $(\ell_2'.\,\vec \nu\,.\,(\ldots, z_{j-1},c_1,c_2, x_n,x_n\ldots) )\!\!$\\ \vspace{0.2cm}
$\quad \downarrow\downarrow^* _{{\cal N}_{\cal M}}$\hspace{4.6cm}$ \downarrow\downarrow^* _{{\cal N}_{\cal M}}$\\\vspace{0.2cm} 
$(\textcolor{blue}{\ell_1}\,.\,\vec \nu\,.\,(\ldots, z_{j-1},x_n,x_n\ldots)) $\qquad \qquad \,$(\textcolor{blue}{\ell_2}\,.\,\vec \nu\,.\,(\ldots, z_{j-1}, x_n,x_n\ldots) )$\quad\\
\end{tabular}}
\caption{A simulation of the execution of $\textcolor{red}{\ell}:\, ${\sf goto \textcolor{blue}{$\ell_1$} or goto \textcolor{blue}{$\ell_2$}}\label{Ableitung}}
\end{figure}

\noindent The simulation of the branching process is represented in Figure \ref{Ableitung} by a derivation or transformation tree. The symbol $\downarrow\downarrow_{{\cal N}_{\cal M}}$ serves to illustrate $\ \genfrac{}{}{0pt}{2}{\longrightarrow} {\longrightarrow}_ {{\cal N}_{\cal M}}$. $\downarrow\downarrow_{{\cal N}_{\cal M}}^*,\, \swarrow \!\swarrow_{{\cal N}_{\cal M}}^*, \ldots $ serve to illustrate the transitive closure $\genfrac{}{}{0pt}{2}{\longrightarrow} {\longrightarrow}^*_ {{\cal N}_{\cal M}}$ of $ \genfrac{}{}{0pt}{2}{\longrightarrow} {\longrightarrow}_ {{\cal N}_{\cal M}}$. \qed

\vspace{0.2cm}

The following statement is more general. 
\begin{theorem}[Recognizable constants and oracle instructions]\label{TheoRecBina}\hfill

\noindent Let ${\cal A}$ be in ${\sf Struc}_{c_1,c_2}$ and $Q=\{c_1,c_2\}^2$. Moreover, let $\{c_1\}\in {\rm DEC}_{\cal A}$ or let both $\{c_1\}\in {\rm SDEC}_{\cal A}$ and $\{c_2\}\in {\rm SDEC}_{\cal A}$ hold. Then, there holds

\vspace{0.25cm}

\quad $\{{\rm Res}_{\cal M}\mid {\cal M}\in {\sf M}_{\cal A}^{\rm NDB}\} \subseteq \{{\rm Res}_{\cal M}\mid {\cal M}\in {\sf M}_{\cal A}^\nu(Q)\}$.
\end{theorem}
{\bf Proof.} 
The execution of the program segment \fbox{$\ell:\, ${\sf goto $\ell_1$ or goto $\ell_2;$}} by a BSS RAM ${\cal M}\in {\sf M}_{\cal A}^{\rm NDB}$ can be simulated by a 3-tape $\nu$-oracle BSS RAM ${\cal N}_{\cal M}$ using the oracle $Q$  as follows.  The  first both tapes of ${\cal N}_{\cal M}$  are used only for the simulation of the non-deterministic branching processes  and each of these  processes is  replaced by a pseudo-parallel simulation of two machines as demonstrated in Figure \ref{ParaSim}. Here, let $P'$ and $P''$ be semi-decidable sets satisfying  $P'=\{c_1\}$ and $\{c_2\}\subseteq P''\subseteq U_{\cal A}^\infty\setminus \{c_1\}$. Thus the cases $P''=\{c_2\}$ and $P''=U_{\cal A}^\infty \setminus \{c_1\}$ are included. Moreover, let ${\cal M}'\in {\sf M}_{\cal A}$ and ${\cal M}''\in {\sf M}_{\cal A}$  semi-decide $\{c_1\}$ and $P''$, respectively. For simulating these machines, the 3-tape BSS RAM ${\cal N}_{\cal M}$  uses the  program given  in Overview \ref{Para2} and the  instructions that are  defined as  described in Overviews \ref{CharactchiP2} and \ref{CharactchiP3} where, additionally,   $;\, {\sf init}(Z_{1,I_{1,k_{{\cal M}'+1}}})$ is  added to any instruction of the form \fbox{$I_{1,j}:=I_{1,j}+1$} and  $;\, {\sf init}(Z_{2,I_{2,k_{{\cal M}''+1}}})$ is added to any instruction of the form \fbox{$I_{2,j}:=I_{2,j}+1$}.

\begin{overview}[To the next label by semi-deciding $\{c_1\}$ and $P''$]\label{Para2}

\hfill

\nopagebreak 

\noindent \fbox{\parbox{11.8cm}{{\sf \begin{tabular}{rl}

\quad $ \textcolor{red}{\ell} : $&$\,Z_{1,1}:= c_1^0; $ \,\, $ Z_{1,1}:={\nu}[{\cal O}](Z_{1,1});$ \,\, $Z_{2,1}:=Z_{1,1};$\\

&$I_{1,1}:=1; \quad \ldots ; \,\, I_{1,k_{{\cal M}'+2}}\,:=1;$ \\

& $I_{2,1}:=1; \quad \ldots ; \,\, I_{2,k_{{\cal M}''+2}}:=1;$\\

$\tilde 1\,':$ & If $I_{2,k_{{\cal M}'+3}}=i$ then goto $i'';$\\

$\tilde 1\,'':$ & If $I_{1,k_{{\cal M}''+3}}=i$ then goto $i';$\\

$1\,':$ & {\sf inst}$'_1; \quad \,\ldots;\quad ( \ell _{{\cal P}_{{\cal M}'}}-1)\,'\,\,\,$: {\sf inst}$'_{ \ell _{{\cal P}_{{\cal M}'}}-1};\,$\quad $ \ell\,' _{{\cal P}_{{\cal M}'}}\,\,$: {\sf goto} \textcolor{blue}{$\ell _1$}$;$\\

$1\,'':$ & {\sf inst}$''_1; \quad \,\ldots;\quad ( \ell _{{\cal P}_{{\cal M}''}}-1)\,''$: {\sf inst}$''_{ \ell _{{\cal P}_{{\cal M}''}}-1};$\quad $ \ell\,'' _{{\cal P}_{{\cal M}''}}$: {\sf goto} \textcolor{blue}{$\ell _2$}$;$\vspace{0.1cm}\\
\end{tabular}}}}
\end{overview}
The execution of all instructions of the types (1) to (8) by ${\cal M}$ can be simulated by means of the third tape and by executing the  instructions ${\sf inst}_k$ that result from ${\sf instruction}_k$ in ${\cal P}_{\cal M}$ by replacing all \fbox{$I_j$} by $I_{3,j}$ and \fbox{$Z_{j}$} by $Z_{3,j}$. For copying the input of ${\cal N}_{\cal M}$ onto the third tape, the subprogram $(Z_{3,1},\ldots,Z_{3,I_{3,1}}):=(Z_{1,1},\ldots,Z_{1,I_{1,1}})$ is added in front of all. Consequently, there is also a 1-tape machine in ${\sf M}_{\cal A}^\nu(Q)$ that can simulate $\!{\cal M}$. \qed

\vspace{0.1cm}

By this and Theorem \ref{TheoDigitOracle}, we get the following.
\begin{corollary}[Binary non-deterministic semi-decidability]\label{CorIdenRecBina} \hfill

\noindent 
 Let ${\cal A}$ be a structure in $ {\sf Struc}_{c_1,c_2}$ and $Q=\{c_1,c_2\}^\infty$ or $\,Q=\{c_1,c_2\}^2$. Moreover, let $\{c_1\}\in {\rm DEC}_{\cal A}$ or let both $\{c_1\}\in {\rm SDEC}_{\cal A}$ and $\{c_2\}\in {\rm SDEC}_{\cal A}$ hold. Then, there holds

\vspace{0.2cm}

\quad $({\rm SDEC}_{\cal A}^\nu)^Q={\rm SDEC}_{\cal A}^{\rm NDB}$.
\end{corollary}

\newpage
\section*{Summary}
\addcontentsline{toc}{section}{\bf Summary}
\markboth{SUMMARY}{SUMMARY}

\begin{overview}[Basic implications: Identity and constants]\label{BasicImplicI}

\hfill

\nopagebreak 

\noindent \fbox{\parbox{11.8cm}{

\hfill {\small Let ${\cal A}\in {\sf Struc}_{c_1,c_2}$.}

\qquad $(a)$ \, ${\cal A}$ contains the identity ${\rm id}_{U_{\cal A}}$.

\qquad $(b)$ \, The identity is decidable by a BSS RAM in ${\sf M}_{\cal A}$.

\qquad $(c)$ \, The identity is semi-decidable by a BSS RAM in ${\sf M}_{\cal A}$. 

\qquad $(d)$ \, $\{c_1\}$ and $\{c_2\}$ are semi-decidable by two BSS RAMs in ${\sf M}_{\cal A}$. 

\qquad $(e)$ \, $\{c_1,c_2\}^\infty$ is semi-decidable by a BSS RAM in ${\sf M}_{\cal A}$.

\vspace{0.2cm}

There hold

\vspace{0.1cm}

\hspace{1.5cm}$(a) \,\, \textcolor{blue}{\Rightarrow} \,\, (b) \,\, \textcolor{blue}{\Rightarrow} \,\, (c) \,\, \textcolor{blue}{\Rightarrow} \,\, (d) \,\, \textcolor{blue}{\Rightarrow} \,\, (e)$
\hfill{\scriptsize (Prop. \ref{BasicProp})}

\vspace{0.1cm}
{\small where, $(a) \Rightarrow (b)$ stands for the implication that $(a)$ implies $(b)$, and so on,

and we have}

\vspace{0.2cm}

\hspace{0.53cm }\begin{tabular}{ccl}
\qquad $(b)$ & ${\rm id}_{U_{\cal A}}$&$\in {\rm DEC}_{\cal A}$,\\
\qquad $(c)$ & ${\rm id}_{U_{\cal A}}$&$\in {\rm SDEC}_{\cal A}$,\\
\qquad $(d)$ & $\{c_1\}, \{c_2\}$&$\in {\rm SDEC}_{\cal A}$,\\
\qquad $(e)$ & $\{c_1, c_2\}^\infty $&$\in {\rm SDEC}_{\cal A}$,\\
\end{tabular}

\vspace{0.2cm}

and there are ${\cal A}\in {\sf Struc}_{c_1,c_2}$ with 

\vspace{0.1cm}

\hspace{1.5cm}$(a) \, \textcolor{blue}{\hspace{0.13cm}\not\hspace{0.03cm} \!\! \Leftarrow} \,\, (b) \mbox{\quad or \quad} (c) \textcolor{blue}{\hspace{0.13cm}\not\hspace{0.03cm} \!\! \Leftarrow} \,\, (d)\mbox{\quad or \quad} (d) \,\, \textcolor{blue}{\hspace{0.13cm}\not\hspace{0.03cm} \!\! \Leftarrow} \,\, (e)$.
\hfill{\scriptsize (Expls. \ref{CounterEx1}, \ref{CounterEx2}, \ref{CounterEx3})}

\vspace{0.1cm}
}}
\end{overview}

\begin{overview}[Some basic properties and some equivalences]\label{BasicImplicII}

\hfill\nopagebreak 

\noindent \fbox{\parbox{11.8cm}{
\hfill {\small Let ${\cal A}\in {\sf Struc}_{c_1,c_2}$.}

\vspace{0.2cm}

\qquad $(d) $ $ \quad \textcolor{blue}{\Leftrightarrow} \quad$ $\{c_1\}, \{c_2\}\in {\rm DEC}_{\cal A}$. 
\hfill{\scriptsize (cf.\,Prop. \ref{SemiDeci})}

\vspace{ 0.2cm}

\hspace{1.17cm} $ \quad \textcolor{blue}{\Leftrightarrow} \quad$ $\chi_{\{c_1\}}, \, \chi_{\{c_2\}}\in \{{\rm Res}_{\cal M}\mid {\cal M}\in {\sf M}_{\cal A}\}$. 
\hfill{\scriptsize (cf.\,Rem. \ref{SemiChar})}

\vspace{ 0.2cm}

\qquad $(e) $ $ \quad \textcolor{blue}{\Leftrightarrow} \quad$ {$\{c_1,c_2\}\in {\rm DEC}_{\cal A}$.} \hfill{\scriptsize (cf.\,Rem. \ref{SemiRem})}

\vspace{ 0.2cm}

\qquad $(b) $ $ \quad \textcolor{blue}{\Leftrightarrow} \quad$ $(c) $.
\hfill{\scriptsize (Theo. \ref{IDTheo})}

\vspace{0.3cm}

{\small Here, $(a) \Leftrightarrow (b)$ means that $(a) \Rightarrow (b)$ and $(b) \Rightarrow (a)$ hold, and so on.}}}
\end{overview}

\begin{overview}[Decidability of $P$ and computability of $\chi_P$]\label{DecAndChar}

\hfill

\nopagebreak 

\noindent \fbox{\parbox{11.8cm}{

\hfill {\small Let ${\cal A}\in {\sf Struc}_{c_1,c_2}$ and $P\subseteq U_{\cal A}^\infty$.}

\vspace{0.2cm}

\qquad $\chi_P\in \{{\rm Res}_{\cal M}\mid {\cal M}\in {\sf M}_{\cal A}\}$ and $\,(a)$ $ \quad \textcolor{blue}{\Rightarrow} \quad$ $P\in {\rm DEC}_{\cal A}$. 
\hfill{\scriptsize (Prop. \ref{Propos1})}

\qquad $\chi_P\in \{{\rm Res}_{\cal M}\mid {\cal M}\in {\sf M}_{\cal A}\}$ and $\,(d)$ $ \quad \textcolor{blue}{\Rightarrow} \quad$ $P\in {\rm DEC}_{\cal A}$.
\hfill{\scriptsize (Prop. \ref{Propos2})}

\qquad $\chi_P\in \{{\rm Res}_{\cal M}\mid {\cal M}\in {\sf M}_{\cal A}\}$ and $\,(b)$ $ \quad \textcolor{blue}{\Rightarrow} \quad$ $P\in {\rm DEC}_{\cal A}$.
\hfill{\scriptsize (Coro. \ref{CorollPropos2})}

\vspace{0.4cm}

\qquad $P\in {\rm DEC}_{\cal A}$ $ \quad \textcolor{blue}{\Rightarrow} \quad$ $\chi_P\in \{{\rm Res}_{\cal M}\mid {\cal M}\in {\sf M}_{\cal A}\}$.
\hfill{\scriptsize (Prop. \ref{Propos3})}

\vspace{0.2cm}}}
\end{overview}

\vspace{0.2cm}

\begin{overview}[Non-determinisms discussed in Part II\,a]

\hfill

\nopagebreak 
\noindent \fbox{\parbox{11.8cm}{\hfill {\small Let ${\cal A}\in {\sf Struc}$.} 

Non-deterministic computability and semi-decidability 

\vspace{0.1cm}

{\small

\hspace{0.8cm} Let $Q \,\,\textcolor{blue}{\subseteq}\,\, U_{\cal A}^\infty$ and $Q\in {\rm SDEC}_{\cal A}$.

\vspace{0.1cm}\hspace{1cm} $\textcolor{blue}{\Rightarrow}$ 

\hspace{3.23cm}$({\rm SDEC}_{\cal A}^\nu)^Q\hspace{0.15cm} \,\,\textcolor{blue}{\subseteq}\,\,\hspace{0.15cm}{\rm SDEC}_{\cal A}^{\rm ND}$, \hfill{\scriptsize (Theo.\,4.4, Rem.\,4.5)}

\hspace{1.6cm}$ \{{\rm Res}_{\cal M}\mid {\cal M}\in {\sf M}_{\cal A}^\nu(Q)\}\hspace{0.15cm} \,\,\textcolor{blue}{\subseteq}\,\,\hspace{0.15cm}\{{\rm Res}_{\cal M}\mid {\cal M}\in {\sf M}_{\cal A}^{\rm ND}\}$. 
\hfill{\scriptsize (Theo.\,4.6)}

\vspace{0.4cm}

\hspace{0.8cm} Let $Q \,\,\textcolor{blue}{=}\,\, U_{\cal A}^\infty$.

\vspace{0.1cm}\hspace{1cm} $\textcolor{blue}{\Rightarrow}$ 

\hspace{3.23cm}$ ({\rm SDEC}_{\cal A}^\nu)^Q\hspace{0.15cm} \,\,\textcolor{blue}{=}\,\,\hspace{0.15cm}{\rm SDEC}_{\cal A}^{\rm ND}$, \hfill{\scriptsize (Theo.\,4.7)}

\vspace{0.1cm}

\hspace{1.6cm}$ \{{\rm Res}_{\cal M}\mid {\cal M}\in {\sf M}_{\cal A}^\nu(Q)\}\hspace{0.15cm} \,\,\textcolor{blue}{=}\,\,\hspace{0.15cm}\{{\rm Res}_{\cal M}\mid {\cal M}\in {\sf M}_{\cal A}^{\rm ND}\}$. 
\hfill{\scriptsize (Theo.\,4.7)}

\vspace{0.2cm}}}}
\end{overview}

\begin{overview}[Non-determinisms --- further relationships] \hfill

\nopagebreak

\noindent \fbox{\parbox{11.8cm}{

\hfill {\small Let ${\cal A}\in {\sf Struc}_{c_1,c_2}$.} 

Non-deterministic computability 

{\small
\quad $(e)$ 
 $\,\,\textcolor{blue}{\Rightarrow}$ \quad$\{{\rm Res}_{\cal M}\mid {\cal M}\in {\sf M}_{\cal A}^{\rm DND}\}\hspace{0.3cm} \,\,\textcolor{blue}{\subseteq}\,\,\hspace{0.15cm} \{{\rm Res}_{\cal M}\mid {\cal M}\in {\sf M}_{\cal A}^{\rm ND}\}$. \hfill{\scriptsize (Theo.\,\ref{TheoDND_ND})}

\vspace{0.2cm}
}}}

\nopagebreak

\noindent \fbox{\parbox{11.8cm}{

\hfill {\small Let ${\cal A}\in {\sf Struc}_{c_1,c_2}$ and $Q_0 \,\,\textcolor{blue}{=}\,\,\{c_1,c_2\}^2$.}

Digitally non-deterministic computability

\vspace{0.1cm}

{\small
\quad $(e)$ 
 $\,\,\textcolor{blue}{\Rightarrow}$ \quad$\{{\rm Res}_{\cal M}\mid {\cal M}\in {\sf M}_{\cal A}^\nu(Q_0)\} \,\,\textcolor{blue}{=}\,\,\hspace{0.15cm}\{{\rm Res}_{\cal M}\mid {\cal M}\in {\sf M}_{\cal A}^{\rm DND}\}$. \hfill{\scriptsize (Coro.\,\ref{CoroDigitOracle})} 

\vspace{0.2cm}}

Binary non-deterministic computability

{\small
\vspace{0.1cm}

\hspace{1.6cm}$ \{{\rm Res}_{\cal M}\mid {\cal M}\in {\sf M}_{\cal A}^\nu(Q_0)\} \,\,\textcolor{blue}{\subseteq}\,\,\hspace{0.15cm} \{{\rm Res}_{\cal M}\mid {\cal M}\in {\sf M}_{\cal A}^{\rm NDB}\}$. \hfill{\scriptsize (Theo.\,\ref{BinaryOracle})}

\vspace{0.2cm}}

Digitally non-deterministic $\nu$-computability

\vspace{0.1cm}

{\small
\quad $(d)$ or $\{c_1\}\in {\rm DEC}_{\cal A}$ 

\vspace{0.1cm}

\hspace{0.8cm} $\textcolor{blue}{\Rightarrow}$ \quad$\{{\rm Res}_{\cal M}\mid {\cal M}\in {\sf M}_{\cal A}^{\rm NDB}\} \hspace{0.3cm} \,\,\textcolor{blue}{\subseteq}\,\,\hspace{0.15cm} \{{\rm Res}_{\cal M}\mid {\cal M}\in {\sf M}_{\cal A}^\nu(Q_0)\}$. \hfill{\scriptsize (Ths. \ref{IdentBinary}, \ref{TheoRecBina})}

\vspace{0.2cm}}
}}
\end{overview}

\begin{overview}[Semi-decidability by several types of machines]

\hfill

\nopagebreak 

\noindent \fbox{\parbox{11.8cm}{

{\small \hfill Let ${\cal A}\in {\sf Struc}$.

\hspace{2cm}${\rm SDEC}_{\cal A} \,\,\textcolor{blue}{=}\,\,{\rm SDEC}_{\cal A}^{\rm NDB}$. \hfill{\scriptsize (Theo.\,\ref{BinaryBranchingDeterm})}}

\vspace{0.1cm}
}}

\noindent \fbox{\parbox{11.8cm}{
\hfill {\small Let ${\cal A}\in {\sf Struc}_{c_1,c_2}$ and $Q_0 \,\,\textcolor{blue}{=}\,\,\{c_1,c_2\}^2$. 

\quad $(e)$ $\,\,\textcolor{blue}{\Rightarrow}$ \quad${\rm SDEC}_{\cal A} \,\,\textcolor{blue}{\subseteq}\,\, ({\rm SDEC}_{\cal A}^\nu)^{Q_0} \,\,\textcolor{blue}{=}\,\,{\rm SDEC}_{\cal A}^{\rm DND} \,\,\textcolor{blue}{\subseteq}\,\,{\rm SDEC}_{\cal A}^{\rm ND}.$
\hfill{\scriptsize (Coro.\,\ref{CoroDND_ND})} 

\vspace{0.3cm}

\quad $(d)$ or $\{c_1\}\in {\rm DEC}_{\cal A}$ 

\vspace{0.1cm}

\hspace{0.8cm} $\textcolor{blue}{\Rightarrow}$ \quad ${\rm SDEC}_{\cal A} \,\,\textcolor{blue}{=}\,\,({\rm SDEC}_{\cal A}^\nu)^{Q_0} \,\,\textcolor{blue}{=}\,\,{\rm SDEC}_{\cal A}^{\rm NDB} \,\,\textcolor{blue}{=}\,\,{\rm SDEC}_{\cal A}^{\rm DND} \,\,\textcolor{blue}{\subseteq}\,\,{\rm SDEC}_{\cal A}^{\rm ND}.$}

\hfill{\scriptsize (Coro.\,\ref{CoroDND_ND}, Coro.\,\ref{CorIdenRecBina})}

\vspace{0.1cm}}}
\end{overview}

\section*{Outlook}\label{Outlook: Computation paths and complexity}
\addcontentsline{toc}{section}{\bf Outlook and Acknowledgment}
\markboth{OUTLOOK AND ACKNOWLEDGMENT}{OUTLOOK AND ACKNOWLEDGMENT}

In the next part, we will consider various types of universal machines and again use properties resulting from Consequences 2.2 in Part I which say that it is also possible to simulate machines over ${\cal A}_{\mathbb{N}}$ by executing index instructions of types (5) to (7). Later, we will also discuss the complexity of machines, functions, and decision problems for first-order structures with identity or without identity and we will consider machines using deterministic Moschovakis operators. 

The description of program paths of $\sigma$-programs by $\sigma$-terms for inputs of a fixed length and the evaluation of the corresponding systems  of {\sf conditions} expressed by {\sf atomic $\sigma$-formulas} and {\sf literals} are important for many proofs of statements on the undecidability of certain decision problems, statements on the uncomputability of some functions, and statements on the complexity of problems. An introduction to program and computations paths will be covered in a special paper that will be made available soon.

\section*{Acknowledgment}

I would like to thank the participants of my lectures on Theory of Abstract Computation and, in particular, Patrick Steinbrunner and Sebastian Bierba\ss{} for useful questions and the discussions. My thanks go also to all participants of meetings in Greifswald and Kloster and Pedro F. Valencia Vizcaíno. In particular, I would like to thank my co-authors Arno Pauly and Florian Steinberg for the discussions on operators related to several models of computation and Vasco Brattka, Philipp Schlicht, and Rupert H\"olzl for many interesting discussions.

{\small

\noindent For this article, we also used translators such as those of Google and DeepL and, for questions about English grammar, we used Microsoft Copilot.}

\newpage

\section*{Flowcharts and a solution}

\markboth{FLOWCHARTS AND A SOLUTION}{FLOWCHARTS AND A SOLUTION}

{\bf A flowchart for the program ${\cal P}_{{\cal M}_1}$ of the ${\cal A}_1$-machine ${\cal M}_1$}

\noindent{\unitlength0.9cm
\begin{picture}(13,14.4) \thicklines

{\sf
\put(4,14.18){\vector(0,-1){0.6}}
\put(4,13.28){\makebox(0,0){$^{\ell_1:}$\fbox{$I_2:=I_2+1$}}} 
\put(4,12.98){\vector(0,-1){0.6}}

 
\put(4,11.69){
\put(-2.9,0){\line(1,0){1.5}} \put(-1.8,0.15){\makebox(0,0){\it no}} 
\put(-0.7,0.6){\makebox(0,0){ $^{\ell_1+1 :}$ }}
\put(-2,-0.25){\makebox(0,0){\scriptsize goto $\ell_1+1$ }}
\put(-2.9,1.2){\line(1,0){2.4}} \put(-0.52,1.17){\vector(1,-1){0.5}}
\put(1.4,0){\line(2,0){1.6}} \put(1.8,0.1){\makebox(0,0){\it yes}}
\put(2.45,-0.25){\makebox(0,0){\scriptsize goto $\ell_2$ }}
\put(-1.4,0){\line(2,1){1.4}}
\put(-1.4,0){\line(2,-1){1.4}}
\put(0,-0.7){\line(2,1){1.4}}
\put(0,0.7){\line(2,-1){1.4}}
\makebox(0,0){$I_1=I_2${\it ?}}}


\put(1.1,11.68){\line(0,1){1.2}}
\put(7,11.68){\vector(0,-1){0.6}}
\put(6.5,10.75){\makebox(0,0){\hspace{0.33cm}$^{\ell_2 :}$\fbox{$Z_5:= c_3^0$}}}
\put(7,10.4){\vector(0,-1){0.6}}
\put(6.5,9.48){\makebox(0,0){\textcolor{red}{\hspace{0.33cm}$^{\ell_3 :}$}\fbox{\textcolor{red}{$Z_3:=f_1^1(Z_1)$}}}}
\put(7,9.13){\vector(0,-1){0.6}}
\put(6.5,8.18){\makebox(0,0){\textcolor{red}{$^{\ell_3+1 :}$}\fbox{\textcolor{red}{ $ Z_4:=f_1^1(Z_2)$}}}}
\put(7,7.85){\vector(0,-1){0.6}}
\put(6.5,6.93){\makebox(0,0){\textcolor{red}{\hspace{0.33cm}$^{\ell_4 :}$}\fbox{\textcolor{red}{$ Z_1:=f_2^1(Z_1)$}}}}
\put(7,6.6){\vector(0,-1){0.6}}
\put(6.5,5.65){\makebox(0,0){\textcolor{red}{$^{\ell_4+1 :}$}\fbox{\textcolor{red}{$ Z_2:=f_2^1(Z_2) $}}}}
\put(7,5.3){\vector(0,-1){0.6}}


\put(7,4){
\put(-1.4,0){\line(-2,0){1.53}} \put(-1.8,0.15){\makebox(0,0){\it no}} \put(-2.9,0){\line(0,1){1.2}} 
\put(-2.9,1.2){\line(1,0){2.4}} \put(-0.52,1.17){\vector(1,-1){0.5}}
\put(-0.55,0.6){\makebox(0,0){ $^{\ell_5 :}$}}
\put(-2,-0.25){\makebox(0,0){\scriptsize goto $\textcolor{blue}{\ell_5}$ }}
\put(1.4,0){\line(2,0){2}} \put(1.8,0.1){\makebox(0,0){\it yes}}\put(3.4,0){\vector(0,-1){1.3}} 
\put(2.85,-0.25){\makebox(0,0){\scriptsize goto $\textcolor{blue}{\ell_6}$ }}
\put(-1.4,0){\line(2,1){1.4}}
\put(-1.4,0){\line(2,-1){1.4}}
\put(0,-0.7){\line(2,1){1.4}}
\put(0,0.7){\line(2,-1){1.4}}
\makebox(0,0){$r_1^{2}(Z_3,Z_4)${\it ?}}}


\put(10.4,2){
\put(-1.4,0){\line(-2,0){4.9}} \put(-1.8,0.15){\makebox(0,0){\it yes}}
\put(-6.3,0){\vector(0,-1){0.7}}
\put(-0.55,0.6){\makebox(0,0){ $^{\ell_6 :}$ }}
\put(-2.2,-0.2){\makebox(0,0){\scriptsize goto $\textcolor{blue}{\ell_7}$ } }
\put(1.4,0){\line(2,0){1.4}} \put(1.8,0.15){\makebox(0,0){\it no}}
\put(2.3,-0.2){\makebox(0,0){\scriptsize goto $\textcolor{blue}{\ell_3}$ }} 
\put(2.8,0){\line(0,1){7.5}} \put(2.8,7.5){\vector(-1,0){5.1}}
\put(-1.4,0){\line(2,1){1.4}}
\put(-1.4,0){\line(2,-1){1.4}}
\put(0,-0.7){\line(2,1){1.4}}
\put(0,0.7){\line(2,-1){1.4}}
\makebox(0,0){$r_1^{2}(Z_4,Z_5)${\it ?}}}


\put(3.8,1){\makebox(0,0){ $^{\ell_7 :\,}$\fbox{\sf stop} }}
}
\end{picture}}

\noindent{\sf 
\begin{tabular}{lll}
${\cal P}_{{\cal M}_1}$: &$\ell_1: \, \, I_2:=I_2+1;$ \,\, $\ell_1+1: $\, if $I_1=I_2$ then goto $\ell_2$ else goto $\ell_1+1;$\\
& $\ell_2: \, \, Z_5:= c_3^0;$ \\ 
&\textcolor{red}{$\ell_3 : \,\,Z_3:=f_1^1(Z_1);$ \,\, $ Z_4:=f_1^1(Z_2); $} \\
&\textcolor{red}{$\ell_4 : \,\,Z_1:=f_2^1(Z_1);$ \,\, $ Z_2:=f_2^1(Z_2); $} \\
&$\ell_5 : \,\,$ if $r_1^{2}(Z_3,Z_4)$ then goto \textcolor{blue}{$\ell_6$} else goto $\textcolor{blue}{\ell_5};$ \\
&$\ell_6 : \,\,$ if $r_1^{2}(Z_4,Z_5)$ then goto $\textcolor{blue}{\ell_7}$ else goto $\textcolor{blue}{\ell_3};$ \\
&$\ell_7: \, \,$ stop. \\
\end{tabular}}

\newpage
{\bf A flowchart for the program ${\cal P}_{{\cal M}_2}$ of the ${\cal A}_1$-machine ${\cal M}_2$}

\noindent{\unitlength0.9cm
\begin{picture}(13,14.4) \thicklines

{\sf
\put(4,14.18){\vector(0,-1){0.6}}
\put(4,13.28){\makebox(0,0){$^{\ell_1:}$\fbox{$I_2:=I_2+1$}}} 
\put(4,12.98){\vector(0,-1){0.6}}

 
\put(4,11.69){
\put(-2.9,0){\line(1,0){1.5}} \put(-1.8,0.15){\makebox(0,0){\it no}} 
\put(-0.7,0.6){\makebox(0,0){ $^{\ell_1+1 :}$ }}
\put(-2.2,-0.25){\makebox(0,0){\scriptsize goto $\ell_7$ }}
\put(1.4,0){\line(2,0){1.6}} \put(1.8,0.1){\makebox(0,0){\it yes}}
\put(2.5,-0.25){\makebox(0,0){\scriptsize goto $\ell_2$ }}
\put(-1.4,0){\line(2,1){1.4}}
\put(-1.4,0){\line(2,-1){1.4}}
\put(0,-0.7){\line(2,1){1.4}}
\put(0,0.7){\line(2,-1){1.4}}
\makebox(0,0){$I_1=I_2${\it ?}}}


\put(1.1,11.68){\line(0,-1){7.7}} \put(1.1,3.98){\line(1,0){2.7}} \put(3.8,3.98){\vector(0,-1){2.7}}
\put(7,11.68){\vector(0,-1){0.6}}
\put(6.5,10.75){\makebox(0,0){\hspace{0.33cm}$^{\ell_2 :}$\fbox{$Z_5:= c_3^0$}}}
\put(7,10.4){\vector(0,-1){0.6}}
\put(6.5,9.48){\makebox(0,0){\textcolor{red}{\hspace{0.33cm}$^{\ell_3 :}$}\fbox{\textcolor{red}{$Z_3:=f_1^1(Z_1)$}}}}
\put(7,9.13){\vector(0,-1){0.6}}
\put(6.5,8.18){\makebox(0,0){\textcolor{red}{$^{\ell_3+1 :}$}\fbox{\textcolor{red}{ $ Z_4:=f_1^1(Z_2)$}}}}
\put(7,7.85){\vector(0,-1){0.6}}
\put(6.5,6.93){\makebox(0,0){\textcolor{red}{\hspace{0.33cm}$^{\ell_4 :}$}\fbox{\textcolor{red}{$ Z_1:=f_2^1(Z_1)$}}}}
\put(7,6.6){\vector(0,-1){0.6}}
\put(6.5,5.65){\makebox(0,0){\textcolor{red}{$^{\ell_4+1 :}$}\fbox{\textcolor{red}{$ Z_2:=f_2^1(Z_2) $}}}}
\put(7,5.3){\vector(0,-1){0.6}}


\put(7,4){
\put(-1.4,0){\line(-2,0){1.53}} \put(-2.9,0){\vector(0,-1){2.7}} \put(-1.8,0.15){\makebox(0,0){\it no}} 
\put(-0.55,0.6){\makebox(0,0){ $^{\ell_5 :}$ }}
\put(-2,-0.25){\makebox(0,0){\scriptsize goto $\textcolor{blue}{\ell_7}$ }}
\put(1.4,0){\line(2,0){2}} \put(1.8,0.1){\makebox(0,0){\it yes}}\put(3.4,0){\vector(0,-1){1.3}} 
\put(2.85,-0.25){\makebox(0,0){\scriptsize goto $\textcolor{blue}{\ell_6}$ }}
\put(-1.4,0){\line(2,1){1.4}}
\put(-1.4,0){\line(2,-1){1.4}}
\put(0,-0.7){\line(2,1){1.4}}
\put(0,0.7){\line(2,-1){1.4}}
\makebox(0,0){$r_1^{2}(Z_3,Z_4)${\it ?}}}


\put(10.4,2){
\put(-1.4,0){\line(-2,0){1.5}} \put(-1.8,0.15){\makebox(0,0){\it yes}} \put(-2.9,0){\line(0,1){1.2}} 
\put(-2.9,1.2){\line(1,0){2.4}} \put(-0.52,1.17){\vector(1,-1){0.5}}
\put(-0.55,0.6){\makebox(0,0){ $^{\ell_6 :}$}}
\put(-2.2,-0.2){\makebox(0,0){\scriptsize goto $\textcolor{blue}{\ell_6}$}}
\put(1.4,0){\line(2,0){1.4}} \put(1.8,0.15){\makebox(0,0){\it no}}
\put(2.3,-0.2){\makebox(0,0){\scriptsize goto $\textcolor{blue}{\ell_3}$ }} 
\put(2.8,0){\line(0,1){7.5}} \put(2.8,7.5){\vector(-1,0){5.1}}
\put(-1.4,0){\line(2,1){1.4}}
\put(-1.4,0){\line(2,-1){1.4}}
\put(0,-0.7){\line(2,1){1.4}}
\put(0,0.7){\line(2,-1){1.4}}
\makebox(0,0){$r_1^{2}(Z_4,Z_5)${\it ?}}}

\put(3.8,1){\makebox(0,0){ $^{\ell_7 :\,}$\fbox{\sf stop} }}
}
\end{picture}}

{\sf 
\begin{tabular}{lll}
${\cal P}_{{\cal M}_2}$:&\quad $\ell_1: \, \, I_2:=I_2+1;$ \,\, if $I_1=I_2$ then goto $\ell_2$ else goto $\ell_7;$\\
&\quad $\ell_2: \, \, Z_5:= c_3^0;$ \\ 
&\textcolor{red}{\quad $\ell_3 : \,\,Z_3:=f_1^1(Z_1);$ \,\, $ Z_4:=f_1^1(Z_2); $} \\
&\textcolor{red}{\quad $\ell_4 : \,\,Z_1:=f_2^1(Z_1);$ \,\, $ Z_2:=f_2^1(Z_2); $} \\
&\quad $\ell_5 : \,\,$ if $r_1^{2}(Z_3,Z_4)$ then goto $\textcolor{blue}{\ell_6}$ else goto $\textcolor{blue}{\ell_7};$ \\
&\quad $\ell_6 : \,\,$ if $r_1^{2}(Z_4,Z_5)$ then goto $\textcolor{blue}{\ell_6}$ else goto $\textcolor{blue}{\ell_3};$ \\
&\quad$\ell_7: \, \,$ stop. \\
\end{tabular}}

\newpage

\subsection*{A solution referring to Example \ref{semi_co_semi}} The complement $U_{{\cal A}_4}^\infty\setminus {\rm id}_{U_{{\cal A}_4}}$ is semi-decidable by a BSS RAM using the following program.

\vspace{0.3cm}

{\sf 
\begin{tabular}{lll}
&\quad $\ell_1: \, \, I_2:=I_2+1;$ \,\, if $I_1=I_2$ then goto $\ell_{\textcolor{magenta}{3}}$ else goto $\ell_7;$\\
&\textcolor{red}{\quad $\ell_3 : \,\,Z_3:=f_1^1(Z_1);$ \,\, $ Z_4:=f_1^1(Z_2); $} \\
&\textcolor{red}{\quad $\ell_4 : \,\,Z_1:=f_2^1(Z_1);$ \,\, $ Z_2:=f_2^1(Z_2); $} \\
&\quad $\ell_5 : \,\,$ if $r_1^{2}(Z_3,Z_4)$ then goto $\textcolor{blue}{\ell_{\textcolor{magenta}{3}}}$ else goto $\textcolor{blue}{\ell_7};$ \\
&\quad$\ell_7: \, \,$ stop. \\
\end{tabular}}

\noindent{\unitlength0.9cm
\begin{picture}(13,14.4) \thicklines
{\sf
\put(4,14.18){\vector(0,-1){0.6}}
\put(4,13.28){\makebox(0,0){$^{\ell_1:}$\fbox{$I_2:=I_2+1$}}} 
\put(4,12.98){\vector(0,-1){0.6}}

 
\put(4,11.69){
\put(-2.9,0){\line(1,0){1.5}} \put(-1.8,0.15){\makebox(0,0){\it no}} 
\put(-0.7,0.6){\makebox(0,0){ $^{\ell_1+1 :}$ }}
\put(-2.2,-0.25){\makebox(0,0){\scriptsize goto $\ell_7$ }}
\put(1.4,0){\line(2,0){1.6}} \put(1.8,0.1){\makebox(0,0){\it yes}}
\put(2.5,-0.25){\makebox(0,0){\scriptsize goto $\ell_{\textcolor{magenta}{3}}$ }}
\put(-1.4,0){\line(2,1){1.4}}
\put(-1.4,0){\line(2,-1){1.4}}
\put(0,-0.7){\line(2,1){1.4}}
\put(0,0.7){\line(2,-1){1.4}}
\makebox(0,0){$I_1=I_2${\it ?}}}


\put(1.1,11.68){\line(0,-1){7.7}} \put(1.1,3.98){\line(1,0){2.7}} \put(3.8,3.98){\vector(0,-1){2.7}}
\put(7,11.68){\vector(0,-1){1.85}}
\put(6.5,9.48){\makebox(0,0){\textcolor{red}{\hspace{0.33cm}$^{\ell_3 :}$}\fbox{\textcolor{red}{$Z_3:=f_1^1(Z_1)$}}}}
\put(7,9.13){\vector(0,-1){0.6}}
\put(6.5,8.18){\makebox(0,0){\textcolor{red}{$^{\ell_3+1 :}$}\fbox{\textcolor{red}{ $ Z_4:=f_1^1(Z_2)$}}}}
\put(7,7.85){\vector(0,-1){0.6}}
\put(6.5,6.93){\makebox(0,0){\textcolor{red}{\hspace{0.33cm}$^{\ell_4 :}$}\fbox{\textcolor{red}{$ Z_1:=f_2^1(Z_1)$}}}}
\put(7,6.6){\vector(0,-1){0.6}}
\put(6.5,5.65){\makebox(0,0){\textcolor{red}{$^{\ell_4+1 :}$}\fbox{\textcolor{red}{$ Z_2:=f_2^1(Z_2) $}}}}
\put(7,5.3){\vector(0,-1){0.6}}


\put(7,4){
\put(-1.4,0){\line(-2,0){1.53}} \put(-2.9,0){\vector(0,-1){2.7}} \put(-1.8,0.15){\makebox(0,0){\it no}} 
\put(-0.55,0.6){\makebox(0,0){ $^{\ell_5 :}$ }}
\put(-2,-0.25){\makebox(0,0){\scriptsize goto $\textcolor{blue}{\ell_7}$ }}
\put(1.4,0){\line(1,0){2.8}} \put(1.8,0.1){\makebox(0,0){\it yes}}
\put(2.85,-0.25){\makebox(0,0){\scriptsize goto $\textcolor{blue}{\ell_{\textcolor{magenta}{3}}}$ }}
\put(-1.4,0){\line(2,1){1.4}}
\put(-1.4,0){\line(2,-1){1.4}}
\put(0,-0.7){\line(2,1){1.4}}
\put(0,0.7){\line(2,-1){1.4}}
\makebox(0,0){$r_1^{2}(Z_3,Z_4)${\it ?}}}


\put(10.4,2){
\put(0.8,2){\line(0,1){5.5}} \put(0.8,7.5){\vector(-1,0){3.1}}
}


\put(3.8,1){\makebox(0,0){ $^{\ell_7 :\,}$\fbox{\sf stop} }}
}
\end{picture}}

\newpage 

\subsection*{A 3-tape BSS RAM ${\cal N}_{\cal M}^{(3)}$ for computing $\chi_{{\rm id}_{\cal A}}$}

\begin{overview}[The program ${\cal P}_{{\cal N}_{\cal M}^{(3)}}$ incl. \textcolor{blue}{${\cal P}_1$} and \textcolor{blue}{${\cal P}_2$}]\label{CharactchiId}

\hfill

\nopagebreak 

\noindent \fbox{\parbox{11.8cm}{\begin{tabular}{rlrl}

\textcolor{blue}{$\tilde 1\,':$}&$\!\!\!${\sf if $I_{3,2}=i$ then goto \textcolor{blue}{$i''$}$;$}\\
\textcolor{blue}{$\tilde 1\,'':$}&$\!\!\!${\sf if $I_{3,1}=i$ then goto \textcolor{blue}{$i'$}$;$}\\
\textcolor{blue}{$1\,' :$}&$\!\!\!${\sf inst}$'_1;\, \,\ldots;\,\hspace{0.34cm}$\textcolor{blue}{$( \ell _{{\cal P}_{{\cal M}}}-1)\,'$}$\, :$ {\sf inst}$'_{ \ell _{{\cal P}_{{\cal M}}}-1};$ \\
\textcolor{blue}{$1\,'':$}&$\!\!\!${\sf inst}$''_1; \,\,\ldots;\,\,\,\textcolor{blue}{( \ell _{{\cal P}_{{\cal M}}}-1)\,''}: $ {\sf inst}$''_{ \ell _{{\cal P}_{{\cal M}}}-1};$ \, $\ell\,'' _{{\cal P}_{{\cal M}}}:\,I_{1,1}=1;  \, $ {\sf goto $\tilde  \ell;$}\\\\
$\tilde \ell:$ &$\!\!\!${\sf if $I_{ 3,1}=I_{3,2}$ then goto $\tilde \ell^*$ else goto $\tilde \ell^{**};$}\vspace{0.2cm}\\
$\tilde \ell^*:$      &$\!\!\!$$ Z_1:= c_1^0; \, \,\, $ {\sf goto $\ell _{{\cal P}_{_{{\cal N}_{\cal M}^{(3)}}}};$}\\
$\tilde \ell^{**}: $&$\!\!\!$$Z_1:= c_2^0; $\vspace{0.2cm}\\
$\ell _{{\cal P}_{_{{\cal N}_{\cal M}^{(3)}}}}:$&$\!\!\!${\sf stop}.
\end{tabular}

}}
\end{overview}

For the program segments \textcolor{blue}{${\cal P}_1$} and \textcolor{blue}{${\cal P}_2$}, see Overviews \ref{CharactchiP2} and \ref{CharactchiP3}. 

Note, that we here use  ${\cal M}'=_{\rm df}{\cal M}$ and  ${\cal M}''=_{\rm df}{\cal M}$.

\vfill

\hfill {\scriptsize\sf C $\cdot $ H $\cdot$ T}

\addcontentsline{toc}{section}{\bf References}

\begin{thebibliography}{57}
\bibitem{AHU74} Alfred V. Aho, John E. Hopcroft, and Jeffrey D. Ullman {\sf The Design and Analysis of Computer Algorithms}, Addison-Wesley 1974.

\bibitem{A75} G\"unter Asser {\sf Einf\"uhrung in die mathematische Logik. Teil 2, Pr\"adikatenkalk\"ul der ersten Stufe}, BSB B.\,G.\,Teubner Verlagsgesellschaft, Leipzig 1975.

\bibitem{BSS89} Lenore Blum, Michael Shub, and Steve Smale {\sl On a theory of computation and complexity over the real numbers; {NP}-completeness, recursive functions and universal machines}, {\sf Bulletin of the American Mathematical Society 21} 1989, pp. 1--46. 

\bibitem{BCSS98} Lenore Blum, Felipe Cucker, Michael Shub, and Steve Smale {\sf Complexity and Real Computation}, Springer 1998. 

\bibitem{B92} Egon B\"orger {\sf Berechenbarkeit, Komplexit\"at, Logik}, Vieweg 1992. (In English: {\sf Computability, Complexity, Logic}, Elsevier 1989)

\bibitem{BP18} Olivier Bournez, and Amaury Pouly {\sl A survey on analog models of computation}, {\sf arXiv:1805.05729}, 2018.

\bibitem{BH98} Vasco Brattka, and Peter Hertling {\sl Feasible real random access machines}, {\sf Journal of Complexity 14}, 1998, pp. 490--526.

\bibitem{CUCKER92} Felipe Cucker {\sl The arithmetical hierarchy over the reals}, {\sf Journal of Logic and Computation 2} (3), 1992, pp. 375--395.

\bibitem{GASS96} Christine Ga\ss ner {\sl An NP-complete problem for linear real machines}, Workshop on Computability, Complexity and Logic (WCCL'96) Zinnowitz, 1996, pp. 37--38. 

\bibitem{GASS97} Christine Ga\ss ner {\sl On NP-completeness for linear machines}, {\sf Journal of Complexity} 13, 1997, pp. 259--271. 

\bibitem{GASS01} Christine Ga\ss ner {\sl The P-DNP problem for infinite Abelian groups}, {\sf Journal of Complexity} 17, 2001, pp. 574–583. 

\bibitem{GASS07} Christine Ga\ss ner {\sl P = NP for expansions derived from some oracles}, S. Barry Cooper, Thomas F. Kent, Benedikt Löwe, and Andrea Sorbi (eds.): Computation and Logic in the Real World, CiE 2007, in the series Quaderni del Dipartimento di Scienze Matematiche e Informatiche "Roberto Magari" of the University of Siena. Technical report no. 487, June 2007, pp. 161--169. 

{https://www.researchgate.net/publication/385686638} 

\bibitem{GASS08A} Christine Ga\ss ner {\sl A hierarchy below the halting problem for additive machines}, {\sf Theory of Computing Systems 43}, 2008, pp. 464--470.

\bibitem{GASS08B} Christine Ga\ss ner {\sl On the power of relativized versions of P, DNP, and NP for groups}, {\sf Preprint-Reihe Mathematik, Greifswald}, Preprint 8/2008. 

\bibitem{GASS08C} Christine Ga\ss ner {\sl On relativizations of the P =? NP question for several structures},
 {\sf Electronic Notes in Theoretical Computer Sciene 221}, 2008, pp. 71--83.

\bibitem{GASS08D} Christine Ga\ss ner {\sl Computation over groups}, Arnold Beckmann, Costas Dimitracopoulos, and Benedikt L\"owe (eds.): {\sf Logic and Theory of Algorithms}, 2008, pp. 147--156.

 {https://www.researchgate.net/publication/385686634}

\bibitem{GASS09A} Christine Ga\ss ner {\sl Oracles and relativizations of the P =? NP question for several structures}, {\sf Journal of Universal Computer Science 15} (6), 2009, pp. 1186--1205. 

\bibitem{GASS09B} Christine Ga\ss ner {\sl Relativizations of the P =? DNP question for the BSS model}, Andrej Bauer, Peter Hertling, Ker-I Ko (eds.): The Sixth International Conference on Computability and Complexity in Analysis, Proceedings 2009, pp. 141--148. DOI: 10.4230/OASIcs.CCA.2009.2266

\bibitem{GASS10} Christine Ga\ss ner {\sl The separation of relativized versions of P and DNP for the ring of the reals}, {\sf Journal of Universal Computer Science 16} (18), 2010, pp. 2563--2568.

\bibitem{GASS11} Christine Ga\ss ner {\sl Gel\"oste und offene P-NP-Probleme \"uber verschiedenen Strukturen}, Habilitations\-schrift, Greifswald 2011.

\bibitem{GASS13} Christine Ga\ss ner {\sl Strong Turing degrees for additive BSS RAM's}, {\sf Logical Methods in Computer Science 9} (4:25), 2013, pp. 1--18. {\sf arXiv:1312.3927}.

\bibitem{GASS15} Christine Ga\ss ner {\sl Hierarchies of decision problems over algebraic structures defined by quantifiers}, CCC 2015. DOI: 10.13140/RG.2.2.35627.14882

\bibitem{GaVV15} Christine Ga\ss ner, and Pedro F. Valencia Vizcaíno {\sl Operators for BSS RAM's}, Martin Ziegler and Akitoshi Kawamura (eds.): The Twelfth International Conference on Computability and Complexity in Analysis, Proceedings 2015, pp.\,\,24--26.
 {https://www.kurims.kyoto-u.ac.jp/$\sim$kawamura/cca2015-local/cca2015proceedings.pdf}

\bibitem{GASS16} Christine Ga\ss ner {\sl BSS RAM's with $\nu$-oracle and the axiom of choice}, Colloquium Logicum 2016.

{https://www.math.uni-hamburg.de/spag/ml/CL2016/contributedpapers.html}

\bibitem{GASS17} Christine Ga\ss ner {\sl Computation over algebraic structures and a classification of undecidable problems}, {\sf Mathematical Structures in Computer Science 27} (8), 2017, pp. 1386--1413.

\bibitem{GaPaSt17} Christine Ga\ss ner, Arno Pauly, and Florian Steinberg {\sl Computing measures as a primitive operation}, Daniel Graça and Martin Ziegler (eds.): The Fourteenth International Conference on Computability and Complexity in Analysis, Proceedings 2017, pp. 26--27.

\bibitem{GaPaSt18} Christine Ga\ss ner, Arno Pauly, and Florian Steinberg {\sl Comparing integrability and measurability in different models of computation}, Florin Manea, Russell G. Miller, Dirk Nowotka (eds): 14th Conference on Computability in Europe, Local proceedings 2018. HAL Id: hal-02019195.

\bibitem{GaPaSt21} Christine Ga\ss ner, Arno Pauly, and Florian Steinberg {\sl Computing measure as a primitive operation in real number computation}, {\sf Leibniz International Proceedings in Informatics Vol. 183}. DOI: 10.4230/LIPIcs.CSL.2021.22

\bibitem{GASS20} Christine Ga\ss ner {\sl An introduction to a model of abstract computation: The BSS-RAM model}, Adrian Rezus (ed.): Contemporary Logic and Computing, College Publications [{\sf Landscapes in Logic 1}], London 2020, pp. 574--603. 

\bibitem{GASS19Preprint} ---,--- Preprint 2019.

{https://www.researchgate.net/publication/343039844}

\bibitem{KleinEnz} W. Gellert, H. K\"ustner, M. Hellwich, and H. K\"astner (eds.): {\sf Kleinen Enzyklop\"adie Mathematik}, VEB Interdruck Leipzig 1968.

\bibitem{H98} Armin Hemmerling {\sl Computability of string functions over algebraic structures}, {\sf Mathematical Logic Quarterly 44}, 1998, pp. 1--44.

\bibitem{Hennie} Fred C. Hennie, and Richard E. Stearns {\sl Two-tape simulation of multitape Turing machines}, {\sf Journal of the Association for Computing Machinery 13} (4), 1966, pp. 533--546.

\bibitem{HU79} John E. Hopcroft, and Jeffrey D. Ullman {\sl Einf\"uhrung in die Automatentheorie, formale Sprachen und Komplexit\"atstheorie}, Addison-Wesley 1993. (In English: {\sf Introduction to Automata Theory, Languages, and Computation}, Addison-Wesley 1979.)

\bibitem{KLEENE} Stephen C. Kleene {\sf Introduction to Metamathematics}, Groningen, Wolters-Noordhoff et al. 2000.

\bibitem{WL83} Nikolai I. Kondakow {\sf W\"orterbuch der Logik} (The editors of the German edition: Erhard Albrecht and G\"unter Asser), Bibliographisches Institut, Leipzig 1983.

\bibitem{KOZEN} Dexter C. Kozen {\sf Theory of Computation}, Springer 2006.

\bibitem{MEERZIEGLER06} Klaus Meer, and Martin Ziegler {\sl Uncomputability below the real halting problem}, CiE 2006, {\sf Lecture Notes in Computer Science 3988}, 2006, pp. 368--377.

\bibitem{MOSCHO} Yiannis N. Moschovakis {\sl Abstract first order computability I}, {\sf Transactions of the American Mathematical Society 138}, 1969, pp. 427--464. 

\bibitem{P95} Bruno Poizat {\sf Les Petits Cailloux:$\!$ Une Approche Modèle-Théorique De L'Algorithmie}, Al\'eas 1995. 

\bibitem{PS85} Franco P. Preparata, and Michael I. Shamos {\sf Computational Geometry. An Introduction}, Springer 1985.

\bibitem{R67} Hartley Rogers {\sf Theory of Recursive Functions and Effective Computability}, McGraw-Hill 1967.

\bibitem{S67} Dana Scott {\sl Some definitional suggestions for automata theory}, {\sf Journal of Computer and System Sciences 1}, 1967, pp. 187--212.

\bibitem{ShS63} John C. Shepherdson, and Howard E. Sturgis {\sl Computability of recursive functions}, {\sf Journal of the Association for Computing Machinery 10}, 1963, pp. 217--255.

\bibitem{SOARE} Robert I. Soare {\sf Recursively Enumerable Sets and Degrees: A Study of Computable Functions and Computably Generated Sets}, Springer, 1987.

\bibitem{TURING37} Alan M. Turing {\sl On computable numbers, with an application to the Entscheidungsproblem},
{\sf Proceedings of the London Mathematical Society} 42 (1), 1937,
pp. 230--265.

\bibitem{WEIHRAUCH} Klaus Weihrauch {\sf Computable analysis. An Introduction}, Springer 2000.

\bibitem{Zerjatke} Thomas Zerjatke {\sf Möglichkeiten und Grenzen des DNA-Computings im Hinblick auf PSPACE}, Greifswald 2009.
{https://math-inf.uni-greifswald.de/storages/uni-greifswald/fakultaet/mnf/mathinf/gassner/ZerjatkeDiplomarbeit.pdf}

\vspace{0.5cm}



\bibitem{Prun02} Mihai Prunescu {\sl A model-theoretic proof for $P\not= NP$ over all infinite Abeian groups}, {\sf The Journal of Symbolic Logic}, 67 (1), 2002, pp. 235--238.


\vspace{0.5cm}

$\!\!\!\!\!\!\!\!\!\!\!\!\!$Part I, Part II\,a:

\bibitem{GASS25A} 
Christine Ga\ss ner {\sl Abstract computation over first-order structures. Part I: Deterministic and non-deterministic BSS RAMs}, Preprint 2025. {\sf arXiv:2502.17539}.

\bibitem{GASS25B}
Christine Ga\ss ner {\sl Abstract Computation over First-Order Structures. Part II\,a: Moschovakis' Operator and Other Non-Determinisms}, Preprint 2025. {\sf arXiv:2504.21192} 

\vspace{0.5cm}


$\!\!\!\!\!\!\!\!\!\!\!\!\!$Further books and papers:

\bibitem{CSS94}Felipe Cucker, Michael Shub, and Steve Smale {\sl Separation of complexity classes in Koiran's weak model}, {\sf Theoretical Computer Science} 133, 1994, pp 3--14. 

\bibitem{Mendelson} 
{\sc Mendelson,\,\, E.},\,\,\,{\sf Introduction to Mathematical Logic}, The University series in undergraduate mathematics. D.van Nostrand Campany, inc. Princeton-New Jersey-Toronto-New York- London, 1964, reprinted 1965. Volume 100, North Holland 1980.

$\!\!\!\!\!\!\!\!\!\!\!\!\!$Further textbooks:

\bibitem{BDGI} José Luis Balc\'azar, Josep D\'iaz, and Joaquim Gabarr\'o {\sf Structural Complexity I}, Springer 1988. 

\bibitem{BlumN} Norbert Blum {\sf Theoretische Informatik: eine anwenderorientierte Einf\"uhrung}, Oldenbourg 2001.

\bibitem{Hoff09} Dirk W. Hoffmann {\sf Theoretische Informatik}, Hanser 2009.

\bibitem{Papad} Christos C. Papadimitriou {\sf Computational Complexity}, Addison-Wesley 1994.

\bibitem{Reischuk99} K. Rüdiger Reischuk {\sf Komplexit\"atstheorie - Band I: Grundlagen: Maschinenmodelle, Zeit- und Platzkomplexität, Nichtdeterminismus}, Teubner 1999.

\bibitem{Schoening08} Uwe Sch\"oning {\sf Theoretische Informatik - kurz gefasst}, Spektrum 2008.

\bibitem{TuschikWolter} Hans-Peter Tuschik, and Helmut Wolter {\sf Mathematische Logik - kurz gefasst}, Spektrum 2002.

 \end{thebibliography}
\end{document}